\documentclass{amsart}

\usepackage{amsrefs} 

\usepackage{epsf,color}
\usepackage{amsmath}
\usepackage{amscd}
\usepackage{amssymb}
\usepackage{enumerate}

\usepackage{graphicx}
\usepackage[all]{xy}
\usepackage{mathrsfs}

\usepackage{bbold}
\usepackage{epsfig}

\numberwithin{equation}{section}

\newcommand{\into}{\hookrightarrow}

\newcommand{\bC}{{\mathbb C}}

\newcommand{\bP}{{\mathbb P}}

\newcommand{\bR}{{\mathbb R}}

\newcommand{\bT}{{\mathbb T}}

\newcommand{\bZ}{{\mathbb Z}}

\newcommand{\cp}{{\mathfrak p}}
\newcommand{\ft}{{\mathfrak t}}
\newcommand{\cA}{{\mathcal A}}
\newcommand{\cB}{{\mathcal B}}

\newcommand{\cE}{{\mathcal E}}
\newcommand{\cF}{{\mathcal F}}

\newcommand{\cH}{{\mathcal H}}

\newcommand{\cL}{{\mathcal L}}
\newcommand{\cM}{{\mathcal M}}

\newcommand{\cQ}{{\mathcal Q}}

\newcommand{\cS}{{\mathcal S}}
\newcommand{\cT}{{\mathcal T}}

\newcommand{\sE}{{\mathscr E}}
\newcommand{\sF}{{\mathscr F}}

\newcommand{\sH}{{\mathscr H}}

\newcommand{\sL}{{\mathscr L}}

\newcommand{\sN}{{\mathscr N}}

\newcommand{\sP}{{\mathscr P}}
\newcommand{\sQ}{{\mathscr Q}}
\newcommand{\sR}{{\mathscr R}}
\newcommand{\sS}{{\mathscr S}}
\newcommand{\sT}{{\mathscr T}}

\newcommand{\sV}{{\mathscr V}}

\newcommand{\bfJ}{{\mathbf J}}

\newcommand{\vf}[1][\!]{\vec{f}^{\, #1}}
\newcommand{\vp}[1][\!]{\vec{p}^{\, #1}}

\newcommand{\vC}[1][\!]{\vec{C}^{\, #1}}

\newcommand{\vb}[1][\!]{\vec{b}^{\, #1}}

\newcommand{\vt}[1][\!]{\vec{t}^{\, #1}}
\newcommand{\vs}[1][\!]{\vec{s}^{\, #1}}
\newcommand{\vE}[1][\!]{\vec{E}^{\, #1}}
\newcommand{\vnu}[1][\!]{\vec{\nu}^{\, #1}}

\newcommand{\Log}{\operatorname{Log}}

\newcommand{\Fuk}{\operatorname{Fuk}}

\newcommand{\Coh}{\operatorname{Coh}}
\newcommand{\grad}{\operatorname{grad}}

\newcommand{\CF}{\operatorname{CF}}
\newcommand{\HM}{\operatorname{HM}}
\newcommand{\CM}{\operatorname{CM}}

\newcommand{\Crit}{\operatorname{Crit}}
\newcommand{\Diff}{\operatorname{Diff}}
\newcommand{\SL}{\operatorname{SL}}
\newcommand{\Ob}{{\mathcal Ob}}
\newcommand{\eva}{{\mathrm{ev}}}
\newcommand{\Rel}{\operatorname{Simp}}
\newcommand{\coll}{\operatorname{coll}}
\newcommand{\inte}{{\operatorname{int}}}
\newcommand{\tr}{{\operatorname{tr}}}
\newcommand{\Stasheff}{\cT}
\newcommand{\Shrub}{\cS}
\newcommand{\codim}{\operatorname{codim}}
\newcommand{\virdim}{\operatorname{virdim}}

\newcommand{\Morse}{\operatorname{Morse}}
\newcommand{\Cell}{\operatorname{Cell}}
\newcommand{\Cech}{\operatorname{\check{C}ech}}
\newcommand{\Part}{\operatorname{Part}}

\def\co{\colon\thinspace}

\newcommand{\ctorus}[1]{(\mathbb{C}^{\star})^{#1}}

\newtheorem{thm}{Theorem}[section]
\newtheorem*{HMS1}{HMS Conjecture for CY}
\newtheorem*{HMS2}{HMS Conjecture for Fano}
\newtheorem*{SYZ}{SYZ Conjecture}
\newtheorem{cor}[thm]{Corollary}
\newtheorem{lem}[thm]{Lemma}
\newtheorem{prop}[thm]{Proposition}
\newtheorem{defin}[thm]{Definition}
\newtheorem{conj}{Conjecture}

\theoremstyle{remark}
\newtheorem{rem}[thm]{Remark}
\newtheorem{rems}[thm]{Remarks}
\newtheorem*{claim}{Claim}

\newcommand{\noproof}{\hfill \qedsymbol}

\newcommand{\comment}[1]{}

\title[Morse Homology and Mirror Symmetry]{Morse Homology, Tropical Geometry, and Homological Mirror Symmetry for Toric Varieties}

\author[M.~Abouzaid]{Mohammed Abouzaid}
\address{Massachusetts Institute of Technology \\
77 Massachusetts Ave, \\
Cambridge, MA, USA}

\begin{document}

\begin{abstract}
Given a smooth projective toric variety $X$, we construct an $A_\infty$ category of Lagrangians with boundary on a level set of the Landau-Ginzburg mirror of $X$.  We prove that this category is quasi-equivalent to the $DG$ category of line bundles on $X$.
\end{abstract}

\maketitle

\tableofcontents

\section{Introduction}
Mirror symmetry entered the mathematical realm when, in \cite{COGP}, physicists were able to enumerate curves on the quintic threefold by studying the complex moduli space of an algebraic variety which is most often referred to as the ``mirror quintic.'' The fact that such a computation was impossible using then available mathematical techniques contributed to the development of an entirely new field of mathematics (Gromov-Witten theory), which eventually resulted in Givental's proof in \cites{givental1,givental2} of a mirror conjecture for complete intersections in toric varieties. One important feature apparent in Givental's work is that, whereas the mirror of a Calabi-Yau variety (one whose canonical bundle is trivial) is another Calabi-Yau variety, the mirrors of Fano varieties (essentially, those with positive anti-canonical bundle) are Landau-Ginzburg models; i.e. affine varieties equipped with a map to $\bC$ called the superpotential.

Meanwhile, there have been attempts at explaining the mirror phenomenon.  Most relevant to us are categorical ideas of Kontsevich and geometric ideas of Strominger, Yau, Zaslow.  In his 1994 ICM address, Kontsevich sought to explain mirror symmetry for pairs of mirror projective Calabi-Yau manifolds $M$ and $\check{M}$ as a consequence of an equivalence of categories which he called homological mirror symmetry:
\begin{HMS1}[\cite{kontsevich}]
The derived category of coherent sheaves on $M$ is equivalent to the derived Fukaya category of $\check{M}$.
\end{HMS1}
A few words of explanations are in order.  Assuming that $M$ is smooth, the derived category $D^b \Coh(M)$ of coherent sheaves on $M$ is generated by (algebraic or, equivalently, holomorphic) vector bundles on $M$ and the morphisms in this category are essentially determined by the higher cohomology groups of such vector bundles.  In particular, $D^b \Coh(M)$ depends only on the complex structure of $M$, and not on any projective embedding.

On the other hand, the derived Fukaya category depends on the symplectic form (essentially the imaginary part of the K\"ahler form) induced by a projective embedding, and is independent of the complex structure.  In essence, the study of the Fukaya category is a problem in symplectic topology, and not in algebraic geometry.  Although current technology makes a satisfactory definition of the Fukaya category possible, setting up all the preliminaries would take us far beyond what would be appropriate for an introduction.  Let us simply say that the objects of such a category are Lagrangian submanifolds (i.e. manifolds on which the symplectic form vanishes) and that the definition of the spaces of morphisms and of the composition entails the study of holomorphic curves with respect to an almost complex structure which is compatible with the symplectic form.  Lest the reader object that the Fukaya category is supposed to depend only the symplectic structure, one should note that, while the choice of an almost complex structure is necessary to give a definition, the quasi-isomorphism class of the resulting category is independent of this choice.  This is essentially the same phenomenon which makes the Gromov-Witten invariants of an algebraic variety deformation invariants.

In his 1999 ENS lectures \cite{kont-ENS}, Kontsevich extended his homological mirror symmetry to Fano varieties.  Since this is the main subject of this paper, we shall return to it in more detail.  Before continuing with homological mirror symmetry, let us discuss the SYZ explanation of mirror symmetry.  As is clear from the title of their paper \cite{SYZ}, their main conjecture is:
\begin{SYZ}
If $M$ and $\check{M}$ are a mirror pair of Calabi-Yau threefolds, then there exists (singular) special Lagrangian torus fibrations $M \to B$ and $\check{M} \to B$ which are dual to each other.
\end{SYZ}
To understand what is meant by dual fibrations, it is best to consider an example which generalizes to describe the phenomenon whenever there are no singularities.  Consider any smooth manifold $B$.  It is well known that the cotangent bundle $T^{*}B$ is a symplectic manifold, whereas an infinitesimal neighbourhood of the zero section in the tangent bundle $TB$ admits a canonical complex structure.  In nice situations, this complex structure extends to the entire tangent space, and for simplicity let us assume that this is the case.  Upon choosing a metric on $B$, we induce a complex structure on $T^*B$ and a symplectic structure on $TB$.  In order to obtain torus fibrations, one should consider a manifold equipped with an integral affine structure, i.e. a reduction of the structure group of $B$ from $\Diff(\bR^{n})$ to the group of diffeomorphisms of $\bR^{n}$ whose differentials preserve the standard lattice.  Such a manifold comes equipped with lattices in its tangent and cotangent bundles.  The corresponding quotients are the desired (prototypical) dual torus fibrations.  For discussion of mirror symmetry for non-singular torus fibrations, see \cite{leung}.  

The only part of the SYZ conjecture relevant to this paper is the following important observation which is by now folkloric:
\begin{conj} \label{section=line}
The mirror of a Lagrangian section of $M \to B$ is a holomorphic line bundle on $\check{M}$.
\end{conj}
While the SYZ conjecture has only been stated for Calabi-Yau manifolds, we shall nonetheless use the intuition that it provides.  We note that toric varieties are, in a sense, equipped with Lagrangian torus fibrations, so the appearance of ideas from SYZ in a discussion of homological mirror symmetry for toric varieties should not be surprising.

\subsubsection*{Statement of results}
The main result concerns half of Kontsevich's homological mirror symmetry conjecture:
\begin{HMS2}[\cite{kont-ENS}]
If the anti-canonical bundle of $X$ is ample (i.e. $X$ is Fano), the bounded derived category of coherent sheaves on $X$, $D^{b}\Coh(X)$, is equivalent, as an $A_{\infty}$ category, to the Karoubi completion of the derived Fukaya category of the mirror of $X$.
\end{HMS2}
The conjecture does not include a description of this conjectural mirror, although a wide class of mirrors is described in \cite{HV}.  In fact, we will be working with smooth projective toric varieties even when their anti-canonical bundles are not ample.  However, as we shall see in the discussion following the statement of Theorem \ref{main}, we cannot expect an equivalence of categories without the Fano condition.

The fact that the homological mirror symmetry conjecture is a statement about $A_{\infty}$ categories means that we are seeking equivalences at the chain level.  The topologically inclined reader should have in mind secondary products on cohomology (for example the Massey products) which could not be detected by someone who simply thinks of cohomology as a functor satisfying the Eilenberg-Steenrod axioms; these are essentially chain-level phenomena.  The algebraically inclined reader should have in mind the fact that derived categories, while extremely useful, have been increasingly supplanted in recent years by DG categories, which remember more information about the original object which is often lost by deriving, see \cite{drinfeld}.  In an appropriate homotopy category, every (unital) $A_{\infty}$ category is equivalent to a DG category.

Let us briefly describe what we mean by the derived Fukaya category of the mirror of $X$.  Recall from the first part of the introduction that the mirrors of Fano varieties are superpotentials 
\[ W \co Y \to \bC .\]
Kontsevich suggested that correct notion of a Fukaya category which takes into account the superpotential has, as objects, Lagrangian submanifolds $L \subset Y$ whose image under $W$ agrees with a line  $\gamma \co [0,\infty) \to \bC$ away from a compact subset.  There is no loss of generality in assuming that the curve $\gamma$ starts at the origin.  Assuming that $W$ is a symplectic Lefschetz fibration (this is a generic situation), it is not hard to see that the subset of $L$ which lies over $\gamma$ can be obtained by parallel transport, along $\gamma$, of the intersection of $L$ with $W^{-1}(0)$.  In particular, we lose no information if we work with compact Lagrangians whose boundary lies on $W^{-1}(0)$.  This is the point of view taken in this paper.

\begin{rem}
Defining the Fukaya category using only holomorphic curves (and no perturbations) does not yield a category because of issues of transversality.  If one is willing to make choices of inhomogeneous perturbations, it is often possible to extract an honest $A_{\infty}$ category from the geometry of the Fukaya category.  We make an alternative choice in this paper, working with the notion of pre-categories, which we briefly discuss in Appendix \ref{precat}.  The same remark applies to the Morse pre-categories which we shall be using.
\end{rem}
\begin{figure}
\begin{center}
$\begin{array}{c@{\hspace{.5in}}c}
\epsfxsize=2in
\epsffile{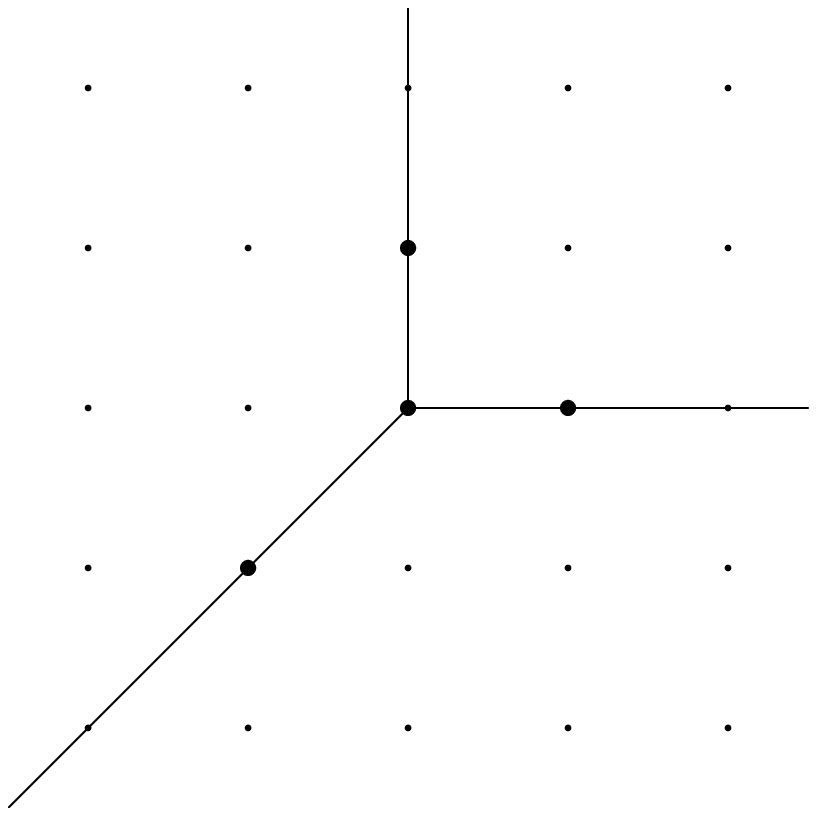} &
\epsfxsize=2in
\epsffile{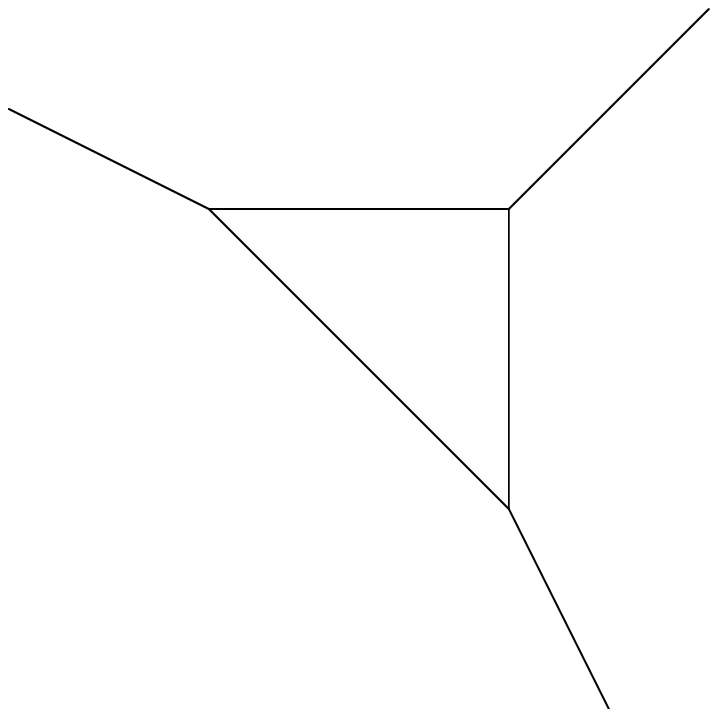}
\end{array}$

\end{center}
\caption{The fan defining $\bC \bP^2$ and a corresponding tropical amoeba.} 
 \label{right}
\end{figure}
From now on, we consider the case where $X$ is a smooth projective toric variety.  The mirror of such a variety is expected to be a Laurent polynomial $W$ on $\ctorus{n}$ which can be explicitly obtained from the fan $\Delta$ which defines $X$ as follows:  

Consider the $1$-cones of $\Delta$.  These are half rays starting at the origin which have rational slope, so they intersect the lattice at a unique primitive (non-divisible) vector.  Let $A$ be the set which consists of the primitive vertices of these $1$-cones and of the origin.  The mirror of $X$ is the family of Laurent polynomials
\[ \sum_{\alpha \in A} c_{\alpha} z^{\alpha},\]
where $z^{\alpha} = z^{a_1} \cdots z^{a_{n}}$ whenever  $\alpha=(a_1, \ldots, a_n)$.  From our point of view the specific coefficients $c_{\alpha}$ are irrelevant as long as they are chosen generically; the corresponding pairs $(\ctorus{n}, W^{-1}(0))$ are symplectomorphic.  In \cite{abouzaid} we consider a family of Lefschetz fibrations $W_{t,s}$ which are not given by Laurent polynomials for $s \neq 0$, but for which the pairs $  (\ctorus{n}, W_{t,s}^{-1}(0))$  are still symplectomorphic for different values of $s$ and $t$.  The purpose of these deformations away from Laurent polynomials is to obtain precise control for the convergence of the hypersurfaces to certain piecewise linear symplectic hypersurfaces constructed from tropical geometry.  From now on, we write $M$ for any of these hypersurfaces $W_{t,s}^{-1}(0)$.  

Given such a pair, we define a (relative) Fukaya pre-category $\Fuk(\ctorus{n},M)$ consisting of compact Lagrangians on $\ctorus{n}$ whose boundary lies in $M$ as described above.  Some details are given in Section \ref{fuk-LG}. In this paper, we prove
\begin{thm} \label{main}
There exists a full sub-pre-category of $\Fuk(\ctorus{n},M)$ which is quasi-equivalent as an $A_{\infty}$ pre-category to the category of line bundles on $X$.
\end{thm}

In order for the statement to make sense, we can equip the category of line bundle with the structure of a $DG$ category by using the (canonical) cover by affine toric subvarieties and passing through \v{C}ech cohomology in order to compute the cohomology of holomorphic line bundles.

As one might expect from Conjecture \ref{section=line} the relevant sub-pre-category of $\Fuk(\ctorus{n},M)$ is obtained by studying Lagrangian sections of a torus fibration
\[ \ctorus{n} \to \bR^{n} .\]
However, the sections are only taken over a subset of the base $\bR^{n}$.  For example, the Lagrangians which are mirror to line bundles on $\bC \bP^2$ are essentially sections over the triangle which appears in the picture on the right of Figure \ref{right}.  The reader might recognize this triangle as the moment polytope of $\bC \bP^2$, but this is not how it appears in our construction.  Rather, following Mikhalkin \cite{mikhalkin}, the graph on the right of Figure \ref{right} is the tropical amoeba of the one parameter family of Laurent polynomials $W_t(x,y) = 1 +t^{-1} x + t^{-1}y + \frac{t^{-1}}{xy}$ which is the Landau-Ginzburg mirror of $\bC \bP^{2}$.

\begin{figure}
\begin{center}

$\begin{array}{c@{\hspace{.5in}}c}
\epsfxsize=2in
\epsffile{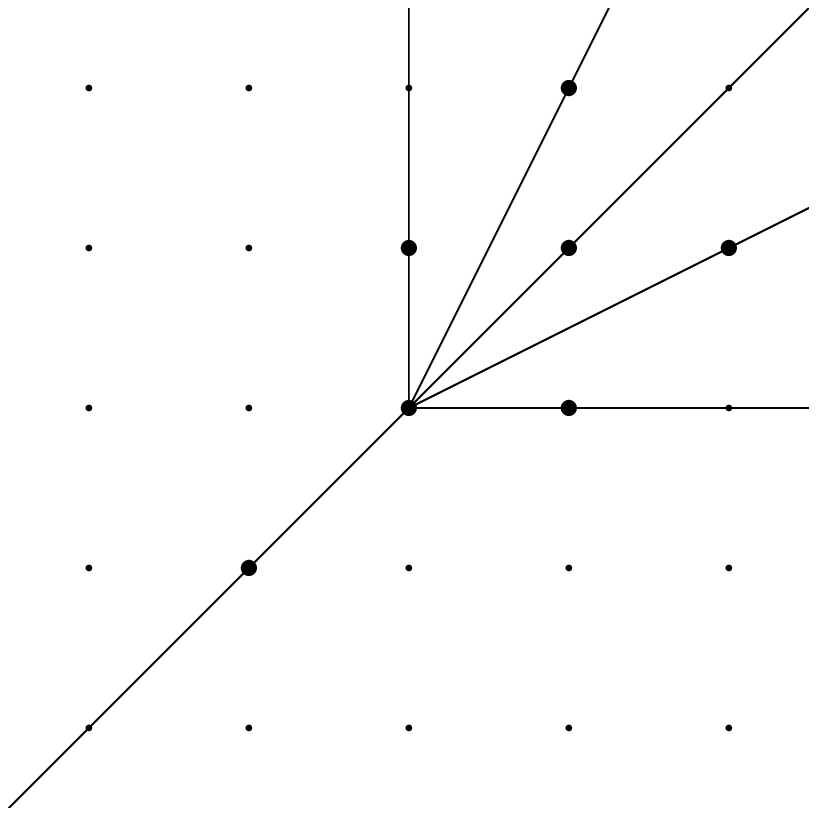} &
\epsfxsize=2in
\epsffile{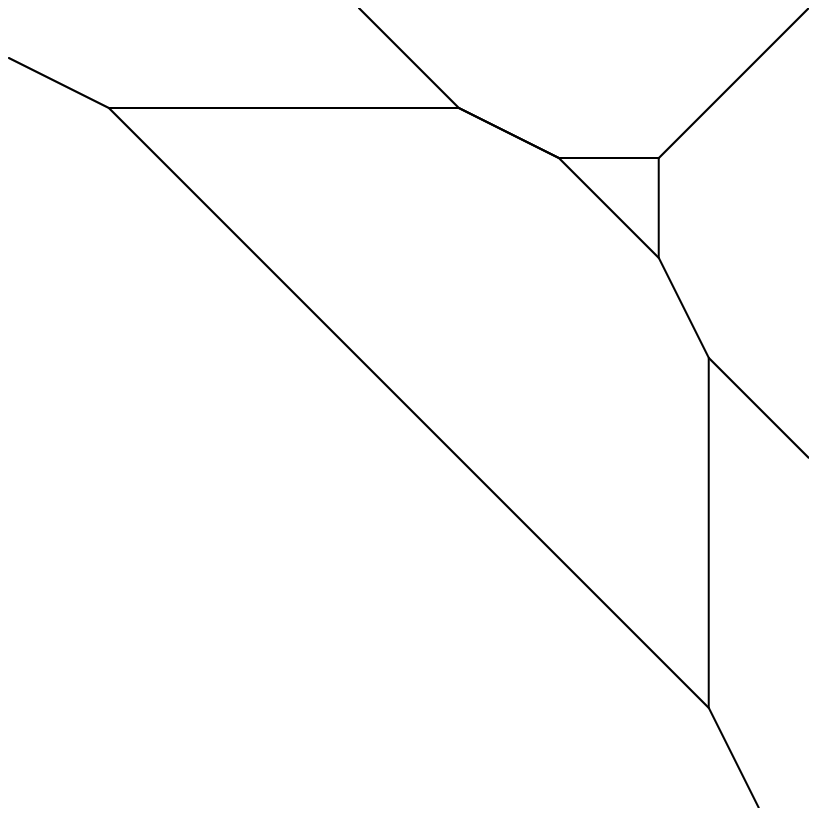}
\end{array}$
\end{center}
\caption{The fan defining an iterated blowup of $\bC \bP^2$ and a corresponding tropical amoeba.}
 \label{wrong}
\end{figure}

The differences with Kontsevich's conjecture are that we prove the result before passing to any derived category and that the toric variety $X$ is not required to have ample anti-canonical bundle.  As Figure \ref{wrong} illustrates, one of the consequences of dropping the Fano condition is that, while the moment polytope of the smooth toric variety still appears as a component of the complement of the tropical amoeba, there may be other bounded polytopes.  In Figure \ref{wrong}, Lagrangian sections over the small triangle form a sub-pre-category of $\Fuk(\ctorus{n},M)$ whose objects do not correspond to any line bundles on $X$ (in fact, with a bit more work, one can see that such sections do not correspond to any object of $D^b\Coh(X)$).

In any case, Seidel observed that one can use a theorem of Thomason to prove the following result:
\begin{prop}
If $X$ is a smooth toric variety, then $D^b \Coh(X)$ is generated by line bundles.
\end{prop}
\begin{proof}
In \cite{thomason}*{Theorem 2.1}, Thomason proves that there is a bijective correspondence between subgroups of $K_0(X)$ (the algebraic $K$-theory of $X$), and full triangulated subcategories $\cA$ of $D^b \Coh(X)$ whose idempotent closure is $D^b \Coh(X)$ (i.e. such that every object of $D^b \Coh(X)$ is a summand of an object in $\cA$).  By an observation of Kontsevich, the full triangulated category generated by line bundles satisfies this second condition if $X$ is smooth.  Indeed, we know that every coherent sheaf has an (infinite) resolution by direct sums of lines bundles.  Consider the truncation of such a resolution at the $n+1$'st position.  Since the category of coherent sheaves on a smooth variety of dimension $n$ has cohomological dimension $n$, there are no non-trivial extensions of length greater than $n$, so that this extension splits (in the derived category).

To complete the proof of the proposition, it suffices to know that $K_0(X)$ is generated by line bundles.  This is a classical result in toric geometry. A stronger result on equivariant $K$-theory of toric varieties is, for example, proved in \cite{morelli}*{Proposition 3}.
\end{proof}

We can therefore use Theorem \ref{main} to conclude

\begin{cor}
$D^{b}\Coh(X)$ is equivalent to a full subcategory of $D^{\pi}\Fuk(\ctorus{n},M)$.  \noproof
\end{cor}

\begin{rem}
This paper is essentially a continuation of \cite{abouzaid} whose results and techniques we use extensively.  The major difference is that, in \cite{abouzaid}, we were able to explicitly identify certain moduli spaces of discs when we restricted ourselves to a fixed ample line bundle and its powers, whereas we are unable to do this in the general case.  Instead, as explained in the next section, we rely on the equivalence between Floer and Morse homology in order to prove the results of this paper.
\end{rem}

The expert reader will notice that we do not rely on the notion of vanishing cycles of Lefschetz fibrations to construct objects of the Fukaya category.  This distinguishes our approach from that taken by Auroux, Katzarkov, and Orlov \cites{AKO1, AKO2} in their proof of the mirror conjecture for Fano surfaces.  On the other hand, using techniques somewhat closer to ours in the sense that they rely on Morse theory, Bondal and Ruan have announced a similar result for weighted projective spaces, see \cite{BR}.   More recently, Fang, Liu, Treumann, and Zaslow \cite{FLTZ} have given a different proof of Homological Mirror symmetry relying on the correspondence established by Nadler and Zaslow between constructible sheaves on the base of a compact manifold and Lagrangians in its cotangent bundle.  Their approach is in a sense dual to the one we've taken.

\subsubsection*{Acknowledgments}
I would first like to thank Paul Seidel for guidance and encouragement.  Conversations with Mitya Boyarchenko, Ludmil Katzarkov, Conan Leung, Peter May, and Ivan Smith helped me connect ideas together and avoid dangerous pitfalls.  Discussions with Kevin Costello, and his comments on an early draft, were particularly helpful.  Finally, I received extensive and helpful comments from Denis Auroux, Mark Gross, and an anonymous referee.

This work was supported by the 2006 Highbridge Research Prize and was partially conducted during the period the author served as Clay Research Fellow.

\section{Outline and discussion of the proof}

In logical order, the proof proceeds as follows:
\begin{enumerate}
\item Identify the appropriate sub-pre-category of $\Fuk(\ctorus{n},M)$.
\item Prove that this sub-pre-category is quasi-isomorphic to a pre-category with essentially the same objects, but whose morphisms can be computed by studying Morse theory on the moment polytope $Q$ of $X$.  We denote this pre-category $\Morse^{\bZ^{n}}(\sQ)$, where $\sQ$ is described below.
\item Prove that the pre-category  $\Morse^{\bZ^{n}}(\sQ)$ is itself equivalent to a pre-category $\Rel^{\bZ^{n}}(Q_b)$ whose morphisms are certain subcomplexes of the cochain complex of $Q_b$, the barycentric subdivision of $Q$.
\item Identify $\Rel^{\bZ^{n}}(Q_b)$ with the category $\Cech(X)$ of line bundles on $X$.
\end{enumerate}

Implicit in the above description is the fact that each of the maps we construct are in fact $A_{\infty}$ functors.  Let us forget for a moment about the $A_{\infty}$ structure, and remember only the structures at the level of homology.  The previous steps become
\begin{enumerate}
\item The appropriate sub-pre-category of $\Fuk(\ctorus{n},M)$ was already implicit in \cite{abouzaid}, particularly in Section 5.  It essentially consists of sections of a torus bundle over a smoothing $\sQ$ of the moment polytope, with appropriate boundary conditions.  The universal cover of this symplectic manifold is the cotangent bundle of $\sQ$.
\item Floer proved in \cite{floer} that Morse homology agrees with Floer homology of exact sections of the cotangent bundle.  Extending the result to the non-compact setting is only a matter of having appropriate definitions for Floer homology (see again \cite{abouzaid} where we follow ideas of Kontsevich and give the relevant example of an appropriate definition of Floer homology of Lagrangians with boundary on a complex hypersurface).
\item The correspondence between Morse homology and ordinary homology is due to Morse.  The earliest appearance of such an idea in the case of manifolds with boundary is probably \cite{BM}.  Because of the specific Morse complex considered, the author was unable to find a reference which proved exactly the desired statement at the level of homology, although nothing original needs to be performed in order to achieve this, so it seems appropriate to simply credit Morse and Baiada.
\item The fact that one can interpret the \v{C}ech cohomology of line bundles on a toric variety in terms of the ordinary homology of the corresponding fan is well known among toric geometers and appears in standard textbooks such as \cite{fulton}.  We are taking the dual point of view by working with a moment polytope.
\end{enumerate}

In particular, it is clear that at the level of homology, every step after the first is a generalization of results that are well established.  We now describe, in more detail, how to proceed in the $A_{\infty}$ setting.
\begin{enumerate}   \setcounter{enumi}{-1}
\item To formalize the issues of transversality (rather, lack thereof), we use unital $A_{\infty}$ pre-categories as a technical tool.  Conjecturally, see \cite{KS}, the homotopy category of such objects is equivalent to that of $A_{\infty}$ categories, although we prove, and use, no result of this type.
\item Consider the Laurent polynomial 
\[ W \co \ctorus{n} \to \bC \]
which is conjecturally the mirror of $X$, and the map
\[\Log \co \ctorus{n} \to \bR^{n} \]
which takes an $n$-tuple of complex numbers to the $n$-tuple of logarithms of their respective norms.  Using the tropical techniques of Mikhalkin, we can pretend that the amoeba of $M$, its image under the $\Log$ map, is a smooth manifold with boundary which is close to a polyhedral complex $\Pi$ called the tropical amoeba.  It was observed in \cite{abouzaid}*{Section 3.2} that, for appropriate choices of tropical degenerations, the moment polytope $Q$ of $X$ appears  as a component of the complement of $\Pi$.  In particular, there is a smoothing $\sQ$ of $Q$ which is a component of the complement of the amoeba of $M$.

Let $\sT \Fuk(\ctorus{n}, M)$ be the sub-pre-category of $\Fuk(\ctorus{n}, M)$ whose objects are Lagrangian sections of the restriction of $\Log$ to $\sQ$ such that their boundary is a subset of $M$, and which are admissible (see Definition \ref{admissibility}).
\item Since our pre-category  $\sT \Fuk(\ctorus{n}, M)$ consists of sections of a Lagrangian torus fibrations, we can compute Floer cohomology groups by passing to the universal cover and taking all possible lifts, remembering the grading by $\bZ^{n} = \pi_1(T^{n})$ .  Since the universal cover is essentially the cotangent bundle of the disc, we can use the results of  \cite{FO} which prove that for compact manifolds without boundary, the moduli spaces of holomorphic discs and of gradient trees agree (exactly!) once the number of marked points is fixed, and the complex structure on the cotangent bundle is appropriately rescaled.  It suffices therefore to develop the correct versions of Floer and Morse theory in the case of manifolds with boundary, and to address the fact that Fukaya and Oh do not prove the existence of an almost complex structure which makes the moduli spaces of holomorphic discs and gradient trees agree regardless of the number of marked points (the author knows of no reason for such an almost complex structure to exist).  The relevant Morse and Floer theories for manifolds with boundary are developed in this paper, while we use some techniques of homological algebra in Section \ref{floer-morse} to construct $A_{\infty}$ quasi-isomorphisms relating the Fukaya pre-category and the Morse pre-category  $\Morse^{\bZ^{n}}(\sQ)$ despite the lack of an almost complex structure which would make the two pre-categories isomorphic.

\item Every smooth function on $\sQ$ yields a decomposition of the boundary of $\sQ$ into a subset $\partial^{+} \sQ$ where the function is increasing, and another where it is decreasing. The Morse cohomology of such a function is in fact computing the relative cohomology group
\[ H^{*}(\sQ, \partial^{+} \sQ) .\]
Further, we define moduli spaces of trees which allow us to relate the $A_{\infty}$ structure on Morse functions to the cup product on the cellular chains of an appropriate triangulation.  A key property of the functions which appear in  $\Morse^{\bZ^{n}}(\sQ)$ is that the subset $\partial^{+} \sQ$ essentially contains a union of cells of the polytope $\partial Q$ as a deformation retract.  This means that this cohomology group can also be understood combinatorially.   By interpreting the chains of the dual barycentric subdivision of the natural cell decomposition of $Q$ as cochains on $Q_b$, and studying moduli spaces of gradient trees with boundary conditions on cells of this barycentric subdivision, we construct an $A_{\infty}$ quasi-isomorphism from a combinatorial pre-category $\Rel^{\bZ^{n}}(Q_b)$ to  $\Morse^{\bZ^{n}}(\sQ)$.

\item We strengthen the classical result of toric geometry by proving that the category $\Rel^{\bZ^{n}}(Q_b)$ is equivalent as a DG category to the category $\Cech(X)$ of line bundles on $X$ with morphisms given by the \v{C}ech complexes and composition given by the cup product.
\end{enumerate}

\section{The category of tropical Lagrangian sections}

\subsection{The Fukaya pre-category of a superpotential} \label{fuk-LG}

Strictly speaking we only defined the Donaldson-Fukaya pre-category of admissible Lagrangians in \cite{abouzaid}; i.e. we only studied the cup product at the level of homology. Before proceeding to the case of Landau-Ginzburg potentials, let us recall the definition of the Fukaya pre-category in the absence of a superpotential.  We will work with manifolds whose symplectic form  is exact which greatly simplifies the technical difficulties; for an overview of the theory, see \cite{seidel-book}*{Chapter 2}.  Let $N$ be a complete Weinstein manifold, in other words, $N$ is equipped with both a symplectic form $\omega$ and a $1$-form $\lambda$ such that
\[ d \lambda = \omega ,\]
and the corresponding Liouville flow is a gradient flow for some proper Morse function on $N$.  For details, see \cite{eliashberg-gromov}. Assume further that the first Chern class of $N$ vanishes.  Choose an almost complex structure on $N$ which is convex at infinity, and a complex volume form, i.e. a trivialization of the (complex) canonical bundle. The reader should have in mind the example where $N$ is the cotangent bundle of a Riemannian manifold $V$, equipped with canonical $1$-form $p dq$ and the complex volume form induced by the real volume form on $V$.

The Fukaya category of $N$ should have, as objects, Lagrangian submanifolds of $N$.  However, even though $N$ is exact, all difficulties do not disappear, so we will restrict to Lagrangians $L$ which satisfy the following conditions:
\begin{itemize}
\item The submanifold $L$ is equipped with a spin structure.  This assumption is necessary in order to obtain natural orientations on the moduli spaces of holomorphic polygons which will be used in the definition of the Fukaya category.  In particular, we can only part with this condition if we are willing to work only in characteristic $2$,
\item The restriction of $\lambda$ to $L$ is exact.  By a standard argument involving Stokes' theorem, this implies that there are no holomorphic discs with boundary on $L$, which in turn implies that the differential defined on Floer complexes is an honest differential; i.e. its square vanishes.
\item We are given a choice of grading on $L$.   Gradings are necessary in order for the Floer chain complex to be, as one might guess, $\bZ$-graded.  Without such a choice, they are only naturally $\bZ/2\bZ$ graded.  A complex volume form on $N$ induces a real structure on its canonical bundle.  At every point on $L$, the tangent space $TL$ determines a real line in the canonical bundle.  A grading is a choice of an $\bR$ valued function whose exponential gives the $S^1$ valued angle between this line and the $x$-axis.  Of course, if $H_1(L) = 0$, such a grading always exists.   For a thorough discussion of grading issues in Floer theory, see \cite{seidelGL}. 
\end{itemize}
\begin{defin}
A {\bf Lagrangian brane} $L \subset N$ is an embedded compact Lagrangian submanifold equipped with the above choices of additional structures.
\end{defin}

Given transverse Lagrangian branes $L_1$ and $L_2$ in $N$, let $\CF^*(L_1, L_2)$ be the graded free abelian group generated by the intersection points $L_1 \cap L_2$.  Again, the reader can refer to \cite{seidelGL} for how to assign gradings to intersection points.  If $L_1$ and $L_2$ are graphs of exact $1$-forms in the cotangent bundle of a manifold $V$ of dimension $n$, an intersection point $p$ corresponds to a critical point of the difference between primitives $f_1$ and $f_2$ for the $1$-forms defining $L_1$ and $L_2$.  In this case, 
\[ \deg(p) = n - \mu(p) ,\]
where $\mu(p)$ is the Morse index of $p$ as a critical point of $f_2 - f_1$.

Given an almost complex structure $J$ on $N$, Floer considered the moduli space of holomorphic maps
\[ u \co \bR \times [0,1] \to N \]
which, for all $t$, take $(t,0)$ to $L_1$ and $(t,1)$ to $L_2$.  If $u$ has finite energy, he proved that $u(t,s)$ must converge to intersection points $p$ and $q$ between $L_1$ and $L_2$ as $t$ converges to $\pm \infty$.  We denote the moduli space of such maps by
\[ \cM(q,p) .\]
Translation in the $\bR$ factor of the domain defines a free $\bR$ action on this space.  We write
\[ \cM(q,p)/\bR \]
for the quotient.   If $L_1$ and $L_2$ are generic (a $C^{\infty}$ small Hamiltonian perturbation of one of the Lagrangians suffices, see \cite{oh-perturb}) we can now define an endomorphism of  $\CF^*(L_1, L_2)$ by the formula
\[ dp = \sum_{q  \ni \dim(\cM(q,p)) = 1} |  \cM(q,p)/\bR| \cdot q. \]
In the above formula, $ |  \cM(q,p)/\bR|$ is a signed count of the number of components of $\cM(q,p)$.  Since we have discussed neither orientations on $\cM(q,p)$ nor gradings on the Floer complex, we refer to \cite{seidel-book} for the sign considerations which are needed to prove of the following Lemma.   Without signs, this result is due to Floer \cite{floer}, and is by now a standard application of Gromov's compactness theorem \cite{gromov}.
\begin{lem}
The map $d$ is a differential on  $\CF^*(L_1, L_2)$. \noproof
\end{lem}

One would like to define compositions
\[\CF^*(L_2, L_3) \otimes \CF^*(L_1, L_2) \to \CF^*(L_1, L_3)\]
by using moduli spaces of curves with boundary on $L_1$, $L_2$ and $L_3$.  The resulting composition will, however, only be associative after passing to homology.  Fukaya observed that even at the chain level, one can construct a strongly homotopy associative composition by considering all moduli spaces of holomorphic polygons.  We briefly describe his construction:

Let $(L_0, \cdots, L_n)$ be a sequence of Lagrangian branes in general position, and let $S$ be a collections of points $(s_{0,1}, s_{1,2}, \ldots, s_{d-1,d}, s_{0,d})$ on the boundary of the closed disc $D^2$ ordered counter-clockwise. We consider finite energy holomorphic maps
\[ u: D^2 - S \to N \]
such that the image of the interval from $s_{i-1,i}$ to $s_{i,i+1}$ lies on $L_i$.  The finite energy condition implies that $u$ extends to a continuous map on $D^2$ taking $s_{i,i+1}$ to an intersection point $p_{i,i+1}$ between $L_i$ and $L_{i+1}$. We define
\[ \cM(p_{0,d}; p_{0,1}, \ldots,  p_{d,d-1} ) \]
to be the moduli space of such maps.  Again, the choices of spin structures on $L_i$ induce an orientation on this moduli space.  Further, assuming the Lagrangians are in generic position, these moduli spaces are smooth manifolds of dimension
\begin{equation} d -2 + \deg(p_{0,d}) - \sum_{i} \deg(p_{i,i+1}) .\end{equation}
The genericity condition used here is the assumption that moduli spaces of holomorphic polygons are regular, see \cite{FHS}.

Fukaya defined an operation $\mu_d$ of degree $2-d$ 
\[ \mu_d \co \CF^{*}(L_{d-1}, L_{d}) \otimes \cdots \otimes \CF^{*}(L_0,L_1) \to \CF^{*}(L_0, L_d).\]
On generators, $\mu_d$ is given by the formula
\begin{equation} p_{d-1,d} \otimes \cdots \otimes p_{0,1} \mapsto \sum_{\deg(p_{0,d}) = 2-d +  \sum_{i} \deg(p_{i,i+1})} |\cM(p_{0,d};  p_{0,1}, \ldots,  p_{d,d-1})| p_{0,d},\end{equation}
where $|\cM(p_{0,d}; p_{0,1}, \ldots,  p_{d,d-1} )|$ is the signed number of components of the zero-dimensional manifold $\cM(p_{0,d};  p_{0,1}, \ldots,  p_{d,d-1})$.  Gromov's compactness theorem implies the following result, whose proof with signs is given in \cite{seidel-book}*{Chapter 2}.
\begin{lem} \label{lem_a_infty_fuk}
The operations $\mu_d$ satisfy the $A_{\infty}$ equation.   Explicitly, we have
\begin{equation} \label{A_infty-equation} \sum_{i+d_2 \leq d} (-1)^{\maltese_i} \mu_{d-d_2+1}(\mathbb{1}^{d-d_2-i} \otimes \mu_{d_2} \otimes \mathbb{1}^{i}) =0 ,\end{equation} 
where $\maltese_{i}$ depends on the degree of the first $i$ inputs as follows:
\[ \maltese_i = i + \sum_{j=0}^{i-1} \deg(p_{j,j+1}) .\] \noproof
\end{lem}

We can summarize the above discussion in the following definition.  We refer the reader to Appendix \ref{precat} for the definition of $A_{\infty}$ pre-categories.
\begin{defin}
The {\bf Fukaya pre-category} of a complete Weinstein manifold $N$ is the $A_{\infty}$ pre-category $\Fuk(N; J)$ given by the following data.
\begin{itemize}
\item Objects are Lagrangian branes in $N$.
\item A sequence $(L_0, \cdots, L_d)$ is transverse whenever all moduli spaces of holomorphic polygons with boundary on a subsequence are regular.
\item The space of morphisms between a transverse pair $(L_0, L_1)$ is the Floer complex $\CF^*(L_0, L_1)$.
\item The higher compositions are given by the maps $\mu_d$.
\end{itemize}
\end{defin}
\begin{rem}
We will often drop $J$ from the notation, and simply write $\Fuk(N)$
\end{rem}
Floer's original applications of Floer homology relied on its invariance under Hamiltonian isotopies.  Using the composition $\mu_2$, one can state a stronger result as follows.  If $L'_0$ is a Hamiltonian deformation of $L_0$, then there exists a cycle $e \in \CF^0(L_0, L'_0)$ such that
\[ \mu_2(e, \_) \co  \CF^*(L_1, L_0) \to \CF^*(L_1, L'_0) \]
induces an isomorphism on cohomology whenever $(L_1, L_0, L'_0)$ is a transverse triple in the Fukaya pre-category.  Reversing the roles of $L_0$ and $L'_0$, we obtain the following result.
\begin{lem}
The pre-category $\Fuk(N)$ is unital. \noproof
\end{lem}
\subsubsection{The Fukaya Category of a Landau-Ginzburg potential}
Consider a map $f:N \to \bC$ which is a symplectic fibration away from its critical values and assume that the zero fibre is a smooth manifold $M$.  Assume further that we have chosen an almost complex structure on $N$ which respect to which $f$ is holomorphic in a neighbourhood of $M$. The following definition appears in \cite{abouzaid}*{Section 2}, and is a reformulation of a definition of Kontsevich in \cite{kont-ENS}.

\begin{defin} \label{admissibility}
An {\bf admissible Lagrangian  brane} is a compact Lagrangian submanifold $\cL \subset N$ with boundary on $M$ together with a choice of grading and spin structure, satisfying the following conditions:
\begin{itemize}
\item The restriction of $\lambda$ to $\cL$ is exact, and
\item There exists a small neighbourhood of $\partial \cL$ in $\cL$ which agrees with the parallel transport of $\partial  \cL$ along a segment in $\bC$.  
\end{itemize}
\end{defin}
Note that the last two conditions are vacuous if $\partial \cL$ is empty, which is a possibility that we do not exclude.

\begin{figure}[h] 
   \centering
 \input{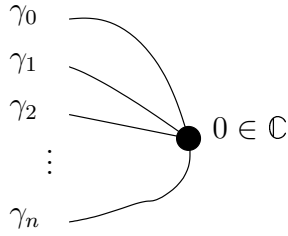}
   \caption{Curves corresponding to positively oriented Lagrangians.}
   \label{boundary_orient}
\end{figure}

Next, we consider the conditions under which the higher products are well behaved.
\begin{defin}
A sequence of admissible Lagrangians $( \cL_0, \ldots, \cL_n )$ is {\bf positive} if the curves $\gamma_k$ which agree with $f(\cL_k)$ near the origin all lie in the left half plane and their tangent vectors $v_k$ are oriented counter-clockwise (See Figure \ref{boundary_orient}).
\end{defin}
Note that if $( \cL_i )_{i=1}^{n}$ is a positive sequence then the angle between $v_i$ and $v_j$ is again positive whenever $i < j$.  In particular, every subsequence of a positive sequence is positive.

Let us assume that the Lagrangians $\cL_i$ intersect transversally in the interior.  Throughout, we will only work with interior intersection points.  If $p_{i} \in \cL_{i-1} \cap \cL_i$ are such interior intersection points between successive Lagrangians, and $p_0 \in \cL_{d} \cap \cL_0$, we can define
\[ \cM(p_0; p_1, \ldots, p_d) \]
to be the moduli space of holomorphic polygons with boundary on $\bigcup_{i} \cL_i $ and corners at $p_i$.   If we perturb each Lagrangian $\cL_i$ near the intersections by a $C^{\infty}$ small amount, we can guarantee that all moduli space of holomorphic polygons are regular, see \cite{FHS}.  This allows us to define the Gromov bordification
\[ \overline{\cM}(p_0; p_1, \ldots, p_d) \]
by including ``cusp-curves''.
\begin{lem}  \label{admissible_implies_compactness}
If $\{ \cL_i \}$ is positive, then the Gromov bordification is compact.
\end{lem}

\begin{proof}[Proof (c.f. Lemma \ref{lem:straight_near_circle_compact})]
It suffices to prove that curves in \[ \cM(p_0; p_1, \ldots, p_d) \] cannot escape to the boundary.  The proof given in \cite{abouzaid}*{Lemma 2.8} can be used here. The key ideas is to project by $f$ to $\bC$, and observe that if $u_s$ is a family of holomorphic maps parameterized by $[0,1]$ such that the image of $u_0$ does not intersect a neighbourhood of the origin, then the same property holds for $u_1$.
\end{proof}
If we define the higher products $\mu_d$ by counting the number of isolated holomorphic polygons as in the case where the Lagrangians have no boundary, we conclude
\begin{cor}
If all moduli spaces are regular, then the $A_{\infty}$ equation holds. \noproof
\end{cor}

We are now ready to define the main object which we shall be studying.
\begin{defin}
The {\bf Fukaya  pre-category} of the Landau-Ginzburg potential $f$ is the $A_{\infty}$ pre-category $\Fuk(N,M)$ given by the following data:
\begin{itemize}
\item Objects are admissible Lagrangian branes.
\item Transverse sequences are positive sequences of Lagrangian branes such that all holomorphic polygons with boundary values on a sub-sequence are regular.
\item Morphisms are the Floer complexes with the differential counting holomorphic bi-gons and higher products counting holomorphic polygons.
\end{itemize}
\end{defin}

\subsection{Tropical Lagrangians for mirrors of toric varieties} \label{tropical-lagr}
Let $X$ be a smooth toric variety defined by a fan $\Delta$.  We assume $X$ to be projective, so $|\Delta| = \bR^{n}$ and there exists an integral piece-wise linear function 
\[ \psi: \bR^{n} \to \bR\]  
which is linear on each cone of $\Delta$, and which is strictly convex.  Let $Q$ be the moment polytope of the line bundle which corresponds to $\psi$.  In \cite{abouzaid}*{Definition 4.3}, we constructed a cover $\{O_{\tau} \}$ of $\bR^{n}$ labeled by the faces of $Q$ such that the complement in $\tau$ of a small neighbourhood of $\partial \tau$ is contained in $O_{\tau}$.

\begin{rem}
We are about to extensively use results from \cite{abouzaid}.  The notation has been kept relatively consistent, with one major difference.  In \cite{abouzaid}, we used the greek letters $\sigma$ and $\tau$ to label cones in the fan $\Delta$.  In this paper, we use the same letter to label the faces of the polytope $Q$.  The open sets $O_{\tau}$ therefore correspond  to $O_{\check{\tau}}$ in the notation of \cite{abouzaid}.
\end{rem}

\begin{figure}
\begin{center}
$\begin{array}{c@{\hspace{.1in}}c@{\hspace{.1in}}c} \epsfxsize=1.5in
\epsffile{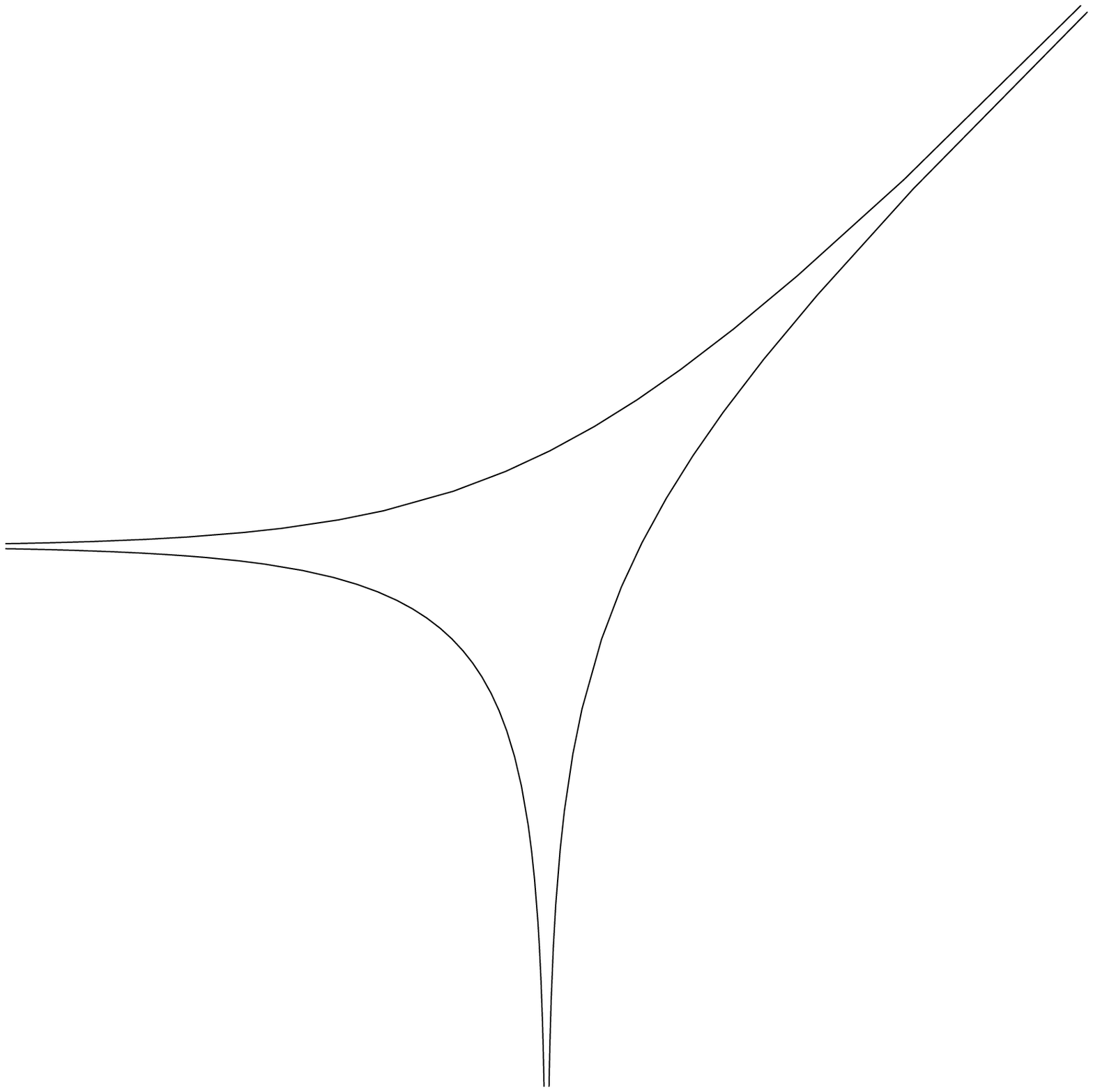} & \epsfxsize=1.5in
\epsffile{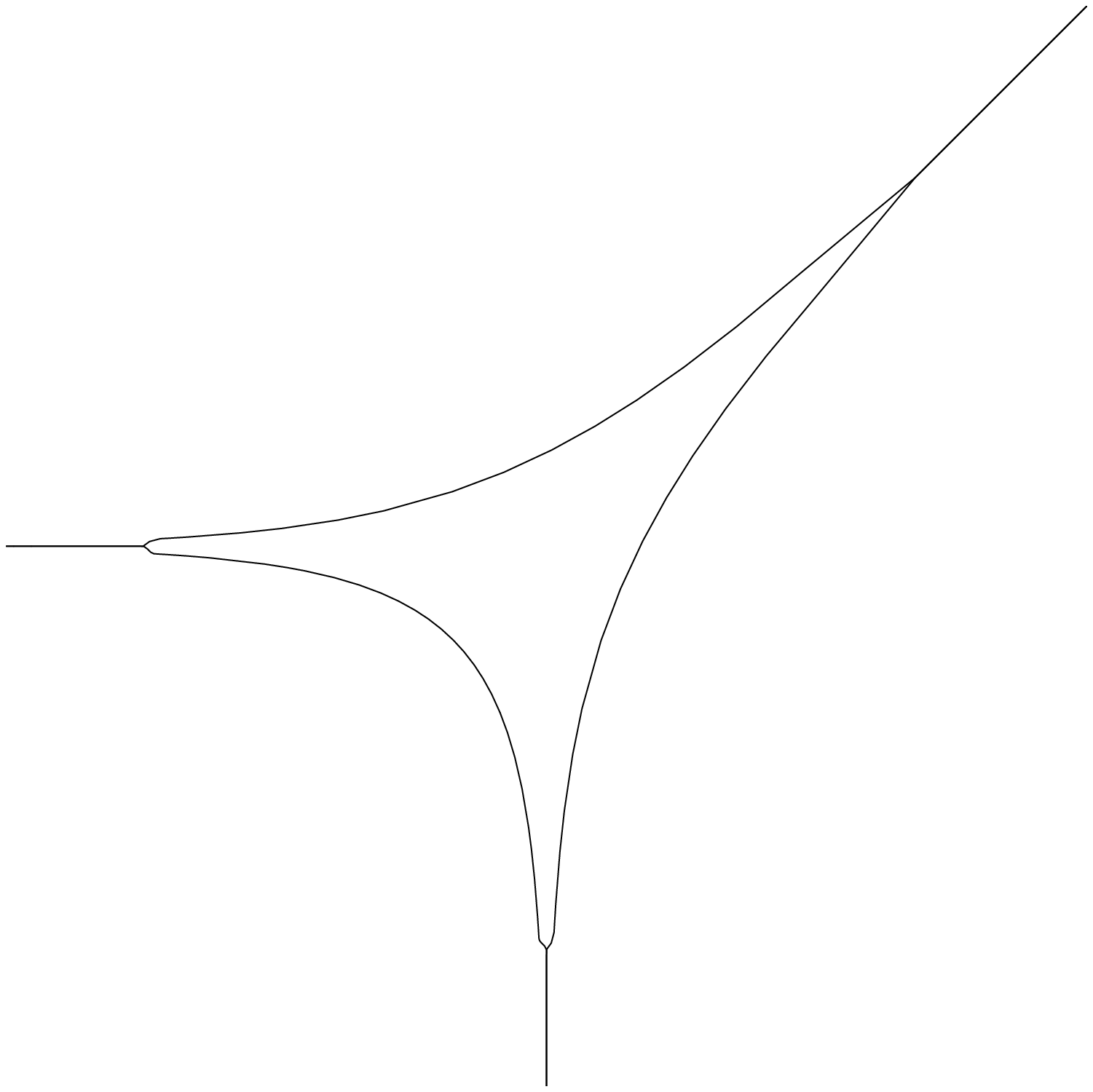} & \epsfxsize=1.5in
\epsffile{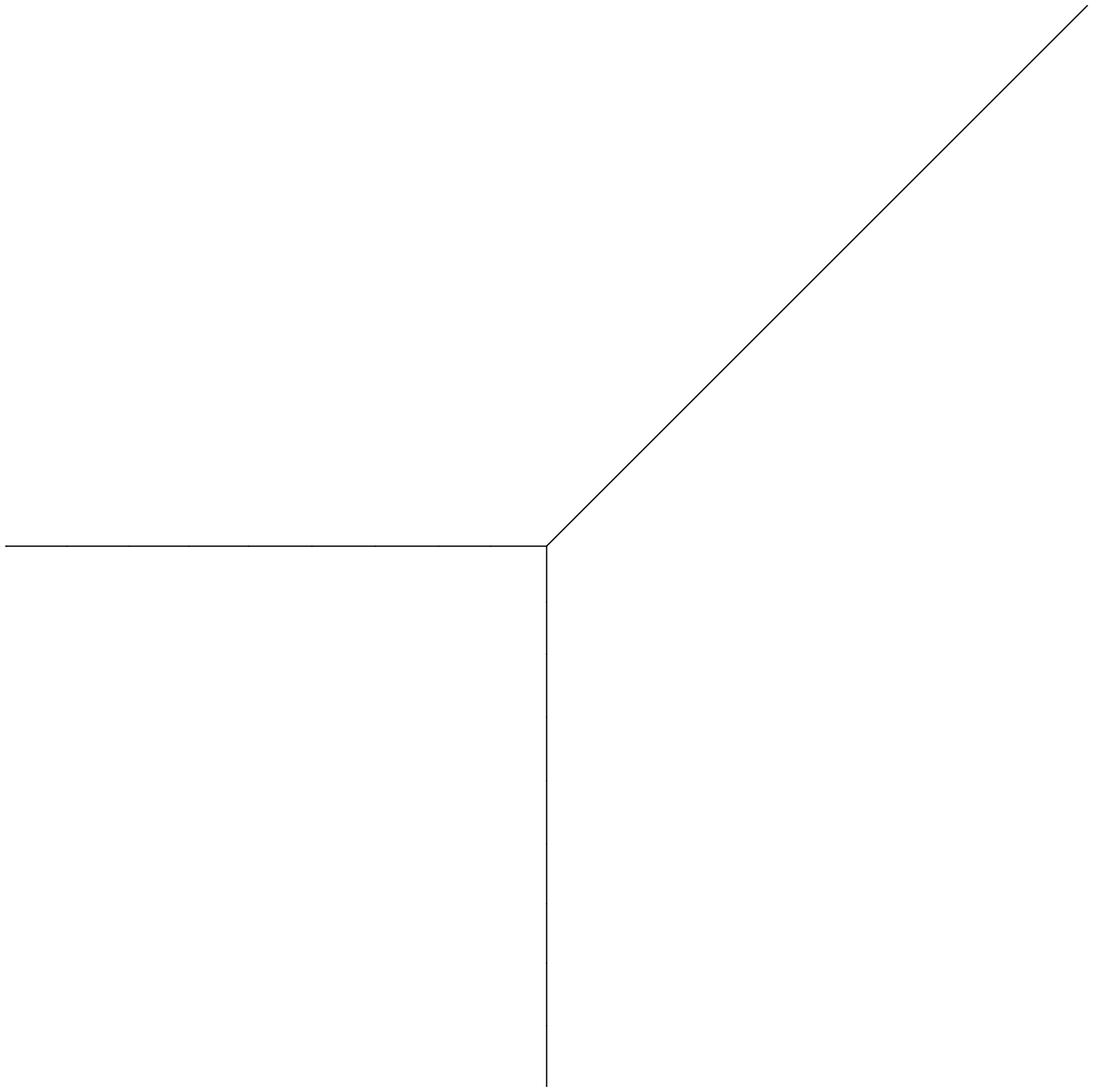}
\end{array}$
\end{center}
\caption{The amoeba of the standard hyperplane in $\ctorus{2}$, its
``tropical localization," and its tropical amoeba}
\label{amoebas-fig}
\end{figure}

We observed that the choice of $\psi$ induces a tropical degeneration of the mirror of $X$.  In particular, we have a family of Laurent polynomials $W_t$, whose zero fibre are (for $t \gg 0$) are smooth complex (and hence symplectic) hypersurfaces which we denote by $M_t \subset \ctorus{n}$.  We considered the amoebas $\cA_t$ of $M_t$, i.e the projection of $M_t$ to $\bR^{n}$ under the logarithmic moment map.  Further, we constructed a ``tropical localization of $M_t$" by deforming $W_t$ to a (non-algebraic) map $W_{t,1}$ so that the resulting hypersurface $M_{t,1}$ is particularly well described by its amoeba.  To understand the relation between the amoebas of $M$, $M_{t,1}$, and the tropical amoeba, it is best to consider the case of the standard hyperplane in $\ctorus{2}$ which is illustrated in Figure \ref{amoebas-fig}.

By construction, the tangent space of any cell $\tau$ is an integral subspace of $\bR^{n}$, i.e. it is defined as the vanishing locus of finitely many integral covectors. In particular, if $T_{\bZ}\tau$ is the intersection of this tangent space with the lattice $\bZ^n$, then we have an identification
\[ T \tau =T_{\bZ}\tau \otimes_{\bZ} \bR .\]

Let us write $\check{\tau}$ for the linear span of the covectors which vanish on $T \tau$, and  $\ctorus{\check{\tau}}$ for the corresponding complex torus.  The natural map
\[ \ctorus{n} \to \ctorus{\check{\tau}} \]
has fibre which can be identified with $\ctorus{\tau} = T \ker(\check{\tau})  / T_{\bZ}  \ker(\check{\tau})$ (the quotient of the tangent space of the linear plane tangent to $\tau$ by the lattice).  Note that the choice of a splitting of the inclusion $T_{\bZ}\tau \into \bZ^n$, gives a product decomposition
\[  \ctorus{n} \cong \ctorus{\check{\tau}} \times \ctorus{\tau} \]
of complex manifolds.  We can now restate, in a slightly enhanced form, a results  which appeared in \cite{abouzaid}*{Lemma 4.4}:
\begin{lem}  \label{local_product}
The intersection of the hypersurface $M_{t,1}$ with the inverse image of $ \log(t) \cdot O_{\tau}$ in $\ctorus{n}$, splits as a product
\[ M_{t,1}^{\check{\tau}} \times (\tau \cap O_{\tau}) \times \bT^{\tau} \]
where 
\begin{itemize}
\item $ M_{t,1}^{\check{\tau}}  $ is the tropical localization of a hypersurface in a torus $\ctorus{\check{\tau}}$, and
\item  $\bT^{\tau} = T \tau / T_{\bZ}\tau$ is the real subtorus of $\ctorus{\tau}$.
\end{itemize}
\noproof
\end{lem}

Unfortunately, this is not a product decomposition as symplectic manifolds.  In order to obtain such a decomposition, we would essentially have to pick the splitting of the inclusion $T_{\bZ}\tau \into \bZ^n$ compatibly with the euclidean metric, which is of course not always possible.  In general, the orthogonal complement of $T_{\bZ} \tau$ in $\bZ^n$ projects to a finite index sub-lattice of the quotient $\bZ^n / T_{\bZ} \tau$.  If we write $\ctorus{\tau^{\perp}}$ for the corresponding torus, we obtain a finite cover
\[ \ctorus{\tau} \times \ctorus{\tau^{\perp}} \to \ctorus{n} .\]
The inverse image of $M_{t,1}$ in the cover splits as a product in both the symplectic and complex categories.  The factor living in $\ctorus{\tau^{\perp}}$ now corresponds to a higher degree hypersurface (rather than a hyperplane).

In \cite{abouzaid}*{Lemma 4.16} we observed that one of the components $\sQ_{t,1}$ of the complement of $\cA_{t,1}$ satisfies the following conditions:
\begin{itemize}
\item $\sQ_{t,1}$ is contained in (and is $C^0$ close to) the $\log(t)$ rescaling of the moment polytope $Q$, and
\item the boundary of $\sQ_{t,1}$ is smooth, and its image under the zero section lies in $M$.
\end{itemize}

Note that the boundary $\sQ_{t,1}$ again decomposes as a product in $O_{\tau}$.  Since we are now working in $\bR^{n}$ rather than in the complex torus, we can choose this product decomposition to be compatible with the metric.  More precisely:
\begin{lem}
The restriction of $\sQ_{t,1}$ to $\log(t) \cdot O_{\tau}$ splits as a Riemannian product 
\[ \sQ^{\tau^{\perp}}_{t,1} \times (\tau \cap O_{\tau}) ,\]
where $\sQ^{\tau^{\perp}}_{t,1}$ is the amoeba of the tropicalisation of a hyperplane.  \noproof
\end{lem}
\begin{rem}
In order to see $\sQ^{\tau^{\perp}}_{t,1}$ as the amoeba of a hyperplane (rather than a higher degree hypersurface) we equip the plane orthogonal to $\tau$ with the lattice coming from its identification with the quotient of $\bR^n$ by the tangent space to $\tau$.
\end{rem}

From now on (except for Section \ref{morse=fuk-tropical} and Appendix \ref{ap:second_look}), we shall drop the subscripts from $M_{t,1}$, $W_{t,1}$, and the associated objects.  We will further rescale $\bR^{n}$ by $\frac{1}{\log(t)}$, and write $\sQ$ for the image of $\sQ_{t,1}$ under this rescaling.  Note that $\sQ$ is $C^0$ close to the original polytope $Q$.

We denoted the image of $\sQ$ under the zero section by $\cL$.  It is supposed to correspond under mirror symmetry to the structure sheaf of $X$.  
\begin{defin}
The pre-category $\sT \Fuk(\ctorus{n}, M)$ of {\bf tropical Lagrangian sections}, is the $A_{\infty}$ pre-category of admissible Lagrangians branes in $(\ctorus{n}, M)$ which are sections of the restriction of the moment map to $\sQ$.
\end{defin}
\begin{rem}
It would be more consistent to use the terminology ``semi-tropical Lagrangian section,'' hinting that we did not take the tropical degeneration to its limit, but rather stopped just before the limit point, with the smooth symplectic structure still present, yet the tropical features apparent.
\end{rem}

Given $\cL_1$ and $\cL_2$ two tropical Lagrangian sections, the space of morphisms $\CF^{i}(\cL_1, \cL_2)$ is naturally graded by an affine space over $\bZ^{n}$.  This is a grading which is independent of the cohomological grading.  Indeed, since we have fixed a zero section $\cL$ and a fixed lift $\tilde{\cL}$ to the universal cover, we can consider the difference $\cL_2 - \cL_1$, which yields a new section of the fibration.  The lifts of $\cL_2 - \cL_1$ to the universal cover form an affine space over $\bZ^{n}$, which we identify with the grading of $\CF^{i}(\cL_1, \cL_2)$.  Given an intersection point $p$ between $\cL_2$ and $\cL_1$, we obtain an intersection point between $\cL_2 - \cL_1$ and the zero section, and the corresponding lift to the universal cover is uniquely determined by the fact that $\tilde{p}$, the lift of $p$, lies on $\tilde{\cL}$.

A priori, there are uncountably many tropical Lagrangian sections.  As we shall see, there are in fact only countably many Hamiltonian isotopy classes, naturally in bijective correspondence with isomorphism classes of line bundles on $X$.  Recall that 
\[ \Log: \ctorus{n} \to \bR^{n} \]
is the moment map taking an $n$-tuple of non-zero complex numbers to the logarithms of their absolute values.  We abusively write   $\Log|_{\sQ}$ for the restriction of this torus fibration to the inverse image of the disc $\sQ$.
\begin{lem}
Every section of $\Log|_{\sQ}$ which corresponds to a tropical Lagrangian section agrees with the zero section when restricted to 
\[ \partial \sQ \cap O_{v} ,\]
whenever $v$ is a vertex of $Q$.
\end{lem}
\begin{proof}
It is easy to check that the boundary of the amoeba of the standard hyperplane $\sum_{i=1}^{n} z_i =1$ agrees with its real locus. In particular, the restriction of the Logarithmic moment map to the standard hyperplane in $\ctorus{n}$ is injective when restricted the inverse image of the boundary of its amoeba.  The results of \cite{abouzaid}*{Lemma 4.4} restated above imply that the amoeba of $M$ is modeled, in a neighbourhood of every vertex of $Q$, after the amoeba of the standard hyperplane.  Since the zero section is a tropical Lagrangian, every tropical Lagrangian section must therefore agree with it on this subset of the boundary.
\end{proof}
In particular, any lift of a tropical Lagrangian section to $T^{*} \bR^{n}$ must take $2 \pi $ integral values.   As an extension of the previous lemma, we have
\begin{lem} \label{classification_tropical}
The lift of any tropical Lagrangian section when restricted to
\[ \partial \sQ \cap O_{\tau} ,\]
lies in an affine subspace of $\bR^{n}$ which is a $2 \pi$-integral translate of the tangent space of $\tau$.  A Lagrangian section which satisfies this property for every $\tau$ is isotopic through Lagrangian sections with boundary on $M_{t,1}$ to a tropical Lagrangian section.
\end{lem}
\begin{proof}
We first prove  that any section (not necessarily Lagrangian) of the restriction of $T^* \bR^{n}$ to $\partial \sQ \cap O_{\tau}$ which lands in $M$ must have values in a $2 \pi$-integral translate of the tangent space of $\tau$.  The previous lemma proves the base case when $\tau$ is a vertex and the induction step follows easily from the the product decomposition of Lemma  \ref{local_product}.

To prove the second part, we must achieve admissibility.  Consider a function $f$ defining a Lagrangian section $\Lambda_f$ with boundary $\partial \Lambda_f \subset M_{t,1}$.  Using the parallel transport maps defined by the symplectic connection induced by $W$, we can construct a Lagrangian $\Lambda_g$ whose boundary agrees with that of $\Lambda_f$.  The arguments of  \cite{abouzaid}*{Section 5} imply that $\Lambda_g$ (near $\partial \sQ_{t,1}$) is defined by a function $g$.  By linearly interpolating between $f$ and $g$, we obtain the desired isotopy.  
\end{proof}
We are now ready to prove
\begin{prop}
The Hamiltonian isotopy class of a tropical Lagrangian section is determined by the values of its lifts near the vertices of $Q$ 
\end{prop}
\begin{proof}
Every such lift is the graph of the differential of a function on the disc.  Let $f_1$ and $f_0$ denote two such functions whose differentials are equal near the vertices. Consider the $1$-parameter family
\[ f_t = t f_0 + (1-t) f_1.\]
From the characterization of the previous lemma, it easily follows that the Lagrangians $\Lambda_{f_t} \subset \ctorus{n}$ defined by $f_t$ have boundary on $M$.  In particular, there is a family $g_t$ of functions whose graphs define admissible Lagrangians with boundaries $\partial \Lambda_{f_t}$.  Picking such a family with $f_0 = g_0$ and $f_1 = g_1$ yields the desired isotopy.
\end{proof}

\begin{cor}
There is a bijective correspondence between Hamiltonian isotopy classes of tropical Lagrangian sections of $(\ctorus{n},M)$ and isomorphism classes of line bundles on $X$. 
\end{cor}
\begin{proof}
It is well known that isomorphism classes of line bundles are given by piece-wise linear integral function on $\bR^{n}$, with domains of linearity the maximal cones of $\Delta$, modulo the (global) addition of an integral linear function on $\bR^n$.  The data of a lattice vector for each vertex of $Q$ corresponds to a linear integral function on the maximal cones of $\Delta$.  The fact that the restriction of our section to an edge of $Q$ can only vary by multiples of the integral tangent vector of the edge implies that the two linear functions associated to the maximal cones dual to the endpoints of an edge agree on the hypersurface where they meet.  Starting with a piecewise linear integral function on $\Delta$, it is easy to construct the corresponding section by a patching argument.  By Lemma \ref {classification_tropical} this Lagrangian with boundary on $M$ is isotopic to an admissible one, hence to a tropical Lagrangian section.  Finally, we note that addition of an integral linear function on $\bR^{n}$ corresponds to choosing different lifts.
\end{proof}

In Sections \ref{fukaya=morse-section} and \ref{morse=cech-section} we will prove that this bijective correspondence extends to an equivalence of $A_{\infty}$ pre-categories.  In order to do this, we shall have to develop some general machinery relating the multiplicative structures on the chain complexes defining Morse, Floer, \v{C}ech, cellular, and simplicial cohomology theories.

\section{General framework}

In this section, we construct chain-level $A_{\infty}$ quasi-equivalences between some pre-categories where the spaces of morphism are respectively Morse, Floer, cellular, and \v{C}ech complexes.  We will work with integral coefficients throughout this section since no serious difficulties are simplified by working over a field.  First, we introduce the relevant Morse theory.
\subsection{Morse theory on manifolds with boundary} \label{morse-theory}

Let $(Q,\partial Q)$ be a smooth oriented closed Riemannian manifold. It is well known that one can use Morse theory to study the topology of the pair $(Q, \partial Q)$ without any properness condition.  The reference which closest in spirit to the approach which we need is \cite{handron}.  However, the exact results that we'll use do not seem to be written down anywhere, so we give a brief description.

\begin{defin} \label{boundary_cvx}
A function $f\co Q \to \bR$ is {\bf boundary convex} if no singular gradient trajectory between interior points intersects the boundary.
\end{defin}
\begin{rem}
Recall that if $f$ is not Morse, then a singular gradient trajectory is a concatenation of gradient trajectories and paths contained in the critical locus.
\end{rem}

There is an important local model which justifies our terminology for boundary convex functions:
\begin{lem}
A subset $C$ of $\bR^n$ has convex boundary if and only if the restriction of every linear function $f \co \bR^n \to \bR$ to $C$ is boundary convex. \noproof 
\end{lem}

In general, the boundary convexity of a function is a global property.  However,  we will exploit the existence of sufficient local condition:

\begin{lem} \label{positive_boundary_cvx}
Let $f \co Q \to \bR$ be a function whose interior critical points are isolated from the boundary.   Then $f$ is boundary convex if every point $r \in \partial Q$ at which $\grad(f)$ is not transverse to the boundary satisfies one of the following two conditions:
\begin{enumerate}
\item If $r$ is a critical point of $f$ and $\gamma$ is a path starting at $r$ which is not tangent to the boundary then $f|_{\gamma}$ is strictly decreasing in some neighbourhood of $r$.
\item If $r$ is not a critical point of $f$, then (if it is defined) a gradient segment starting at $r$ is contained in $\partial Q$.
\end{enumerate}
\end{lem}
\begin{proof}
We assume, by contradiction, that $f$ satisfies the above property, and that there are interior points $p$ and $q$ which are connected by a singular gradient trajectory $\gamma$ which is tangent to the boundary at some point.  Let $r$ be the last point at which $\gamma$ intersects $\partial Q$.  Since $f$ is non-decreasing along $\gamma$, the two conditions above preclude the existence of such a point, hence the existence of $\gamma$.
\end{proof}

A function satisfying the conditions of the above Lemma will be called {\bf positively boundary convex}.

\subsubsection{The Morse complex}
The material in this section is familiar in the proper case, see for example \cite{schwarz}, \cite{fukaya-morse-htpy}, and \cite{betz-cohen}.  We will restate some definitions in a slightly non-standard way.  Recall that for each non-degenerate critical point $p$ of a smooth function $f$ we have an ascending (stable) manifold
\[ W^{s}(p)\] 
consisting of points whose image under the negative gradient flow converges to $p$ as $t \to+ \infty$, and a descending (unstable) manifold
\[ W^{u}(p) \] consisting of points whose image under the positive gradient flow converges to $p$ as $t \to + \infty$.
\begin{defin}
The {\bf degree} of $p$, denoted $\deg(p)$, is the dimension of $W^s(p)$.
\end{defin}
Even through $Q$ has been assumed to be oriented, the submanifolds $W^{s}(p)$ do not inherit a natural orientation.
\begin{defin}
The { \bf orientation line} at $p$ is the abelian group generated by the two orientations of  $W^{s}(p)$, with the relation that the two opposite orientations satisfy
\[ [\Omega_{p}] + [-\Omega_p] = 0 .\] 
We denote this group by $|o_p|$.
\end{defin}

\begin{defin}
If $f$ is a boundary convex function on a smooth manifold $Q$,  we define $\partial_f^+ Q$ to be the set of points in $\partial Q$ consisting of limit points of gradient trajectories of $f$ with initial points in the interior of $Q$.  The complement of $\partial_f^+ Q$ in $\partial Q$ is denoted $\partial_f^- Q$.
\end{defin}

We find the following definition convenient, even though it may not agree with some conventions which require Morse-Smale functions to be proper.
\begin{defin}
A function with non-degenerate critical points is called {\bf Morse-Smale} if all ascending and descending submanifolds of interior critical points intersect each other transversally, and $\partial^{+}_f Q$ is a compact codimension $0$ submanifold of the boundary.
\end{defin}
Given a boundary convex Morse-Smale function $f$, we can now define the Morse complex
\[ \CM^{i}(f) = \bigoplus_{\substack{p \in Q - \partial Q \\ \deg(p) = i}} |o_p| .\]
To define the differential, we begin by observing that since $Q$ is oriented, an orientation of any ascending manifold induces an orientation of the corresponding descending manifold.  Further, if 
\[ \deg(q) = \deg(p) +1 ,\]
then $W^{u}(p)$ and $W^{s}(q)$ intersect in a $1$-dimensional submanifold  of $Q$,  consisting of finitely many components whose image under $f$ is a fixed open segment in $\bR$. So if we pick orientations on the ascending manifolds at $p$ and $q$, we can compare the induced orientation on
\[ W^u(p) \cap W^{s}(q)\] 
with the natural orientation pulled back from $\bR$ using $f$.  Let $n(q,p)$ denote the difference between the number of components on which the two orientations agrees, and the ones on which they differ.  This number depends on the choice of orientations, but only up to sign.  We can now define the differential
\[ \mu_1 \co \CM^{i}(f) \to \CM^{i+1}(f)\]
by choosing orientations $\Omega_q$ of each ascending manifold and linearly extending the formula
\[ \Omega_p \mapsto \sum_{q} n(q,p) \Omega_q .\]
Note that even though we chose orientations to write a formula, the differential is independent of all choices.  Let us also note that we are counting gradient lines for $-f$ starting at $p$ and ending at $q$.  This causes some unfortunate sign conventions in the next section.

In general, the differential on the Morse complex of a Morse-Smale function on a manifold with boundary may not square to $0$.  However, we have
\begin{lem}
If $f$ is a boundary convex Morse-Smale function, then $\mu_1^2=0$.
\end{lem}
\begin{proof}
The same proof usually given in the compact case works.  The moduli space of gradient trajectories converging at either ends to critical points whose degree differs by $2$ is a manifold whose boundary consists precisely of the terms of $\mu_1^2$.  Convexity prevents $1$-parameter families of gradient lines from escaping to the boundary, so this manifold is compact, and the signed count of these boundary points vanishes.
\end{proof}

\subsubsection{The $A_{\infty}$ structure}
The $A_{\infty}$ structure on the pre-category of Morse functions was discovered by Fukaya, see \cite{fukaya-garp}.  Let $\{f_i \}_{i=0}^{n}$ be a sequence of functions on $Q$, and assume that each pair $f_i - f_j$ for $i<j$ is boundary convex. We will construct higher products by considering moduli spaces of gradient trees.  For concreteness, and to set the stage for Section \ref{cellular}, we review their construction.
\begin{defin}
A {\bf singular metric ribbon tree} $T$ consists of the following data:
\begin{enumerate}
\item one dimensional finite CW complex,
\item a map $l$ from the set of edges to the interval $(0 , + \infty]$ (the length),
\item a cyclic ordering of the edges adjacent to each vertex,
\item a marked univalent vertex called the outgoing vertex, and
\item a labeling of all vertices as either of finite or infinite type.
\end{enumerate}
The above data must satisfy the following conditions
\begin{enumerate}
\item There are no non-trivial cycles,
\item all vertices of valence greater than $2$ are of finite type,
\item all vertices of valence $2$ are of infinite type,
\item an edge has infinite length if and only is it has an endpoint of infinite type.
\end{enumerate}
\end{defin}

Univalent vertices of $T$ and the edges which include them are called {\bf external} with all other vertices and edges called {\bf internal}.  Note that external edges may have either finite or infinite length.

\begin{rem}
Upon choosing the outgoing vertex, the remaining external vertices (or edges) of $T$ acquire a linear ordering.  We will think of them as {\bf incoming} vertices (similarly for incoming edges).  We equip the edges of a ribbon tree with a natural orientation as follows:
\begin{itemize}
\item The orientation on the incoming external edges points towards the interior of the tree.  The outgoing external edge is oriented away from the interior.
\item The orientations on all but one of the edges adjacent to a fixed vertex $v$ point towards $v$.
\end{itemize}
One can easily prove by induction that an orientation satisfying the above properties exists and is unique.  We shall therefore freely speak of the (unique) ``outgoing edge at a vertex'' (the one oriented away from $v$) or of the incoming vertices of a subtree.  In addition, the internal vertex of a tree which lies on the outgoing edge will be a convenient basepoint for many of our argument.  We will call it the {\bf node}.
\end{rem}

\begin{defin}
A singular ribbon tree is {\bf reducible} if it contains a bivalent vertex.  It is otherwise called irreducible, or, more simply, a {\bf ribbon tree}.  

If $T$ is a singular ribbon tree then a subtree $S \subset T$ is called an {\bf irreducible component} if it is maximal among irreducible subtrees.
\end{defin} 

\begin{lem}
Every singular ribbon tree admits a unique decomposition into irreducible components. \noproof
\end{lem}
In particular, one should think of singular ribbon trees as obtained by taking two ribbon trees with infinite external edges and identifying the univalent vertices of these two edges.

If we remove the univalent vertices of a ribbon tree, we obtain an honest metric space.  If a vertex lies on an infinite edge, the labeling of vertices as either of finite or infinite type determines whether the edge is isometric to $[0 , + \infty)$ or $(- \infty, + \infty)$, with $0$ corresponding to the vertex of finite type.

In particular there is a natural notion of (isometric) automorphisms.  In analogy with the Deligne-Mumford compactification, we have
\begin{defin}
A ribbon tree is {\bf stable} if its group of automorphisms as a metric tree is finite.  A singular ribbon tree is stable if so are all its irreducible components.
\end{defin}

\begin{lem}
All ribbon trees which are not stable are isometric to $( - \infty , + \infty)$.  In particular, a singular ribbon tree is stable unless it includes an edge with two bivalent endpoints. \noproof
\end{lem}

\begin{defin}
The $d$-th {\bf Stasheff moduli space}, $\Stasheff_{d}$, is the moduli space of ribbon trees with $d+1$ external edges all of which have infinite type.
\end{defin}
The topology of $\Stasheff_d$ is such that the operation of contracting the length of an internal edge to zero is continuous.  In \cite{stasheff}, Stasheff showed that this moduli space is homeomorphic to a ball (see also Fukaya and Oh in \cite{FO}).
\begin{lem}[Definition-Lemma]
There exists a compactification of $\Stasheff_d$ into a compact manifold with boundary $\overline{\Stasheff}_d$ called the {\bf Stasheff polyhedron} such that
\[\partial \overline{\Stasheff}_d = \bigcup_{d_1+d_2=d+1}  \bigcup_{1 \leq i \leq d} \overline{\Stasheff}_{d_1} \times  \overline{\Stasheff}_{d_2} .\] 
Moreover, $\overline{\Stasheff}_d$ is naturally stratified by smooth manifolds with corners.
\noproof
\end{lem}
The compactification is simply obtained by adding stable singular ribbon trees all of whose external vertices again have infinite type, and the topology is such that a sequence of trees with fixed topology having an edge whose length becomes unbounded converges to a singular tree with a bivalent vertex precisely at that edge.

Given a fixed tree $T \in \Stasheff_{n}$ and a sequence $\vf=(f_0, \ldots, f_n) $ of functions on $Q$, we label the incoming edges with the functions $\{ f_{i} - f_{i-1} \}_{i=1}^{n}$ and the outgoing edge with $f_n - f_0$.  
\begin{lem}
There exists a unique labeling of the edges of $T$  by functions $f_e \co Q \to \bR$, extending the above labeling of the external edges, such that for every vertex $v$, the sum of the functions labeling incoming edges is equal to the function labeling the (unique) outgoing vertex.  Further, every function $f_e$ is of the form $f_i - f_j$ with $j<i$.\end{lem}
\begin{proof}[Proof following \cite{FO}]
From the data of a cyclic orientation at each vertex we can construct a unique (up to isotopy) embedding $\iota$ of $T$ in the unit disc with external vertices lying on the unit circle.  We can label the components of the complement of $\iota(T)$ by functions $f_i$ such that each external edge is labeled by the difference between the functions associated to the two components which it separates.  To get the correct sign, we assign a positive sign to the function which lies to the left (according to the orientation of the edge).  Note that we can extend this labeling to the internal edges as well, and that it satisfies the desired balancing conditions at every vertex.
\end{proof}

The higher algebraic structure of Morse functions is obtained from the topological properties of the space of maps from trees to $Q$.
\begin{defin} \label{def-gradient-trees}
Given a sequence of functions $\vf$, and critical points $p_{i,i+1}$ of $f_{i+1}-f_i$ and $p_{0,d}$ of $f_d-f_0$,  the { \bf moduli space of gradient trees} 
\[ \Stasheff(p_{0,d}; p_{0,1},p_{1,2}, \ldots, p_{d-1,d}) = \Stasheff(p_{0,d};\vp)  \]
is the space of maps $\phi: S \to Q$ with $S$ an arbitrary element of $\Stasheff_{d}$, which satisfy the following conditions:
\begin{itemize}
\item The image of each edge $e$ is a gradient line for $f_e$.
\item The orientation of $e$ is given by the negative gradient of $f_e$.  Further, $\phi$ identifies the length parametrization of each edge $e$ with the gradient vector field of $f_e$. 
\item If $e$ is an external edge, and $f_e = f_i - f_j$, then the image of $\phi$ converges to $p_{i,j}$ along $e$.
\end{itemize}
\end{defin}
Let us observe that this moduli space depends not only on the functions $f_{i}$, but also on the Riemannian metric of $Q$. We will make this dependence explicit as needed.  For simplicity, we will also write
\[ \Stasheff(\vf) = \coprod_{p_{0,d}, \vp} \Stasheff(p_{0,d}; \vp) .\]

Note that the moduli space of gradient trees is not compact as the length of an internal edge can go to infinity.  However, it admits a bordification
\[ \overline{\Stasheff}(\vf)\]
by observing that, if such a length becomes infinite, the gradient tree converges to a union of gradient trees which share common vertices, i.e. the domain tree converges to a singular ribbon tree equipped with a map to $Q$.  In \cite{fukaya-morse-htpy}, Fukaya proves that $\overline{\Stasheff}(\vf)$ is compact whenever $Q$ is compact.  We have
\begin{lem} \label{morse_compactness}
If each pair $f_i - f_j$ for $j < i$ is boundary convex, then $ \overline{\Stasheff}(\vf) $ is compact.
\end{lem}
\begin{proof} Let $\phi_t$ be such a $1$-parameter family converging to $\phi_1$, such that there is a point $x \in T$ for which
\[ \phi_1(x) \in \partial Q.\] 
If $x$ lies in the interior of an edge $e$, then the image of every point on $e$ converges to the boundary by the boundary convexity of the function $f_i - f_j$.  Since the external edges cannot converge to the boundary, it suffices to prove the following claim.
\begin{claim}
If $v$ is a vertex of $T$ and $\phi_1(v) \in \partial Q$, then every edge adjacent to $v$ lies in $\partial Q$
\end{claim}
\begin{proof}[Proof of Claim]
Assume otherwise.  To simplify the argument, we re-orient the outgoing edge at $v$ so that it now points towards $v$, and take the negative of the corresponding function.  We now have
\[ \sum_{v \in e} f_{e} = 0 .\]
However, if $\nu$ is the outwards pointing normal vector to the boundary at $\phi_1(v)$ then
\[ \langle \nu, \grad(f_e) \rangle \leq 0 ,\]
since tracing back the gradient line of $f_e$ pushes us towards the interior of the manifold.  Comparing the two equations, we conclude that these inner products must all vanish, i.e. that all gradient lines are tangent to the boundary at $v$.  By boundary convexity, the corresponding edges are contained in $\partial Q$.
\end{proof}
The same argument which yields the result when the boundary is empty therefore applies in the boundary convex case to complete the proof of Lemma \ref{morse_compactness}.
\end{proof}
To extract algebraic invariants, we will assume genericity of the functions $f_i$. 
\begin{defin} \label{defin:regularity_moduli_spaces_trees}
A sequence of functions $\vf$ is {\bf Morse-Smale} if all functions $\{ f_{i} - f_{j}\}_{i < j}$ are Morse-Smale, and all finite collections of their ascending and descending manifolds are in ``general position''.  Such a collection is {\bf boundary convex} if all  functions $\{ f_{i} - f_{j}\}_{i < j}$ are boundary convex.
\end{defin}
The precise meaning of general position is a requirement on the regularity of all moduli spaces of gradient tree  as explained in Appendix \ref{mod-space-trees}, in particular Definition \ref{defin:regularity_moduli_trees}.  As before, an orientation of the ascending manifolds of the critical points of $f_{i+1} - f_i$ and of $f_{d} - f_0$ induces an orientation of $ \overline{\Stasheff}(\vf)$; the explanation is again relegated to Appendix \ref{coherent-orient}.  In particular, transversality implies that if the degrees of the various critical points satisfy
\begin{equation} \label{degree_shift} 2-d + \sum_{i} \deg(p_{i,i+1}) = \deg(p_{0,d}) ,\end{equation}
then the corresponding moduli space of gradient trees consists of (signed) points, so we obtain integers
\[n(p_{0,d}, p_{d-1,d}, \ldots, p_{0,1})\]
by choosing orientations of the ascending manifolds, and counting the (signed) number of such points.   We define the above integers to be $0$ whenever Equation \eqref{degree_shift} is not satisfied. We can define a higher product
\[ \mu_d \co \CM^{*}(f_{d-1}, f_{d}) \otimes \cdots \otimes \CM^{*}(f_0,f_1) \to \CM^{*}(f_0, f_d) \]
by linearly extending the formula
\begin{equation}\label{morse_product} \Omega_{p_{d-1,d}} \otimes \cdots \otimes \Omega_{p_{0,1}} \mapsto \sum_{p_{0,d}} n(p_{0,d}, p_{d-1,d}, \ldots, p_{0,1}) \Omega_{p_{0,d}},\end{equation}
where each $ \Omega_{p_{i,j}}$ is a choice of orientation of the appropriate ascending manifold.  As in the case of the differential, it is clear from the construction of Appendix \ref{coherent-orient} that even though we have chosen orientations for the ascending manifolds in order to write a formula, the higher product is independent of all such choices.  Modulo signs, Lemma \ref{morse_compactness} immediately implies the next result.
\begin{prop} \label{signed_A_infty}
If $\vf$ is a set of functions which is  Morse-Smale and boundary convex, then the $A_{\infty}$ equation \eqref{A_infty-equation} holds. \noproof
\end{prop}
The proof of this proposition with signs is given in Appendix \ref{coherent-orient}.

We can now see that Morse functions on a manifold with boundary fit within Kontsevich and Soibelman's framework of pre-categories.
\begin{defin}
The {\bf Morse pre-category} of $Q$ is the pre-category $\Morse(Q)$ defined by the data:
\begin{itemize}
\item Objects are functions on $Q$.
\item  Transverse sequences are sequences of functions which are Morse-Smale and boundary convex.
\item Morphisms between a transverse pair $(f,g)$ are given by the Morse complex
\[ \CM^*(f,g) .\]
\item Higher products are given by Equation \eqref{morse_product}.
\end{itemize}
\end{defin}

\subsection{From Morse chains to cellular chains} \label{cellular}

The proof that ordinary and Morse homology agree extends to our setting.

\begin{lem} \label{morse-ordinary-iso}
\[\HM^{*}(Q,f) \cong H^{*}(Q,\partial^{+}_f Q) \]
\end{lem}
\begin{proof}
Consider the union of all descending manifolds
\begin{equation} \label{eq:union_descending_manifold}  W^s(f) = \bigcup_{p} \overline{W}^{s}(p) .\end{equation}
Since $f$ has no critical points outside $W^s(f)$, the image of any point in $Q - W^s(f)$ under the gradient flow intersects $\partial^{+}_{f} Q$.  In particular, we have a retraction
\begin{equation} Q - W^s(f) \to \partial^{+}_{f} Q  \cap \left(Q - W^s(f) \right) .  \end{equation} 
We conclude that the relative cohomology group
\[  H^{*}(Q - W^s(f), \partial^{+}_f Q \cap (Q - W^{s}(f)))  \] 
vanishes, and hence that the inclusion of $W^s(f) $ induces an isomorphism
\begin{equation} H^{*}(Q,\partial^{+}_f Q) \cong  H^{*}(W^s(f), \partial^{+}_f Q \cap W^{s}(f)) .  \end{equation}
It suffices therefore to prove that there is an isomorphism
\begin{equation}  H^{*}(W^s(f), \partial^{+}_f Q \cap W^{s}(f)) \cong \HM^{*}(Q,f) .\end{equation} 
To prove this, we note that the Morse-Smale condition guarantees that the closure of $\overline{W}^{s}(p)$ intersects $\overline{W}^{s}(q)$ only if the degree of $q$ is smaller than that of $p$.  In particular, we consider the filtration of $W^s(f)$ by the subsets $W^{s,k}(f)$ consisting of the closures of ascending manifolds of critical points whose degree is smaller than or equal to $k$.  Every sub-quotient 
\begin{equation*} \left( W^{s,k}(f) \cup  \left( \partial^{+}_f Q \cap W^{s}(f)\right) \right)  / \left(W^{s,k-1}(f) \cup \left( \partial^{+}_f Q \cap W^{s}(f) \right) \right) \end{equation*}
is homeomorphic to a wedge of spheres indexed by the critical points of degree $k$.  The result follows by induction on $k$, and the inductive step is a classical computation first appearing as Corollary 7.3 in Milnor's book \cite{milnor} that, given a cobordism with critical points only in degree $k$ and $k+1$, the matrix of intersection pairings between ascending and descending manifolds defines a differential on this subquotient of the Morse complex whose cohomology agrees with the cohomology of the cobordism relative a certain subset of its boundary.
\end{proof}

We would like to construct an explicit chain level map inducing the above isomorphism on cohomology.  Further, we would like to relate the cup products at the chain level.  This essentially entails the construction of $A_{\infty}$ functors relating the Morse pre-category to a pre-category whose morphisms are chain complexes that compute ordinary cohomology.  There are many possible options for this putative pre-category, such as singular cocycles, pseudo-cycles, or currents.  If we're willing to work over $\bR$, then currents would be the most convenient category to work with.  However, we would like to keep working over $\bZ$ while we still can.  Unfortunately, pseudo-cycles offer the disadvantage that the intersection product is only partially defined (although, see \cite{mcclure} for an algebraic solution to this problem), while it is not clear to the author how to relate the tree product on the Morse complex to the products on singular cochains without using the abstract techniques of acyclic models as used by Guggenheim to relate the singular cup product with the wedge product on forms \cite{gugenheim}.  

We will choose instead to work with the simplicial cochains of a fixed triangulation, although we will use cellular chains of dual triangulations as an intermediary. We will need to introduce a new moduli space of trees.
\begin{defin} \label{defin-shrubs}
The moduli space of {\bf shrubs} $\Shrub_{d}$ is the moduli space of ribbon trees with $d$ incoming edges all of which have finite type, whose outgoing edge has infinite type, and such that the distances from the node (interior vertex of the outcoming edge) to all the incoming vertices are equal.
\end{defin}
In Appendix \ref{mod-space-trees} we construct a compactification $\overline{\Shrub}_{d}$ of $\Shrub_{d}$ and prove the necessary smoothness and compactness results.

As with the ribbon trees which form the Stasheff polyhedra, we can endow the edges of a tree in $\Shrub_{d}$ with a natural orientation.  Further, if $\vf$ is a sequence of $d+1$ functions, then every edge of a tree $T \in \Shrub_{d}$ can be labeled uniquely by a difference $f_{i} -f_{j}$ satisfying the conditions of Definition \ref{def-gradient-trees}.
\begin{defin} \label{defin:gradient_shrubs}
Let $\vf$ be a sequence of $d+1$ functions on $Q$, $\vC$ a collection of $d$ cycles $(C_1, \ldots, C_d)$, and $p$ a critical point of $f_{d} - f_{0}$.  The moduli space of {\bf gradient trees from $\vC$ to $p$}, denoted $\Shrub(p, \vC)$, consists of all maps
\[ \phi: T \to Q \]
with $T \in \Shrub_{d}$ satisfying the following conditions:
\begin{itemize}
\item The image of the $i$th incoming vertex lies in $C_i$, and
\item the outgoing vertex is mapped to $p$, and
\item each edge $e$ is a gradient line of the corresponding function $f_e$ as in Definition  \ref{def-gradient-trees}.
\end{itemize}
\end{defin}
\begin{rem}
The notation $\Shrub(p, \vC)$ is ambiguous since the function $\vf$ whose gradient trees we are considering are elided.  In practice, however, this should cause no confusion.
\end{rem}

In Appendix \ref{mod-space-trees}, we construct a decomposition of $\Shrub_{d}$ into smooth manifolds with corners corresponding to the different topological types of a shrub, yielding a decomposition of the closure $\overline{\Shrub}_{d}$ into compact manifolds with corners.   To construct an analogous decomposition of  $\Shrub(p,\vC)$, we consider evaluation maps
\[ \eva_{i} \co   W^{s}(p) \times \Shrub_{d} \to Q \]
which are defined for each $1 \leq i \leq d$.  Given a point $(x,T)$, its image under $\eva_i$ is simply the image of $x$ under the composition of gradient flows corresponding to the path in $T$ from the node of $T$ to the input labeled by $i$.  It is easy to see that  $\Shrub(p,\vC)$ is simply the inverse image of $C^1 \times \cdots \times C^{d}$ under the map
\[ \eva \co  W^{s}(p) \times \Shrub_{d} \to Q^{d}\]
whose $i$\textsuperscript{th} component is $\eva_i$.

We will now fix a triangulation $\sP$ of $Q$ such that $\partial Q$ is a subcomplex.  In general, since Morse-Smale functions on a disc can be arbitrarily complicated, there is no hope to recover any information about the Morse homology of an arbitrary function by relating it to the cellular complex of $\sP$.  We therefore fix a subdivision $\check{\sP}$ dual to $\sP$, and we shall include a condition on the behaviour of the Morse flow near strata of $\check{\sP}$ as a restriction on the Morse functions we shall study.

In addition, we fix a family $\{ \kappa_{t} \}_{t\in [0,1]}$ of diffeomorphisms of $Q$ such that $\kappa_{0} $ is the identity, and such that $\check{\sP_{t}}$, the subdivision obtained by applying $\kappa_{t}$  to $\check{\sP}$, is in simplicial position with respect to $\check{\sP_{t'}}$ whenever $t < t'$ (see Appendix \ref{singular-cochains-manifolds}).  In particular, whenever $t_i$ is an increasing sequence, the subdivisions $\check{\sP_{t_i}}$ are mutually transverse.

Given a cell $C$ or a sequence of cells $\vC = (C_1, \ldots, C_d)$ of $\check{\sP}$, we will write $C^{t}$ for the corresponding cell of $\check{\sP}_{t}$ or $\vC(\vt)$ for the sequence 
\[  (C^{t_1}_1, \ldots, C^{t_d}_d)\]
whenever $\vt{}= (t_1, \ldots t_d)$ satisfies $0 \leq t_1 < \cdots < t_d$.

We will further write, as in Appendix \ref{singular-cochains-manifolds}, $ \check{\sP}_{t_1} \Cap \check{\sP}_{t_2} $  for the cellular subdivision of $Q$ consisting of the intersection of all cells in $\check{\sP}_{t_1}$ and $\check{\sP}_{t_2}$, and  
\[ \coll_{t_2} \co  \check{\sP}_{t_1} \Cap \check{\sP}_{t_2}  \to \check{\sP}_{t_2}\]
for the corresponding cellular map.

Given two cells $C^{t_1}_{1}$ and $C^{t_2}_{2}$ we can use the family $\kappa_{t}$ to construct a map
\begin{equation} \label{eq:cell_htpy} \left( C^{t_1}_{1} \cap C^{t_2}_{2}  \right) \times [t_1,t_2] \to Q \end{equation}
whose restriction to $s_1 \in [t_1,t_2] $ identifies $C^{t_1}_{1} \cap C^{t_2}_{2} $ with $C^{s_1}_{1} \cap C^{t_2}_{2} $.  This map realizes a homotopy between the identity and the collapse of $ C^{t_1}_{1} \cap C^{t_2}_{2} $ onto its image under $ \coll_{t_2}$, producing a cell of one higher dimension which we write $K(C_{1},C_{2}, t_1,t_2)$.  More generally, if $\vt{}= (t_1, \ldots, t_{d})$ is a sequence of increasing positive real numbers, let $K(\vt)$ denote the subset of points $(s_1, \ldots, s_{d-1})$ in $[t_1, t_d] \times [t_2, t_d] \times \ldots \times [t_{d-1}, t_{d}]$ satisfying the conditions
\[ s_1 \leq s_2 \leq \cdots \leq s_{d-1}. \]
$K(\vt)$ is a convex polytope; in dimension $2$, for example, this polytope is shown in Figure \ref{fig:K-polytope}.
\begin{figure}[h] 
   \centering
 \input{K-polytope.pstex_t}
   \caption{}
   \label{fig:K-polytope}
\end{figure}

Given a sequence of cells  $(C_1, \cdots, C_{d})$ of  $\check{\sP}$, Lemma \ref{lem:parametrise_family_intersections_cells} implies that we have a diffeomorphism
\begin{equation} C_{1}^{t_1}  \cap \cdots \cap C_{d-1}^{t_{d-1}} \cap C_{d}^{t_d} \cong  C_{1}^{s_1}  \cap \cdots \cap C_{d-1}^{s_{d-1}} \cap C_{d}^{t_d} \end{equation}
that varies continuously with $\vs \in K(\vt)$.  In particular, we obtain a map
\begin{equation} \label{eq:higher_cell_htpy}
K(\vC,\vt) \cong \bigcap_{i} C_{i}^{t_i}  \times K(\vt) \to Q . \end{equation}

We now construct a new pre-category structure on  $\Morse(Q)$:
\begin{defin}  \label{transverse-in-Morse}
A pair of functions $(f_0,f_1)$ is transverse in $\Morse(\check{\sP})$ if and only if
\begin{itemize}
\item  the pair $(f_0, f_1)$  is transverse in $\Morse(Q)$, and 
\item the inclusion
\[ \partial^{+}_f Q \subset \bigcup_{\substack{ C \in \check{\sP} \\ C \cap \partial^{+}_f Q \neq \emptyset}} C \equiv \partial^{+}_{f} \check{\sP} \]
is a deformation retract.
\end{itemize}

A sequence of functions $(f_0, \ldots, f_d)$ is transverse in $\Morse(\check{\sP})$ if and only if 
\begin{itemize}
\item it is transverse in $\Morse(Q)$, and
\item every subsequence consisting of two elements is transverse.
\end{itemize}
\end{defin}

The second condition required for the transversality of a pair of functions essentially forces the dynamics of the gradient flow of $f_1 - f_0$ to be no more (and no less) complex than the inclusion of $ \partial^{+}_{f} \check{\sP}$ into $\check{\sP}$.

In Appendix \ref{singular-cochains-manifolds}, we construct a pre-category $\Cell(\check{\sP})$ with objects cell subdivisions which are $C^1$-close to $\check{\sP}$ (together with extra data at the boundary corresponding to the subsets $\partial^{+}_f \check{\sP}$), and where the cup product of simplicial cochains is realized as a geometric intersection pairing.  We write $\Cell_{\kappa}(\check{\sP})$ for the subcategory of $\Cell(\check{\sP})$ consisting of cell subdivisions obtained by moving $\check{\sP}$ along $\kappa_{t}$ for some $t$.  We now consider a new pre-category $\Cell-\Morse(\check{\sP})$ which will serve as an intermediary between $\Cell_{\kappa}(\check{\sP})$ and $\Morse(\check{\sP})$.  Again, we do this in two steps.

\begin{defin}
The objects of $\Cell-\Morse(\check{\sP})$ are pairs $(\check{\sP}_{t_0}, f_0)$, with $f_0$ a function on $Q$. 

Two objects  $(\check{\sP}_{t_0}, f_0)$ and $(\check{\sP}_{t_1}, f_1)$ are transverse if and only if \begin{enumerate}[(i)]
\item $t_1$ is greater than $t_0$, and
\item $(f_0, f_1)$ are transverse in $\Morse(\check{\sP})$, and 
\item the moduli spaces of shrubs with endpoints on a cell of $\check{\sP}_{t_1}$ are regular.
\end{enumerate}

If  $(\check{\sP_0}, f_0)$ and $(\check{\sP_1}, f_1)$  are transverse objects of $\Cell-\Morse(\check{\sP})$, the {\bf space of morphisms} between them is defined to be
\[ C_{n-*}(\check{\sP_1} -\partial^+_{f_1 - f_0} \check{\sP_1}). \]
\end{defin}

Since our ultimate goal is to pass from $\Cell-\Morse(\check{\sP})$ to $\Morse(\check{\sP})$, we shall require transversality between the gradient flow lines of functions $f_i$ and the cells of various subdivisions $\check{\sP}_{t_i}$.  

To state this formally, given a sequence of cells $\vC = (C_1, \ldots, C_{d})$ of $\check{\sP}$, times $\vt{}= \{  t_1 < \ldots < t_d \}$, and a partition $\vE$ of the set $\{1, \ldots ,d \}$ into subset 
\begin{equation} \label{eq:partition_E} \vE = \{ E_1, \ldots, E_{e} \} = \{ \{1, 2, \ldots, i_2 -1 \}, \{i_2, \ldots, i_3-1 \}, \cdots, \{ i_{e}, \ldots, d \} \} , \end{equation}
we write $\vC_{E_k} = (C_{i_{k}}, \ldots, C_{i_{k+1} -1})$, and $\vt_{E_{k}}$ for the corresponding subsequence of $\vt$.  We consider
\[ K_{\vE}(\vC,\vt) = K(\vC_{E_1},\vt_{E_{1}}) \times \cdots \times K(\vC_{E_e},\vt_{E_{e}}) . \] 

Generalising  Definition \ref{defin:gradient_shrubs}, we have a moduli space
\begin{equation} \Shrub_{e}(p,  K_{\vE{}}(\vC{},\vt{}) ) \end{equation}
obtained as a fiber product
\begin{equation} \Shrub_{e}(p,  Q, \ldots, Q) \times_{Q^{e}}  K_{\vE}(\vC,\vt) \end{equation}
where the map from $K_{\vE}(\vC,\vt)$ to $Q^e$ is given component-wise by Equation \eqref{eq:higher_cell_htpy}.

We incorporate a transversality conditions for  these moduli spaces of shrubs with endpoints on $K_{\vE}(\vC,\vt)$ into our definition of the category $\Cell-\Morse(\check{\sP})$.    We shall need this in order to define a functor to $\Morse(\check{\sP})$.
\begin{defin} \label{end-defin-cell-morse}
A sequence $ \{ (\check{\sP}_{t_i}, f_i) \}_{i=0}^{d}$ is transverse in $\Cell-\Morse(\check{\sP})$ if 
\begin{enumerate}[(i)]
\item the sequence $\{ \check{\sP}_{t_i} \}_{i=0}^{d}$ is transverse in $\Cell(\check{\sP})$, and
\item the sequence $\{ f_{i} \}_{i=0}^{d}$ is transverse in $\Morse(\check{\sP})$, and 
\item if $\vC = (C_1, \ldots, C_{d})$ is any sequence of cells in $\check{\sP}$, $\vt{}'$ is a subsequence of $\vt$ of length $d'$ with corresponding subsequence $\vC{}'$ of $\vC$, and $\vE{}'$ a partition of $\{1 , \ldots, d'\}$ as in Equation \eqref{eq:partition_E}, then the moduli spaces
\[ \Shrub_{e'}(p,  K_{\vE{}'}(\vC{}',\vt{}') ) \]
are regular for any critical point $p$ of $f_d - f_0$.
\end{enumerate}
\end{defin}

Note that whenever $\vC{}' = \vC$, and $\vE$ is the trivial partition $\{\{1\}, \{2\}, \ldots, \{d\} \}$, the last conditions is imposing the regularity of the moduli spaces
\[ \Shrub_{d}(\vC, \vt) .\]

Let us now consider the situation where a moduli space
\[ \Shrub_{e}(p,  K_{\vE}(\vC,\vt) ) \]
has virtual dimension $0$.  If the partition $\vE$ contains a subset $E_{k}$ which consists of more than one element, then the evaluation map
\begin{equation*} K_{\vE}(\vC,\vt) \to Q^{d} \end{equation*}
has image which is contained in a neighbourhood of a cell of $\sP^{d}$ whose dimension is strictly smaller than that of $K_{\vE}(\vC,\vt) $.  In particular, if all parameters $t_i$ are sufficiently small, we expect this moduli space to be empty.  Indeed, assuming regularity, the moduli problem for the cell of lower dimension has no solutions, so there should be no solution for $C^0$ close cells, regardless of their dimension.  In this case, we could avoid introducing the higher dimensional cells $K_{\vE}(\vC,\vt)$ and define a more straightforward transversality condition using only the moduli spaces $\Shrub_{d}(p, \vC(\vt))$ which would be essentially independent of sufficiently small choices of $(t_1, \ldots, t_d)$.  We elaborate further on this point in the next section when it is time to prove the $A_{\infty}$ equation for the functor from the cellular to the Morse category.

The problem is that while the moduli spaces of shrubs with inputs on cells of $\sP$ are regular for generic choices of Morse functions, their compactifications do not satisfy this property because cells of $\sP$ do not necessarily intersect transversally.  It seems likely that the moduli spaces $\Shrub_{e}(p,  K_{\vE}(\vC,\vt) )$ are still empty whenever some $|E_{k}| > 0$ for appropriate choices of diffeomorphisms $\kappa_{t}$, but proving this would require more subtle transversality arguments.

\begin{lem} \label{lem:dimension_indep_of_partition}
The virtual dimension of $\Shrub_{e}(p,  K_{\vE}(\vC,\vt) )$ is
\[ (d-1) + \deg(p) - \sum_{i} \codim(C_{i})   .\] 
In particular, it is independent of the subdivision $\vE$.  \noproof
\end{lem}

\subsubsection{The functor defined by shrubs}
In this section, we will explain the  proof of:
\begin{prop} \label{cech-to-morse}
The pre-categories $\Cell(\check{\sP})$ and $\Morse(\check{\sP})$ are equivalent as $A_{\infty}$  pre-categories.
\end{prop}

Before proceeding with the construction, let us make the auxiliary choice of an outward pointing vector field $\nu$ on $\partial Q$ which is parallel to the normal vector of the boundary.  We let
\[ H(f) =  \langle \grad(f), \nu \rangle \]
and observe that $H(f)$ is positive precisely on $\partial^+_f Q$.  The function $H(f)$ will play a r\^ole analogous to that of the function $H$ in Appendix \ref{singular-cochains-manifolds}.  In particular, we have
\[\partial^+_{H(f_1 - f_2)}  \check{\sP_1} = \partial^+_{f_1 - f_2}  \check{\sP_1} ,\]
giving us two (slightly) incompatible notations.

We are now ready to use $\Cell-\Morse(\check{\sP})$ for its stated purpose.  First, observe that the map on objects
\[ (\check{\sP}_{t_0}, f_0) \to (\check{\sP}_{t_0}, H(f_0) ) \]
induces the obvious equivalence of $A_{\infty}$-pre-categories from  $\Cell-\Morse(\check{\sP})$ to  $\Cell_{\kappa}(\check{\sP})$; the subcategory of $\Cell(\check{\sP})$ consisting of subdivisions obtained by applying $\kappa_{t}$ to $\sP$ for some $\kappa_{t}$.  The inclusion
\[ \Cell_{\kappa}(\check{\sP}) \to \Cell(\check{\sP})\]
is clearly essentially surjective.

We now construct an $A_{\infty}$ functor from $\Cell-\Morse(\check{\sP})$ to $\Morse(\check{\sP})$.  The functor is the forgetful map on objects
\[ (\check{\sP}_{t_0}, f_0) \to f_0 .\]
The linear term of the functor is
\begin{align*} \cF_{cm}:  C_{n-k}(\check{\sP}_{t_1} - \partial^+_{f_1 - f_0} \check{\sP}_{t_1}) &\to \CM^{k}(f_0,f_1) \\
C^{t_1} & \mapsto \sum_{p \in \Crit(f_1 - f_0)} m(p,C^{t_1}) \Omega_p,
\end{align*}
where $m(p,C^{t_1})$ is the signed count of elements of $\Shrub(p,C^{t_1})$.  Note that the transversality assumptions in $\Cell-\Morse(\check{\sP})$ ensure that $\Shrub(p,C^{t_1})$ consists only of points, given that $C^{t_1}$ has codimension $k$ and $p$ has degree $k$.  The sign of a gradient trajectory is simply the sign of the intersection of the oriented chain $C^{t_1}$ with the ascending submanifold $W^{s}(p)$ carrying the orientation $\Omega_p$.  To prove that $\cF_{cm}$ is a chain map, we need to understand $1$-dimensional moduli spaces:
\begin{claim}
The union of the manifolds $\Shrub(q,C^{t_1})$ over all critical points $q$ of degree $k+1$ is a $1$-dimensional manifold admitting a compactification whose boundary points correspond to rigid gradient trajectories from $\partial C^{t_1}$ to $q$,  or to broken gradient trajectories from $C^{t_1}$ to $q$ passing through a critical point $p$ of degree $k$.
\end{claim}

In the absence of a boundary on $Q$, this is a standard result about moduli spaces of gradient flow lines.  The only difference for $\partial Q \neq \emptyset$ is that we must prove compactness of moduli spaces separately.  However, transversality in $\Morse(\check{\sP})$ subsumes the boundary convexity of the function $f_1 - f_0$.  In particular, gradient lines cannot escape to the boundary so the desired result holds.  By interpreting the boundary terms algebraically, we conclude that
\[\cF_{cm} \circ \partial = \mu_1 \circ \cF_{cm}.\]

To see that $\cF_{cm}$ is a quasi-isomorphism, we recall that $m(p,C)$ is the signed intersection number between the oriented chain $C^{t_1}$ and the ascending submanifold $W^{s}(p)$.  Consider the complex of chains $C^{\tr}_{*}(Q, \partial^+_{f_1-f_0}Q)$ which are transverse to cells of $\check{\sP}_{t_1}$.  By assigning to each critical point its stable manifold, we obtain a map
\[ \CM_{*}(f_1 - f_0) \to C^{\tr}_{*}(Q, \partial^+_{f_1-f_0}Q) ,\]
which is a quasi-isomorphism by Lemma \ref{morse-ordinary-iso}.  By construction, $\cF_{cm}$ intertwines the duality pairing between Morse homology and cohomology, with the intersection pairing between $ C^{\tr}_{*}(Q, \partial^+_{f_1-f_0}Q)$ and $C_{n-*}(\check{\sP}_{t_1} - \partial^+_{f_1 - f_0} \check{\sP_{t_1}})$.  Since both pairings are known to be perfect pairings on homology, it follows that $\cF_{cm}$ must also be a quasi-isomorphism.

Let us describe the polynomial extensions of the functor $\cF_{cm}$.  Given a collection $ \{ (\check{\sP}_{t^i}, f_i) \}_{i=0}^{d}$  of $d+1$ transverse objects in $\Cell-\Morse(\check{\sP})$, and cycles $\{C^{t_i}_{i} \}_{i=1}^{d} = \vC(\vt)$ belonging to $C_{*}(\check{\sP}_{t_i} - \partial^+_{f_{i+1} - f_i} \check{\sP}_{t_i})$, and a partition $E$ as in Equation \eqref{eq:partition_E},  we consider a map
\begin{align*} 
\cH_{E}: C_{*}(\check{\sP}_{t_d} - \partial^+_{f_{d+1} - f_d} \check{\sP}_{t_d}) \otimes \cdots \otimes  C_{*}(\check{\sP}_{t_1} - \partial^+_{f_{1} - f_0} \check{\sP}_{t_1}) &\to \CM^*(f_0,f_d) \end{align*}
whose action on generators is given by
\begin{align*}
C^{t_d}_{d} \otimes \cdots \otimes C^{t_1}_{1} & \mapsto \sum_{p \in \Crit(f_d-f_0)} m(p, K_{\vE}(\vC,\vt)) \Omega_p,
\end{align*}
where $m(p, K_{\vE}(\vC,\vt))$ is a signed count of the points in the moduli space
\[\Shrub_{e}(p,  K_{\vE}(\vC,\vt) ) .\]

We define
\begin{equation} \label{eq:defin_morse_simpl_functor} \cF^{d}_{cm}= \sum_{E} \cH_{E} .\end{equation}
An easy application of Lemma  \ref{lem:dimension_indep_of_partition} implies that $\cF^{d}_{cm}$ has degree $1-d$ as expected.

We must prove that the $A_{\infty}$ equation for functors
\begin{align} \label{a_oo-morse-cell} & \sum_{\substack{k \\ \l_1 + \ldots + l_k = d} } \mu_{k} \circ ( \cF_{cm}^{l_1} \otimes \cdots \otimes \cF_{cm}^{l_k}) = \\
& \quad \sum_{i \leq d-1} (-1)^{\maltese_{i}}\cF_{cm}^{d} \circ (\mathbb{1}^{d-i-1} \otimes \partial \otimes \mathbb{1}^{i})  +  \sum_{i\leq d-2} (-1)^{\maltese_{i}} \cF_{cm}^{d-1} \circ (\mathbb{1}^{d-i-2} \otimes \cup \otimes \mathbb{1}^{i}) \notag \end{align}
holds, where $\maltese_{i}$ is as in Lemma  \ref{lem_a_infty_fuk} and $\cup$ stands for the cup product.  In order to do this, we study, as usual, $1$-dimensional moduli spaces of trees.   

Before explaining the general case, we discuss why the first equation holds, ignoring all signs.  This will shed light on the point of view that the term $\cH_{E}$ for which $E$ consists only of singletons is the essential part of the functor, and that the others terms are corrections imposed by the necessity to introduce perturbations of the subdivision $\check{\sP}$.  So we need to prove that:
\begin{multline} \label{eq:A_2-equation} \mu_{2} \circ \left( \cH_{\{1\}} \otimes \cH_{\{2\}} \right) + \mu_{1} \circ \cH_{\{1\}, \{2 \} } + \cH_{\{1\}, \{2 \} } \circ \left( \mathbb{1} \otimes \partial \right) + \cH_{\{1\}, \{2 \} } \circ \left( \partial \otimes  \mathbb{1} \right) \\ + \cH_{\{2\}} \circ \cup    =  \mu_{1} \circ \cH_{\{1, 2 \} } + \cH_{\{1 , 2 \} } \circ \left( \mathbb{1} \otimes \partial \right) + \cH_{\{1 , 2 \} } \circ \left(  \partial \otimes  \mathbb{1}\right)  \end{multline}

Given two cells $C_1$ and $C_2$ such that the moduli space
\[ \Shrub_{2}(p,  K_{\{1\},\{2\}}(\vC,\vt)) = \Shrub_{2}(p,C_{1}^{t_1}, C_{2}^{t_2})\]
is $1$-dimensional, we recall that the boundary strata of the compactification correspond to the following possibilities
\begin{enumerate}
\item   The distance between some internal point of the tree and an incoming vertex goes to infinity.  Because of the restriction that the distance function to the incoming vertex is independent of the choice of incoming vertex, this can only happen if either (i) the outgoing edge converges to a broken gradient trajectory, or (ii) both incoming edges break.   These two possibilities correspond to the terms:
\[ \mu_{2} \circ \left( \cH_{\{1\}} \otimes \cH_{\{2\}} \right)  + \mu_{1} \circ \cH_{\{1\}, \{2 \} }  \]
\item  One of the incoming vertices escapes to the boundary of the cell that it lies on.  This accounts for the terms
\[ \cH_{\{1\}, \{2 \} } \circ \left( \mathbb{1} \otimes \partial \right) + \cH_{\{1\}, \{2 \} } \circ \left(  \partial \otimes  \mathbb{1} \right) . \]
\item The length of an incoming edge can approach zero.  By definition, both incoming edges collapse, so we end up with a gradient trajectory from $ C_{1}^{t_1} \cap C_{2}^{t_2}  $ to $p$. \end{enumerate}
This last type of boundary stratum is related, but not in general equal to, the term
\[ \cH_{\{2\}} \circ \cup .\]
Indeed, in Appendix \ref{singular-cochains-manifolds}, we interpret the cup product as the composition of the intersection pairing with the collapse map 
\[ \coll_{\sP_{t_2}} \co \check{\sP}_{t_1} \Cap \check{\sP}_{t_2} \to \check{\sP}_{t_2}.\]
Note that the image of $ C_{1}^{t_1} \cap C_{2}^{t_2}  $  under the collapse map is $C_{1}^{t_2} \cap C_{2}^{t_2} $.  We must therefore account for the difference between the number of gradient trajectories starting at $ C_{1}^{t_1} \cap C_{2}^{t_2} $ and $ C_{1}^{t_2} \cap C_{2}^{t_2} $, which we claim is precisely given by the remaining terms
\[ \mu_{1} \circ \cH_{\{1, 2 \} } + \cH_{\{1 , 2 \} } \circ \left( \mathbb{1} \otimes \partial \right) + \cH_{\{1 , 2 \} } \circ \left(  \partial \otimes  \mathbb{1} \right). \]

Consider the $1$-dimensional moduli space 
\[ \Shrub_{2}(p,  K_{\{1,2\}}(\vC,\vt))=  \Shrub_{1}(p,K(C_{1},C_{2},t_1,t_2)) \]
of gradient trees with inputs lying on the image of the evaluation map from $K(C_{1},C_{2},t_1,t_2)$  and output the critical point $p$.  By considering the boundary of this moduli space, we will find five terms respectively corresponding to the composition $\cH_{\{2\}} \circ \cup $,  the count of gradient trajectories starting at $ C_{1}^{t_1} \cap C_{2}^{t_2} $, and the right hand side of Equation \eqref{eq:A_2-equation}.

The compactification has two types of boundaries:
\begin{enumerate}
\item either the gradient trajectory breaks at a critical point $q$, corresponding to the term
\[ \mu_{1} \circ \cH_{\{1, 2 \} }, \] 
\item or the gradient trajectory escapes to the boundary of $K(C_{1},C_{2},t_1,t_2)$.
\end{enumerate}

By construction, the boundary of $K(C_{1},C_{2},t_1,t_2)$ can be decomposed geometrically as a union
\[ \left(C_1^{t_1} \cap C_{2}^{t_2}\right) \bigcup K(\partial C_1,C_2,t_1,t_2) \bigcup  K(C_1, \partial  C_2,t_1,t_2)\bigcup  \left(C_1^{t_2} \cap C_{2}^{t_2}\right).\]
The count of gradient trajectories starting on these $4$ different strata completes the proof of the $A_{\infty}$ equation for $d=2$.

We now explain the general idea behind the proof of the $A_{\infty}$ relation for functors which appears in Equation \eqref{a_oo-morse-cell}.  We again consider $\vC$ and $p$ such that
\[  \Shrub_{e}(p, K_{\vE}(\vC, \vt)) \]
has dimension $1$. Let us recall the reasons why this moduli space is not compact:
\begin{enumerate}
\item   The length of some internal edges can go to infinity.  In $1$-parameter families, the only possibility is that every path from an incoming vertex to the outgoing vertex contains precisely one edge whose length becomes infinite.  This corresponds to the left hand side of Equation  \eqref{a_oo-morse-cell}:
\begin{align}  \sum_{\substack{k \\ \l_1 + \ldots + l_k = d} } \mu_{k} \circ ( \cF_{cm}^{l_1} \otimes \cdots \otimes \cF_{cm}^{l_k}) \end{align}
\item The length of an incoming edge can approach zero.  As we are assuming transversality, this occurs for a pair of incoming edges (in particular, the collapse of triples of incoming edges does not occur in one-parameter families, and hence does not affect the differential).  This corresponds to the count of elements in the moduli space
\begin{multline} \label{eq:naive_intersect} \Shrub_{e-1}{\Big{(}} p, K_{E_1}(\vC_{E_1}, \vt_{E_1}), \ldots, K_{E_{i-1}}(\vC_{E_{i-1}}, \vt_{E_{i-1}}) \cap K_{E_{i}}(\vC_{E_{i}}, \vt_{E_{i}}), \ldots   \\
  \ldots, K_{E_d}(\vC_{E_d}, \vt_{E_d}){\Big{)}} \end{multline}
for some fixed $i$.
\item  One of the incoming vertices escapes to the boundary of the cell $K_{E_{j}}(\vC_{E_{j}}, \vt_{E_{j}})$ that it lies on.   This cell has three types of boundary strata.  To simplify our notation, we consider the case where the partition consists of only one set $E$, allowing us to drop the subscript $j$ from the discussion below:
\begin{enumerate}
\item We take the boundary of a cell in $\vC_{E}$.  This accounts for the first term in the right hand side of Equation  \eqref{a_oo-morse-cell}:
\begin{align} \sum_{i \leq d-1} (-1)^{\maltese_{i}}\cF_{cm}^{d} \circ (\mathbb{1}^{d-i-1} \otimes \partial \otimes \mathbb{1}^{i})\end{align}
\item we can reach the boundary of a face in $K( \vt_{E})$ where some $s_i = s_{i+1}$, with the convention that $s_{d} = t_{d}$.  Note that the corresponding cell in $Q$ is the same as the one obtained by applying $K$ to the sequence of cells
\[ (C^{t_1}_{1}, \ldots, C^{t^{i+1}}_{i} \cap C^{t^{i+1}}_{i+1}, \ldots, C^{t^d}_{d}) \] 
whenever $C^{t^{i}}_{i}$ and $C^{t^{i+1}}_{i+1}$ intersect non-trivially.  In this case, Lemma \ref{model-compatible-cup-products} implies that the intersection $C^{t^{i+1}}_{i} \cap C^{t^{i+1}}_{i+1}$ is the result of applying the cup product to $C^{t^{i}}_{i}$ and $C^{t^{i+1}}_{i+1}$, so this boundary stratum corresponds to the remaining term of Equation  \eqref{a_oo-morse-cell}:
\begin{align}  \sum_{i\leq d-2} (-1)^{\maltese_{i}} \cF_{cm}^{d-1} \circ (\mathbb{1}^{d-i-2} \otimes \cup \otimes \mathbb{1}^{i}) \end{align}
\item the remaining boundary stratum correspond to $s_{i} = t_{i}$ for some $i$.  One can check that these boundary strata exactly cancel out with the strata coming from Equation \eqref{eq:naive_intersect}.
\end{enumerate}
\end{enumerate}

The $A_{\infty}$ equation will then hold for the usual argument that the algebraic count of the number of boundary points of a $1$-dimensional manifold vanishes.  To complete the argument, one should check the correctness of the signs.  We omit the long, but straightforward, computation which can be performed along the lines of the arguments given in Appendix \ref{coherent-orient} to prove the correctness of the signs in the $A_{\infty}$-structure of  $\Morse(\check{\sP})$.

Combining this result with Proposition  \ref{cell-sing}, we conclude
\begin{cor} \label{cell-morse}
The pre-categories $\Rel(\sP)$ and $\Morse(\check{\sP})$ are equivalent as $A_{\infty}$  pre-categories. \noproof
 \end{cor}

\subsection{The equivalence of Fukaya and Morse pre-categories} \label{floer-morse}

Given a collection of functions $(f_0, \ldots, f_n)$ on a smooth compact Riemannian manifold we obtain a collection of Lagrangians $L_i$, the graphs of $df_i$, in $T^{*} M$, and the critical points $x_i \in \Crit(f_i - f_j)$ correspond to the intersection points of $L_i$ with $L_j$.  This yields a canonical isomorphism between the underlying graded vector spaces of Floer and Morse complexes.   In \cite{FO}, Fukaya and Oh extended Floer's results that the differentials on Floer and Morse theory agree.  These results have been reworked and amplified in \cites{ekholm,ruan}. Although the results are usually stated for manifolds without boundary, one can easily extend them to the case where the boundary is not empty.  In our notation below, $\epsilon g$ stands for the rescaled metric, and $\epsilon J$ for the induced almost complex structure on $T^*Q$.
\begin{prop}  \label{morse-floer-moduli-equal}
Let $(Q,g)$ be a compact Riemannian manifold possibly with boundary.  For every generic finite collection of functions $(f_0, \cdots, f_n)$ which are transverse in $\Morse(Q)$ and which have no boundary critical points, there exists a constant $\epsilon_0 > 0$ such that for all $\epsilon < \epsilon_0$ and for every collection
\[ \{x_{0,n},x_{0,1}, \ldots, x_{n-1,n} \} \]
 of critical points $x_{i,j}$ of $f_j -f_i$, there exists an orientation preserving homeomorphism between moduli spaces of holomorphic discs and gradient trees:
\[ \cM_{ \epsilon J}(x_{0,n} ; x_{0,1}, \ldots, x_{n-1,n}) \cong   \Stasheff_{\epsilon g}(x_{0,n} ; x_{0,1}, \ldots, x_{n-1,n}).\]
Further, every holomorphic curves is $C^0$ close to the corresponding gradient tree
\end{prop}
\begin{proof}
Extend $Q$ to a smooth compact manifold without boundary $\hat{Q}$, and extend the function $f_i$ to smooth functions on $\hat{Q}$.  By Fukaya-Oh's result, there exists an $\epsilon$, small enough, so that moduli spaces of holomorphic discs and gradient trees on $\hat{Q}$ are in bijective correspondence, and all holomorphic discs are $C^0$ close to gradient trees.  We write $\sN_{\delta} (K)$ for the $\delta$-neighbourhood of a set $K$ in $\hat{Q}$.  Note that the assumption that $(f_0, \cdots, f_n)$  are transverse in $\Morse(Q)$ implies that, for sufficiently small $\delta$,  all gradient trees contained in $ Q$ are contained in $Q - \sN_{\delta} (\partial Q)$.  Choosing $\epsilon$ small enough therefore guarantees that all holomorphic discs which project to  $\sN_{\delta} (Q)$ project to $Q - \sN_{\delta} \partial Q$.  We conclude that for $\epsilon$ small enough, the theorem of Fukaya-Oh holds for the manifold with boundary $Q$.
\end{proof}
\begin{rem}
If the sequence $\vf$ is allowed to have critical points at the boundary, then boundary convexity implies that there is some neighbourhood of the boundary in which gradient trees of interior critical points do no enter.  For sufficiently small $\epsilon$, the argument above implies that count of holomorphic discs in the complement of this neighbourhood agrees with the count of gradient trees.
\end{rem}

Fukaya and Oh also state
\begin{cor}  \label{morse-floer-equivalent}
The Fukaya pre-category of exact Lagrangian sections of $T^{*}Q$ is $A_{\infty}$ quasi-isomorphic to the Morse pre-category of $Q$.
\end{cor}
We will write  $\Fuk(Q)$ for this category of exact Lagrangian sections.  In order for this to make sense when $Q$ is a manifold with boundary, we need to make sure that we have compactness for moduli spaces of holomorphic curves.  In the application we have in mind, this follows from an extra structure imposed at infinity (admissibility with respect to the ``superpotential").  We explain the (completely formal) proof of the above Corollary in the absence of a boundary.

\begin{proof}[Proof of Corollary \ref{morse-floer-equivalent}]
Let $\epsilon_i \to 0$ be a sequence of rescaling parameters, and let
\[ \Fuk(Q; \epsilon_i J) \] 
be the corresponding Fukaya categories.  The homotopy method (see \cites{FOOO,seidel-book}) defines $A_{\infty}$ equivalences of pre-categories
\[ \Phi_{i,i+1} \co \Fuk(Q; \epsilon_i J) \to \Fuk(Q;\epsilon_{i+1}J) ,\]
which is essentially the identity on objects.  In fact, because of transversality problems, this functor is only defined on a full pre-subcategory of $\Fuk(Q; \epsilon_i J)$ consisting of objects for which the family of almost complex structure $\epsilon J$ satisfies an appropriate transversality condition.  For simplicity of notation, we will neglect this problem.

For each $i$, we consider the sub-pre-category $ \Fuk_i(Q; \epsilon_i J) $ where a collection $( L_0, \cdots, L_n )$ of Lagrangians is transverse if and only if the moduli spaces
\[ \cM_{\epsilon J}(x_{0,n} ; x_{0,1}, \ldots, x_{n-1,n}) \cong   \Stasheff_{\epsilon g}(x_{0,n} ; x_{0,1}, \ldots, x_{n-1,n}).\]
are diffeomorphic for all intersection points $x_{i,j} \in L_i \cap L_j$ and for all $\epsilon < \epsilon_i$.  The functors $\Phi_{i,i+1}$ preserve these subcategories in the sense that if a sequence is transverse in $ \Fuk_i(Q;J_{\epsilon_i}) $, its image under $\Phi_{i,i+1} $ is transverse in $ \Fuk_{i+1}(Q; \epsilon_{i+1}J)$.  In particular, we obtain a directed system (in fact, a nested embedding)
\begin{equation} \label{FO-equation}  \Fuk_0(Q; \epsilon_0 J) \to \Phi_{0,1}^{-1}( \Fuk_1(Q; \epsilon_1 J)) \to \Phi_{0,2}^{-1}( \Fuk_2(Q; \epsilon_2  J))  \to \cdots  \end{equation}
whose limit is $\Fuk(Q;J_{\epsilon_0})$.  We can think of this as a filtration of  $\Fuk(Q;J_{\epsilon_0})$.

Denote by $j_i$ the ``identity functor"
\[ \Fuk_{i}(Q;\epsilon_{i+1}J)  \to \Morse_{g}(Q)\]
taking every exact Lagrangian section to its defining function.  By construction, there is an equality of (restricted) functors
\[ j_0 = j_{1} \circ \Phi_{0,1} \]
on $ \Fuk_{0}(Q; \epsilon_0 J) $.  This implies that we can define the right hand side as the extension of $j_0$ to the subcategory $\Phi^{-1}_{0,1}(\Fuk_{1}(Q; \epsilon_1 J))$. Proceeding inductively, we can extend $j_0$ through every step of the direct limit of Equation \eqref{FO-equation}.  Passing to the limit, we obtained the desired $A_{\infty}$ equivalence
\[ \Fuk_{0}(Q; \epsilon_0 J)   \to  \Morse_{g}(Q)  \]

\end{proof}

\section{HMS for toric varieties: Fukaya = Morse} \label{fukaya=morse-section}

Having established the general framework for relating the multiplicative structures on various chain models of ordinary homology, we return to our goal of proving our main result on homological mirror symmetry for toric varieties.  

\subsection{Tropical Lagrangians as a Morse pre-category} \label{tropical-morse}

The observation of Section \ref{tropical-lagr} that all tropical Lagrangian sections lift to graphs of exact $1$-forms allows us to define a Morse pre-category. 

\begin{defin}
The {\bf Morse pre-category of $\sQ$} is the pre-category $\Morse^{\bZ^{n}}(\sQ)$ given by the following data
\begin{itemize}
\item Objects are functions on $\sQ$ the graphs of whose differentials agree with lifts to $T^{*} \bR^{n}$ of tropical Lagrangian sections in $(\ctorus{n}, M)$.
\item Transverse sequences are given by  functions $(f_0, f_1, \ldots, f_d)$ which are Morse-Smale and define lifts of a positive sequence of Lagrangians $(L_0, \ldots, L_d)$.
\item Morphisms between $f_0$ and $f_1$ on $\sQ$ are given by
\[ \CM_{\bZ^{n}}^{*}(f_0,f_1) = \bigoplus_{u \in \bZ^{n}} \CM^{*}(f_0, f_1+u) .\]
\item Higher products are defined by Equation \eqref{morse_product}.
\end{itemize}
\end{defin}
\begin{rems}
As in the previous section, the graded abelian group $\CM^{*}(f_0, f_1+u)$ is generated by interior critical points. Further, when we write $f_1 +  u$, we're thinking of $u$ as a linear function on $\bR^n$
\[ y \mapsto \langle y, u \rangle .\]
\end{rems}

We have chosen this pre-category for its convenience.  Its morphisms are $\bZ^{n}$ graded, and the higher products respect this grading.  If we further required quasi-isomorphisms to live in the $0$ grading, then we would obtain a pre-category which is equivalent to the category of equivariant $\ctorus{n}$ line bundles on the mirror toric variety.  However, we will not impose such a restriction, so that different lifts of the same function lie in the same quasi-isomorphism class.  We could have defined the objects to be tropical Lagrangian sections, but our spaces of morphisms would have only been graded by an affine space over $\bZ^n$.

\begin{lem}
 $\Morse^{\bZ^{n}}(\sQ)$ is a unital $A_{\infty}$ pre-category.
\end{lem}
\begin{proof}
We must prove that the $A_{\infty}$ equation holds for transverse sequences. By the results of Section \ref{morse-theory}, this follows from the
\begin{claim} The lifts of a sequence of admissible Lagrangians satisfy the boundary convexity property of Lemma \ref{morse_compactness}.
\end{claim}
\begin{proof}[Proof of Claim]
It suffices to prove the results in a neighbourhood of each cell of $\partial \sQ$.  We first consider the cells labeled by vertices of $\partial Q$.  On these open sets, the  gradient vector fields of functions defining tropical Lagrangian sections are locally constant when restricted to the part of the boundary which lies in such neighbourhoods.  Since $\sQ$ is a convex subset of $\bR^{n}$, the result follows near the $0$-skeleton.  We now proceed inductively by using, near each $i$ cell, the product decomposition which was defined in \cite{abouzaid}*{Lemma 4.4}.
\end{proof}

It remains to prove unitality, i.e. to produce a sufficient number of quasi-isomorphic copies of each object.  We observed in \cite{abouzaid}*{Lemma 5.8} that if $r_{\partial \sQ}$ is any function which agrees with the square of the distance to the boundary in some neighbourhood of $\partial \sQ$  then, for any real number $K$, the image of the graph of $-K dr_{ \partial \sQ}$ under the map $W: \ctorus{n} \to \bC$ is well approximated, near the origin, by a straight line of slope $K$.  More generally, if $f$ is a function whose differential defines an admissible tropical Lagrangian section, then for every real number $K$,
\[ f - K r_{ \partial \sQ} \]
is $C^2$ close in some neighbourhood of $\partial \sQ$ to a function  $f'_K$ defining another admissible tropical Lagrangian section such that the pair $(f, f'_K)$ is positive (See \cite{abouzaid}*{Lemma 4.15}).  More precisely, if the curve corresponding to $f$ has slope $a$ at the origin, then the curve corresponding to $\gamma'_k$ has slope $a - K$,

Since
\[ \CM^*(f, f'_K) \cong \CM^*(0 ,r) \cong H^{*}(\bR^{n}), \]
we have a canonical element in $e \in CM^{0}(f,f'_K)$ which represents the generator of $H^0(\bR^{n})$.  By Proposition \ref{cech-to-morse}, this generator acts as a cohomological unit.  In particular, if $(f_0, \ldots, f_n)$ is a finite transverse sequence of admissible functions, there exists a constant $K \gg 0$ such that the sequence
\[ (f'_{-K}, f_0, \ldots, f_n, f'_{K}) \]
defines a positive sequence of admissible Lagrangians.  The constant $K$ can simply be chosen to be much larger than the absolute value of the slopes at the origin of all the curves corresponding to the admissible Lagrangians defines by $f_i$ and $f$. We can then choose a $C^{\infty}$-small perturbation of $r$ away from the boundary of $\sQ$ to achieve transversality. 
\end{proof}
\begin{rem}
We were careful to state the results of Section  \ref{morse-theory}, allowing critical points on the boundary in order to be able to apply them directly here.  There is a slightly subtle point concerning the positivity condition.  Consider the case where we are computing the space of morphisms between the zero-section and a Lagrangian defined by a function $f$.  Whenever $f$ has a boundary critical point, the set-up of Section \ref{morse-theory} allows for this point to be included in $\partial^{+}_{f} \sQ$ or $\partial^{-}_{f}  \sQ$  depending on the behaviour of the gradient flow nearby.  However, the positivity assumption precisely eliminates this ambiguity.  Indeed, the gradient vector field of $f$ points inward as soon as we move away from the boundary towards the interior.  This implies that $\partial^{-}_{f} \sQ$ is exactly the subset of $\partial \sQ$ where the gradient flow of $f$ satisfies the (non-strict) inequality
\[  \langle \grad(f), \nu \rangle \leq 0. \]
\end{rem}

\subsection{Equivalence of Morse and Fukaya categories} \label{morse=fuk-tropical}

In order to prove the equivalence of Morse and Fukaya categories, we need to establish more control over the tropical construction.  While both Floer and Morse trajectories do not escape to the boundary, we cannot directly apply the Fukaya-Oh argument for two reasons:
\begin{enumerate}
\item The almost complex structure we're using near $\partial T^{*} \sQ$ does not agree with the complex structure $J$ induced by the flat metric on $\sQ$ (which is the restriction of the complex structure on $\ctorus{n}$), and
\item admissibility relies on the existence of a holomorphic function near the boundary.  Such a function will not remain holomorphic upon ``rescaling" $J$ as in the Fukaya-Oh procedure.
\end{enumerate}

We return to the notation where we use a  semi-tropical degeneration producing a family of functions $W_{t,s}$,  which define symplectic hypersurfaces $M_{t,s} \subset \ctorus{n}$, and hence a family of subsets $\sQ_{t,s} \subset \bR^{n}$ and almost complex structures $J_{t,s}$ on $\ctorus{n}$.  We will also be using the notation of Section \ref{floer-morse} so we write $\epsilon J$ for the almost complex structure induced on $\ctorus{n} = T^* \bR^{n} / \bZ^n$ by rescaling the metric on the base $\bR^n$ by $\epsilon$.  We will write $\cL^{t}$ for the part of the zero section of $\ctorus{n}$ lying over $\sQ_{t,1}$, and $\tilde{\cL}^{t}$ for its lift to the cotangent bundle.

In Appendix \ref{ap:second_look} we explain how to modify the construction given in \cite{abouzaid} in order to prove the following result.  We simply state and prove the case with one Lagrangian, the techniques we use readily extend to establish compactness theorems for finitely many Lagrangians.

\begin{prop} \label{morse_floer_compactness}
Let $L$ be an admissible tropical Lagrangian section with boundary on $M_{t_0,1}$ such that $( \cL^{t_0}, L)$ is a positive pair, and let $f$ be a function defining a lift $\tilde{L}$ of $L$ to $T^{*} \bR^{n}$.  There exists a family $f^t$ for $t \in [t_0, + \infty)$ of functions defining admissible tropical Lagrangian sections $L^t$ with boundary on $M_{t,1}$  such that $(\cL^{t}, L^t)$ is positive, all interior intersection points between $L^t$ and the zero section are transverse, and for sufficiently large $t$ there are positive constants $\delta^t_1$,  $\delta^t_2$ and $\delta^t_3$ such that:
\begin{enumerate}
\item for any $u \in \bZ^{n}$, the Morse gradient flows of $f -u$ connecting interior critical points lie in the complement of the $\delta^{t}_1$ neighbourhood of $\partial \sQ_{t,1}$, and 
\item all $J_{t,s}$ holomorphic discs with marked points on interior intersection points of $L^t$ with the zero section project to the complement of the $\delta^{t}_2$ neighbourhood of $\partial \sQ_{t,1}$, and
\item families of $\epsilon J$ holomorphic discs with marked points on interior intersection points of $L^t$ and $\cL^t$ which for some $\epsilon_0$ lie outside the $2\delta^{t}_3$ neighbourhood of $\partial \sQ_{t,1}$ remain outside the $\delta^{t}_3$ neighbourhood of $\partial \sQ_{t,1}$ for all $\epsilon$, and
\item $2\delta^t_{3}$ is smaller than $\delta^{t}_1$ and $\delta^{t}_2$.
\end{enumerate}
\end{prop}

Assume that we have chosen a family $L^{t}$ and $t$ large enough to satisfy the above conditions.  We will want to produce a chain isomorphism
\[  \CF^{*}(\tilde{\cL}^{t_0}, \tilde{L}; J_{t_0,1}) \to   \CM^{*}(0,f) .\]
While it is possible to do this in one step it will be conceptually easier for us to construct a sequence of chain isomorphisms
\begin{multline*} \CF^{*}(\tilde{\cL}^{t_0}, \tilde{L}; J_{t_0,1}) \to \CF^{*}(\tilde{\cL}^{t}, \tilde{L}^{t}; J_{t,1}) \to \CF^{*}(\tilde{\cL}^{t}, \tilde{L}^{t}; J_{t,0}) \\
\to \CM^{*}(0,f_{t}) \to  \CM^{*}(0,f). \end{multline*}
All these maps can be constructed in essentially the same way via either continuation maps or the homotopy method, so we will focus on constructing a chain isomorphism
\[ \sF_{\bfJ} \co  \CF^{*}(\tilde{\cL}^{t}, \tilde{L}^{t}; J_{t,1}) \to \CF^{*}(\tilde{\cL}^{t}, \tilde{L}^{t}; J_{t,0}). \] 
Part (2) of Proposition \ref{morse_floer_compactness}  ensures that the right hand side is a complex; it states that holomorphic discs with respect to $J_{t,0}$ stay away from the boundary so we can apply the usual compactness arguments.  We now introduce the relevant moduli spaces for the definition of $\sF_{\bfJ}$.  The technique we are describing is originally due to Floer, although it has been extended by Fukaya, Oh, Ohta, and Ono in their work on Lagrangian Floer homology, \cite{FOOO}

If $p$ and $q$ are intersection points of $ \tilde{L}^{t}$ and $ \tilde{\cL}^{t}$, we write $\cM_{\bfJ}(p,q)$  for the union of all moduli spaces $\cM_{J_{t,s}}(p,q)$ for fixed $t$ and varying $s \in [0,1]$.  This space is naturally topologised as the vanishing locus of a section of an infinite dimensional vector bundle over the product of the interval $[0,1]$ with the space of smooth maps from the disc to $T^{*} \bR^{n}$ satisfying appropriate boundary conditions.   As in the usual Floer theoretic setup, this allows us to define regularity as the requirement that the linearised operator is surjective at the zeros; i.e. that the section intersect the zero level set transversely.  Of course, to make this discussion formal, one introduces appropriate Sobolev spaces.  We refer the interested reader to the references \cites{FOOO,seidel-book}, and simply note that regardless of the details, it is easy to see that,
\[ \virdim( \cM_{\bfJ}(p,q)) = \virdim( \cM_{J_{t,s}}(p,q)) +1 .\]

We note that the condition that $\cM_{\bfJ}(p,q)$ is regular does not imply that each moduli space $\cM_{J_{t,s}}(p,q)$ is regular.  Indeed, we will consider moduli spaces $\cM_{\bfJ}(p,q)$ whose virtual dimension is $0$.  Assuming regularity, this space will therefore consist of finitely many $J_{t,s_i}$ holomorphic discs which occur for some $s_i \in (0,1)$.  Since the virtual dimension of $\cM_{J_{t,s}}(p,q)$ is $-1$ for any $s$, this moduli space would be empty if it were regular.  We will often refer to these discs (which are only regular in families) as {\bf exceptional discs}.

Using the usual transversality arguments, we see that for a generic choice of $L^{t}$:
\begin{enumerate}
\item all moduli spaces $\cM_{J_{t,0}}(p,q)$ and $\cM_{J_{t,1}}(p,q)$ are regular, and
\item the parametrised moduli space $\cM_{\bfJ}(p,q)$ is regular.
\end{enumerate}

\begin{figure}     
\input{chain_homotopy.pstex_t}
\caption{} \label{fig:homotopy}
\end{figure}

In this case, we define a new moduli space $\sH(p,q)$ consisting of sequences of $J_{t,s}$-holomorphic discs as in Figure \ref{fig:homotopy}. To be more precise:
\begin{defin} \label{moduli_space_homotopy}
An element of $\sH(p,q)$ consists of the following data:
\begin{enumerate}
\item A collection $\{ p = p_0, p_{1}, \ldots, p_{k} = q \}$ of interior intersection points between $\tilde{L}^{t}$ and the zero section, together with an increasing sequence  $0 < s_1 < \cdots < s_{k} < 1$, and
\item for each $0 < i \leq k$,  a $J_{t,1-s_i}$ holomorphic disc $u_i$ with boundary on $\tilde{L}^{t} \cup \tilde{\cL}^{t}$  with the incoming point mapping to $p_{i-1}$ and outgoing marked point mapping to $p_{i}$.
\end{enumerate}
\end{defin}

We note that the dimension of $\sH(p,q)$ is given by the difference between the degrees of $p$ and $q$.  In particular, $\sH(p,q)$ is discrete precisely whenever $\deg(p) = \deg(q)$ so we can define a degree $0$ map
\begin{align*}
\sF_{\bfJ}  & \co  CF^{*}(\tilde{\cL}^{t}, \tilde{L}^{t}; J_{t,1}) \to CF^{*}(\tilde{\cL}^{t}, \tilde{L}^{t}; J_{t,0}) \\
p & \mapsto p + \sum_{\deg (p) = \deg (q)} |\sH(p,q)| q
\end{align*}
where $| \sH(p,q)|$ is the signed count of the number of points in $\sH(p,q)$.  The following Lemma essentially goes back to Floer:
\begin{lem}
The map $\sF_{\bfJ}$ is a chain isomorphism. \noproof
\end{lem}

We briefly discuss the proof.  First, we recall that the exactness condition implies that we have an energy filtration on Floer cochain groups: in our case the energy of an intersection point $p$ is the value of $f$ at $p$.  By ordering all intersection points according to their action, we obtain a basis for the Floer cochain groups. The usual argument involving Stokes' theorem implies that the matrix representing $\sF_{\bfJ}$ in this basis is upper triangular with $1$ along the diagonal.  In particular $\sF_{\bfJ}$ will necessarily be a chain isomorphism if it happens to be a chain map.

To prove that $\sF_{\bfJ}$ is a chain map, one considers the compactification of the $1$-dimensional moduli spaces $\sH(p,q)$; this corresponds to the condition that $\deg(q) = \deg(p) + 1$.  One can check that the only sequences $\{ p = p_0, p_{1}, \ldots, p_{k} = q \}$ which contribute to $\sH(p,q)$ are those for which there is $i < k$ such that $\deg(p_j) = \deg(p)$ for all $j \leq i$ and all other intersection points have degree equal to $\deg(q)$.  In this case, we have a collection of exceptional solutions together with one (regular) $J_{t,1-s_i}$-holomorphic strip connecting $p_{i}$ and $p_{i+1}$, and it is clear that we can allow $s_i$ to vary, giving us precisely a $1$-dimensional moduli space.

There is a natural bordification $\overline{\sH}(p,q)$ of $\sH(p,q)$ which allows broken curves, as well as the possibility that $s_0 = 0$ or $s_1 =1$.    The proof that $\sF_{\bfJ}$ is a chain map has now been reduced to an analysis of the boundary of the bordification $\overline{\sH}(p,q)$.  Two of the  boundary strata correspond to having $s_0 = 0$, and $\deg(p_1) = \deg(p) +1$ or, alternatively, $s_1 = 1$ and $\deg(p_{k-1}) + 1 = \deg(q)$.  It is rather clear that these terms correspond to
\[ \mu_{1} \circ \sF_{\bfJ} + \sF_{\bfJ} \circ \mu_1 .\]
We must therefore argue that $\overline{\sH}(p,q)$  is compact, and that the contributions of all other boundary strata vanish.  Using part (2) of Proposition \ref{morse_floer_compactness}, we see that $J_{t,s}$ holomorphic curves do not approach the boundary, so we can apply the usual compactness arguments to $\overline{\sH}(p,q)$.  As to the other boundary strata, we observe that they correspond to a broken curve consisting of two components, one exceptional and the other regular, occurring at some fixed time $s_i$.  Note that if $u_{i}$ is a regular curve at time $s_i$, and $u_{i+1}$ is an exceptional curve at time $s_{i+1}$, we obtain exactly such a broken curve if we consider a family where $s_i$ converges to $s_{i+1}$.  However, the usual gluing procedure yields another a family of curves converging to this broken curve.  Taken together, these contributions cancel, proving that $\sF_{\bfJ}$ is a chain map.

The homotopy method was extended by Fukaya, Oh, Ohta, and Ono in their joint work \cite{FOOO} to produce $A_{\infty}$-functors; an account in the exact setting appears in \cite{seidel-book}*{Section 10e}.  The end result is that, using the analogue of Proposition \ref{morse_floer_compactness} in the presence of more than one Lagrangian, we can conclude that the pre-categories 
\[  \sT \Fuk(\ctorus{n},M_{t,1}; J_{t,1}), \, \,   \sT \Fuk(\ctorus{n},M_{t,1}; J_{t,0}) . \]
are $A_{\infty}$ equivalent.

Continuing along the same train of thought, we consider the family of almost complex structures $\epsilon J$.  Part (3) of Proposition \ref{morse_floer_compactness} implies a compactness theorem for families of $\epsilon J$-holomorphic curves which for any $\epsilon_0$ project to the complement of the $\delta^{t}_3$ neighbourhood of $\partial \sQ_{t,1}$.  The assumption that $2\delta^t_{3} < \delta^{t}_2$ implies an equivalence between the chain-level multiplicative structures defined using (i) $J$-holomorphic curves, and (ii) those $\epsilon J$ holomorphic cures whose projecting to $\sQ_{t,1}$ is contained in the complement of the $\delta^{t}_3$ neighbourhood of $\partial \sQ_{t,1}$.

Further, by Fukaya and Oh's theorem, the moduli spaces of $\epsilon J$ holomorphic discs correspond, for $\epsilon$ small enough, to gradient trajectories of $f^t$.  Since $\delta^{t}_{3} < \delta^t_{1}$, it is sufficient to count holomorphic curves which project to the complement of the  $\delta^{t}_3$  neighbourhood of $\partial \sQ_{t,1}$.   Running time backwards (and applying the homotopy method one last time) we can relate the count of gradient trajectories for $f^t$ to the one for the original function $f$.

To implement these ideas in the language of pre-categories, one has to introduce an intermediate pre-category whose objects are families of function $f^{t}$ satisfying the conditions of Proposition \ref{morse_floer_compactness}, and whose morphisms are the Floer complexes at time $t = t_0$.  By forgetting the family, we obtain an $A_{\infty}$ quasi-equivalent to $\sT \Fuk(\ctorus{n},M_{t_0,1})$, while the use of homotopy maps defines a $A_\infty$ functor to $\Morse^{\bZ^{n}}(\sQ_{t_0,1})$ which one can easily see is an $A_{\infty}$ equivalence by an energy filtration argument.   Leaving the details to the reader,  we conclude
\begin{thm}
The categories $\sT \Fuk(\ctorus{n},M_{t_0,1})$ and $\Morse^{\bZ^{n}}(\sQ_{t0,1})$ are equivalent as $A_{\infty}$-pre-categories. \noproof
\end{thm}
We now return to the practice of dropping the subscripts from $M_{t,1}$ and all associated objects.

\section{HMS for toric varieties: Morse = \v{C}ech} \label{morse=cech-section}

In this section, we shall prove the equivalence of the category $\Morse^{\bZ^{n}}(\sQ)$ with a DG category whose objects are holomorphic line bundles on $X$, and whose morphisms are \v{C}ech complexes.   We will use the following version of \v{C}ech theory:

Given a sheaf of vector spaces (or free abelian groups) $\sS$, and an ordered cover $\{ O_{k} \}_{k=0}^{N}$ of a (nice) space, such that $\sS$ is acyclic on all iterated intersections of the elements of our cover, we consider the \v{C}ech complex 
\[  \check{C}^{i}(\sS)  \equiv \bigoplus_{k_0 < \ldots < k_i} \Gamma( O_{k_0} \cap \cdots \cap O_{k_i} ,\sS ), \]
with the differential on $\Gamma( O_{k_0} \cap \cdots \cap O_{k_i} ,\sS ) $ defined by
\begin{equation} d \phi = \sum_{\stackrel{k}{k_0 < \ldots < k_j < k < k_{j+1} < \ldots < k_i}} (-1)^{j+1} \phi|_{ O_{k_0} \cap \cdots O_{ k_j} \cap O_{k} \cap  O_{k_{j+1}} \cap \cdots \cap O_{k_i} } . \end{equation}
Further, we define a cup product
\[ \check{C}^{i}(\sS_1) \otimes  \check{C}^{j}(\sS_2) \to  \check{C}^{i+j}(\sS_1 \otimes  \sS_2) \] 
which, on $\Gamma( O_{k_0} \cap \cdots \cap O_{k_i} ,\sS_1 )  \otimes \Gamma( O_{l_0} \cap \cdots \cap O_{l_j} ,\sS_2 ) $ is defined by
\begin{equation} \label{cech-product} \phi \otimes \psi \mapsto 
\begin{cases} \phi \otimes  \psi|_{ O_{k_0} \cap \cdots \cap O_{k_i} \cap  O_{l_1} \cap \cdots \cap O_{l_j}}  & \textrm{if } k_i = l_0 \\
0 & \textrm{otherwise.}  \end{cases} \end{equation}
It is easy to check that the above cup product is associative.  Since our covers will be labeled by the cells of $Q$, we will assume that we have chosen an ordering of these cells.

\subsection{Another cover of the moment polytope} \label{sec:covers_moment_polytope}
We inductively construct a cover of the moment polytope analogous to the one used in \cite{abouzaid}*{Section 5.2}.  We will write $B_{c}(Z)$ for (closed) balls of radius $c$ about a subset $Z$.

\begin{figure}
\centering
\input{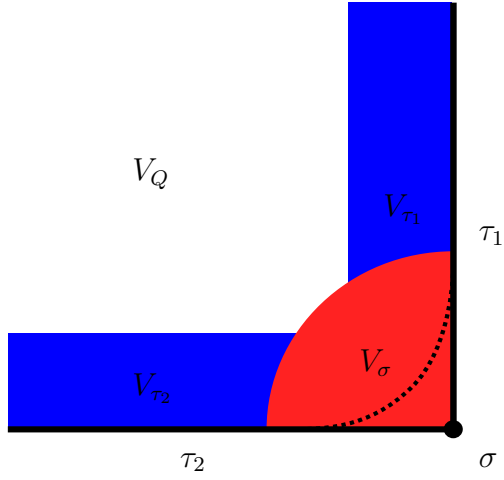}
 \caption{The cover $V_{\tau}$ in dimension $2$. The dashed line represents the boundary of $\partial \sQ$}
 \label{neighbourhood_system}
 \end{figure}

Consider a sequence of strictly positive real numbers $\epsilon_{n-1} < \ldots < \epsilon_{0}$, and inductively define $V_{\tau}$ to be the closure of
\[  B_{\epsilon_{k}} (\tau) - \bigcup_{\sigma \subset \partial \tau} \inte(V_{\sigma}) ,\]
for $\tau$ a cell of dimension $k$.  For intuition, the reader should consult Figure \ref{neighbourhood_system}.

\begin{prop} 
If $\epsilon_{n-1} \ll \ldots \ll \epsilon_{0} \ll 1$, then the cells $V_{\tau}$ define a cellular subdivision of $Q$ which is dual to the barycentric subdivision.
\end{prop}
\begin{proof}[Sketch of proof]
It is clear that the intersections of neighbourhoods of cells of the same dimension are localized near their common boundaries.  In particular, picking $\epsilon_{k} \ll \epsilon_{k-1}$ guarantees that the cells $V_{\tau}$ do not have overlapping interior (the construction only guarantees this if the cells of $Q$ have different dimension).  Once this condition is satisfied, there is a bijective correspondence between the maximal dimensional cells of the cellular subdivision and the original cells of $Q$, so it suffices to understand when two of these cells intersect.  In particular, we must prove that the $n-k$-dimensional cells of the cellular decomposition we constructed are in bijective correspondence with iterated inclusions of cells $ \sigma_0 \subset \ldots \subset \sigma_{k} \subset Q $.  

We illustrate the case $k=1$; no new ideas are need for the higher dimensional cases.  Assume that $\sigma$ and $\tau$ are respectively of dimension $i$ and $j$.  Note that if $\sigma$ and $\tau$  do not have a common boundary cell, the distance between them is so large that any reasonable choice of $\epsilon_0$ guarantees that no intersection occurs.  If they do share a boundary cell $\rho$ of dimension $k$, then the intersection of $\epsilon_i$ and $\epsilon_j$ neighbourhoods of $\sigma$ and $\tau$ occur in a neighbourhood of $\rho$.  Choosing  $\epsilon_i$ and $\epsilon_j$ much smaller than $\epsilon_k$ therefore guarantees that the boundaries of the cells $V_{\sigma}$ and $V_{\tau}$ do not intersect unless, possibly, if $\sigma \subset \tau$.

We must now prove that all such intersections do, in fact, occur.  This is clear if $\sigma$ has codimension $1$ in $\tau$.  By taking a normal slice, we reduce the problem to the case where  $\sigma$ is a vertex and $\tau = Q$.  So we must show that
\[ S_{\epsilon_0}(\sigma) - \bigcup_{\sigma \neq \tau' \subset \partial Q} V_{\tau'} \]
is a non-empty cell, where $S_{\epsilon_0}(\sigma)$ is the sphere of radius $\epsilon_0$ in $Q$.  Again, this can be clearly achieved by picking constants $\epsilon_{n-1} \ll \ldots \ll \epsilon_{0}$.  Indeed, for such a choice of constants, the $n-1$-dimensional volume of the $S_{\epsilon_0}(\sigma)$ is much larger than the volume of its intersection with the cells $V_{\tau'}$.
\end{proof}
\begin{rem}
The specific choice of coefficients depends on geometry of the polytope $Q$.  Even if we restrict to polytopes corresponding to smooth toric varieties, the ``angles" at vertices can still be arbitrarily acute.  This, of course, is a result of the fact that, by using a metric, we are restricting the symmetries to $O(n,\bZ)$, whereas constructions in tropical geometry are usually invariant under the larger group $GL(n,\bZ)$.
\end{rem}

\begin{cor} \label{cellular-cover}
If, for each cell $\tau \in Q$, $W_{\tau}$ is a sufficiently small open neighbourhood of $V_{\tau}$, then non-empty intersections
\[  W_{\sigma_0} \cap \ldots  \cap W_{\sigma_k} \] 
correspond to nested inclusions $\sigma_0 \subset \ldots \subset \sigma_k$. \noproof
\end{cor}

We observe that if we pick $\sQ$ arbitrarily close to $Q$, we can ensure that the above cover yields a cover of $\sQ$.  Further, the construction of $\sQ$ and of tropical Lagrangian sections relied on a choice of a cover $O_{\tau}$ of $\sQ$.

\begin{lem}
If the parameters which appear in the construction of tropical Lagrangian sections are be chosen to be sufficiently small, then 
\[ O_{\tau} \subset \bigcup_{\sigma \subset  \bar{\tau}} V_{\sigma} \] \noproof
\end{lem}

In particular, the cells $\sQ \cap V_{\tau}$ define a cellular subdivision of $\sQ$, which we denote $\check{\sQ}_{b}$.  Abusing notation, we will still refer to its maximal cells as $V_{\tau}$.

\subsection{Tropical Lagrangians and simplicial cochains}

In order to use the theory developed in earlier sections to relate (relative) simplicial cochains to Morse cochains, we must know that the gradient flow of functions defining tropical Lagrangian sections is well-behaved at the boundary.  The technical result that we need is:

\begin{lem} \label{barycentric_subdiv_good}
If $f$ is a smooth function on $\sQ$ defining the lift of a tropical Lagrangian section whose gradient vector field is transverse to the cells of $\check{\sQ}_{b}$, then, for each cell $\tau$, either $ \partial^{+}_f \sQ \cap V_{\tau} $ is empty, or its inclusion in $V_{\tau}$ is a deformation retract.
\end{lem}
\begin{proof}
It suffices to prove that $\partial^{+}_f \sQ \cap V_{\tau}$ is contractible when it is not empty.  Our transversality assumption implies that we may restrict our attention to maximal cells $V_{\tau}$.  

We first prove that $ \partial^{+}_f \sQ \cap  O_{\tau}$ is contractible for each cell $\tau$.  On this intersection, $\partial \sQ$ splits as a product with one factor $\tau$ and, in  $\partial \sQ \cap O_{\tau}$, the gradient flow of $f$ is only allowed to vary in the direction of $\tau$.  By taking a normal slice, the problem is reduced to the case where $\tau$ is a vertex, and the gradient flow of $f$ is constant.  Let us write $v$ for the corresponding vector in $\bR^n$, and $\partial^{+}_{v} Q$ for the (closure) of the set of maximal faces of $Q$ on which the vector $v$ points outwards.

By construction, $\partial \sQ$ is a smoothing of the convex polytope $Q$, and is itself convex.   Observe that $O_{\tau} \cap \partial \sQ$  is a convex subset of $\partial \sQ$ (it was constructed using cutoff functions with radial support).  Assume (for simplicity) that $\partial \sQ$ is strictly convex. In this case, since the Gauss map to $S^n$ has no critical points and $ \partial^{+}_f \sQ $ is the inverse image of a hemisphere, the subset  $O_{\tau} \cap \partial \sQ$  is contractible.  In our situation, the Gauss map is not a diffeomorphism, but can be easily seen to preserve contractibility.

To prove that $ \partial^{+}_f \sQ \cap  V_{\tau}$ is contractible,  we note that for each $\sigma$ with $\tau \subset \partial \sigma$, we have a non-empty intersection $V_{\tau} \cap  O_{\sigma}$.   We can apply the same analysis to these various intersections to conclude that $ \partial^{+}_f \sQ \cap V_{\tau} \cap  O_{\sigma}$ is again contractible.  Therefore $ \partial^{+}_f \sQ \cap  V_{\tau}$ is covered by contractible sets, so it suffices to prove that the geometric realization of this cover is contractible.
However,  $ \partial^{+}_f \sQ \cap V_{\tau}  \cap  O_{\sigma}$ is $C^0$ close to $ \partial^{+}_v Q \cap V_{\tau}  \cap  O_{\sigma}$  so the combinatorics of this cover of 
  $ \partial^{+}_f \sQ \cap V_{\tau} $ are equivalent to those of the cover of $\partial^{+}_v Q \cap V_{\tau} $ by the intersections  $ \partial^{+}_v Q \cap V_{\tau}  \cap  O_{\sigma}$.  Since   $\partial^{+}_v Q \cap V_{\tau} $ is easily seen to be contractible, the result follows.

\end{proof}

We define $\partial^+_f \check{\sQ}_b$ to be the union of the cells of $\check{\sQ}_b$ which intersect $\partial^{+}_f \sQ$ (this is not a subcomplex, but rather the complement of one).  We write  $\partial^+_{f} Q_b$ for the dual subcomplex of $Q_b$, and $\partial^+_{f} Q$ for the cells of $Q$ corresponding to vertices of $\partial^+_{f} Q_b$. As a direct consequence of the above analysis, we have
\begin{cor} \label{outwards-tropical}
If $\tau$ is a maximal cell of $Q$, then $\tau \subset \partial^+_{f} Q$  if and only if there is a vertex $\sigma$ of $\partial \tau$ such that $\grad(f)(\sigma)$ is outwards pointing on $\tau$.  Further, $\partial^+_{f} Q$ is the closure of its maximal cells. \noproof
\end{cor}

Let 
\[ H \co \partial Q  \to \bZ \]
be a function which is constant on each maximal cell, (the function is likely to be discontinuous on the $n-2$ skeleton). One natural way to obtain such a function $H$ is to consider a lattice vector $u \in \bZ^n$. Indeed, since $Q$ is an integral polytope, we have (integral) normal vectors $\nu_{\tau}$ which allows us to define a function
\begin{align*}
u \co \partial Q & \to \bZ \\
u|_{\tau} & = \langle u, \nu_{\tau} \rangle.
\end{align*}

As in Appendix \ref{singular-cochains-manifolds}, we set
\[ \partial^+_{H} Q \] 
to be the closure of the maximal cells where $H$ is positive.  For each such function, we can define a sheaf $\sS(H)$ on $Q$ by
\[ \Gamma(U, \sS(H)) = \begin{cases} 0 &\textrm{if } U \cap \partial^+_{H} Q \neq \emptyset \\
\bC & \textrm{otherwise.} \end{cases} \]
Using the ordering of the cells of $Q$, we write $W_k = W_{\sigma_k}$, we observe that
\[ \Gamma(W_{k_0} \cap \cdots \cap W_{k_i} , \sS(H)) = \begin{cases} 0 &\textrm{if } W_{k_j} \in \partial^+_{H} Q \textrm{ for all } 0 \leq j \leq k_i\\
\bC & \textrm{otherwise.} \end{cases} \]
The fact that
\[ \partial^+_{H_1+H_2} Q  \subset  \partial^+_{H_1} Q  \cup  \partial^+_{H_2} Q,\] 
implies that there is a map of sheaves
\begin{equation}\label{cup-sheaves} \sS(H_1) \otimes   \sS(H_2) \to \sS(H_1+H_2).\end{equation}

We consider the DG pre-category
\[ \Cech^{\bZ^n}(Q) \]
whose objects are such integral valued function on $\partial Q$ and where morphisms between $H_0$ and $H_1$ are given by
\[\bigoplus_{u \in \bZ^n}  \check{C}^*(  \sS(H_1-H_0+u)), \]
and multiplication is given by the composition of the cup product with the map on \v{C}ech complex induced by Equation \eqref{cup-sheaves}.

\begin{prop}
$\Cech^{\bZ^n}(Q) $  and  $\Morse^{\bZ^n}(\sQ)$ are equivalent as $A_{\infty}$ pre-categories.
\end{prop}
\begin{proof}
We pick a normal vector field $\nu_{\partial \sQ}$ on $\partial \sQ$ which agrees with the integral normal vector $\nu_{\tau}$ on every maximal cell $\tau$ on $O_{\tau} \cap \partial \sQ$. Starting with a function $f$ defining a tropical Lagrangian section, we consider two functions
\[ H(f)  \co \partial \sQ \to \bR  \, \, \, \, \, \, \, \,  \, \, \, \,  \hat{H}(f) \co \partial Q \to \bZ.\]
The first function is defined as in Appendix \ref{singular-cochains-manifolds} to take value $ \langle \grad(f)(p),  \nu_{\partial \sQ} \rangle$ at a point $p$.  The second function, on every maximal cell $\tau$, is given by
\[ \langle \grad(f)(p),  \nu_{\tau} \rangle\]
for any $p \in O_{\sigma} \cap \partial \sQ$ with $\sigma \subset \partial \tau$.  The function $\hat{H}(f)$ is well defined since the restriction of the gradient vector field of $f$ to $O_{\sigma} \cap \partial \sQ$ varies only in directions tangent to $\sigma$, hence the inner product with the normal vector of $\tau$ does not depend on $p$ or $\sigma$.  By taking $\sigma$ to be a vertex, we see that $\hat{H}(f)$ is, indeed, integral valued.

To define the functor, we consider an intermediate category $ \Rel^{\bZ^n}(Q_b) $ whose objects are the functions $H(f)$, for $f$ an object of $\Morse^{\bZ^n}(Q)$, and where morphisms between $H_0$ and $H_1$ are given by relative simplicial
cochains 
\[ \bigoplus_{u \in \bZ^n} C^*(Q_b,  \partial^+_{H_1-H_0+u} Q_b).\]
\begin{claim}
 $\Cech^{\bZ^n}(Q)$ and $ \Rel^{\bZ^n}(Q_b) $ are equivalent $DG$ categories.
\end{claim}
\begin{proof}
The functor take $H(f)$ to $\hat{H}(f)$.  We claim that there is a natural isomorphism between the chain complexes on either side.  Indeed, in our construction of $V_{\tau}$ and $W_{\tau}$, we already observed that there is a bijective correspondence between $k$-cells of the barycentric subdivision and intersections of $k+1$ open elements of the cover $\{ W_{\tau} \}_{\tau \subset Q}$.   It is easy to check that taking cochains relative cells in $  \partial^+_{H(f_1)-H(f_0)+u} Q_b $ exactly corresponds to taking coefficients in the sheaf $\sS(\hat{H}(f_1)-\hat{H}(f_0)+u)$, since this sheaf is trivial on open sets which intersect $\partial^+_{H_1-H_0+u} Q_b$.  Both have the effect of excluding generators corresponding to iterated inclusions of cells where all of the cells lie in $\partial^+_{H_1-H_0+u} Q_b$.  That the differential and the composition agree can be seen by the explicit formulae (e.g, Equation \eqref{cech-product} is just rewriting Equation \eqref{relative-cup}  in the dual basis).\end{proof}

Having proved the claim, it now suffices to prove that $ \Rel^{\bZ^n}(Q_b) $ and $\Morse^{\bZ^n}(\sQ)$ are equivalent.    This is done by interpreting ${Q}_{b}$ as the triangulation dual to $\check{\sQ}_b$.  We are then in the context where we can apply Corollary \ref{cell-morse}.  There are two issues to consider
\begin{enumerate}
\item  We're working with functions that have specific behaviour at the boundary (i.e. the tropical condition and the integrality condition on the functions $H$ away from the $n-2$ skeleton).  These two conditions simply restrict the set of objects we're considering, without affecting the morphisms.  In other words we consider subcategories of $\Rel(Q_b)$ and $\Morse(\sQ)$.  An equivalence of categories as in Corollary \ref{cell-morse} restricts to an equivalence of categories between every subcategory and its image.
\item The morphisms in the categories $ \Rel^{\bZ^n}(Q_b) $ and $\Morse^{\bZ^n}(\sQ)$ are $\bZ^n$ graded.  However, by construction, the differential and the composition respect this grading. 
\end{enumerate}
\end{proof}

\subsection{The \v{C}ech category of line bundles}
In this section, we complete the proof of the main theorem by constructing an equivalence from  $\Cech^{\bZ^n}(Q) $  to the \v{C}ech category of line bundles.  We're taking the dual point of view to that of \cite{fulton}*{Chapter 3}. Let $X$ be a projective toric variety, and, as before, let $Q$ denote its moment polytope.  The $i$-cells of $Q$ correspond to invariant, irreducible, $i$ dimensional toric subvarieties of $X$.  

In order to have a chain-level model for these cohomology groups, we will use the  \v{C}ech complexes associated to the canonical cover by $X$ by affine toric varieties.  Equivalently, the polyhedron $Q$ has a natural cell division to which we can associate a cover by open stars of cells.  Define
\[ U_i = U_{\sigma_i} \equiv \bigcup_{\tau | \sigma_i \subset \bar{\tau}}  \tau .\]
This gives a open cover of $Q$, but, considering the inverse image of cells under the moment map, we also obtain a cover of $X$.  Basic techniques in toric geometry, see \cite{fulton}*{Chapter 3}, imply that this is an affine cover (the affine toric variety corresponding to a cell $\sigma_i$ is defined by the normal cone of $\sigma_i$).

Let $L$ be a line bundle on $X$, and let us choose an (algebraic) trivialization over the open dense torus contained in $X$.  Since algebraic functions on $\ctorus{n}$ correspond to Laurent polynomials, we can identify the lattice $\bZ^n$ with a basis of sections of $L|_{\ctorus{n}}$.  Such sections extend to meromorphic sections of $L$ over $X$.  In particular, for each Laurent polynomial $u$, we can record the order of the pole or zero of the extension to every toric divisor.  We take the negative of the order, and obtain a $\bZ$ valued function on the facets of $Q$, which we extend to a (not necessarily continuous) function
$H(L,u) \co \partial Q \to \bZ.$
Let us write $\partial^+_{H(L,u)} Q$ for the closure of the set where $H(L,u)$ is positive.  This is the set of poles of the extension of $u$ to the toric variety $X$.  The \v{C}ech complex $\check{C}^{*}(L)$ splits as a direct sum labeled by lattice vectors $u \in \bZ^n$ with $u$-factor
\begin{equation} \label{cech-u-summand} \check{C}_u^{d}(L) \equiv  \bigoplus_{\substack{j_0 < j_1 < \cdots < j_{d} \\ \sigma_{j_k} \nsubseteq \partial^+_{H(L,u)} Q \textrm{ for all } 0 \leq k \leq d }} \Gamma(U_{j_0} \cap \ldots \cap U_{j_{d}}, \bC)  ,\end{equation}
where $\Gamma(U_{j_0} \cap \ldots \cap U_{j_{d}}, \bC) $ consists of all constant maps from $U_{j_0} \cap \ldots \cap U_{j_{d}}$ to $\bC$. Since in fact all these intersections are connected, we simply obtain a copy of $\bC$ for each such intersection.  

To obtain a category, we set
\[ \check{C}_u^{d}(L_0, L_1) \equiv   \check{C}_u^{d}(L_1 \otimes L_0^{-1}),\] 
and define the cup product as follows:  if $\phi \in \Gamma(U_{j_0} \cap \ldots \cap U_{j_{d}}, \bC)$, and $\psi \in \Gamma(U_{i_0} \cap \ldots \cap U_{i_{d'}}, \bC)$, set
\begin{equation} \label{cech-cup-product-formula} \phi \cup \psi = \begin{cases} \phi \cdot \psi|_{ U_{j_1} \cap \ldots \cap U_{j_{d}} \cap U_{i_1} \cap \ldots \cap U_{i_{d'}}}  & \textrm{if } j_{d}=i_0 \\ 0 & \textrm{otherwise.}
\end{cases} \end{equation}
The above formula can be interpreted as a cup product
\[ \check{C}_{u'}^{d'}(L_1, L_2) \otimes \check{C}_u^{d}(L_0, L_1) \to  \otimes \check{C}_{u+u'}^{d+d'}(L_0, L_2)   .\]

\begin{defin}
The DG category \[ \Cech(X) \] has objects line bundles on $X$, morphisms the \v{C}ech complexes with respect to the affine toric cover of $X$ defined by Equation \eqref{cech-u-summand}, and composition given by the product on \v{C}ech complexes defined by Equation \eqref{cech-cup-product-formula}.
\end{defin}

\begin{prop}
There is an equivalence of categories 
\[ \Cech(X)  \mapsto \Cech^{\bZ^n}(Q)\] 
\end{prop}
\begin{proof}
At the level of objects, the functor takes a line bundle $L$ to the function $H(L,u)$.  Scrutinizing the morphisms, we observe  that Equation \eqref{cech-u-summand} is simply a low-tech method of recording the \v{C}ech complex of $\sS(H(L,u))$ with respect to the cover  $U_i$.  Since both covers are acyclic, the inclusion $W_i \subset U_i$ induces an equivalence of categories.
\end{proof}

\appendix

\section{A second look at semi-tropical degenerations} \label{ap:second_look}
In this Appendix, we will modify the construction of semi-tropical degenerations and admissible Lagrangians  in order to prove Proposition \ref{morse_floer_compactness}.  First, we will prove a sharper convergence result than the one used in \cite{abouzaid}, then we will discuss the three parts of Proposition \ref{morse_floer_compactness}
\subsection{Choosing $M_{t,1}$ $C^{k}$-close to $M_t$}
Fix the partition $O_{\tau}$ of $\bR^n$, and a time $t_0$ for which we can carry through the construction of $M_{t_0,s}$.  We will write $O_{\tau}(t)$ for $\log(t) \cdot O_{\tau}$.  The following subsets of $Q$ will be useful in what follows:
\begin{align*}
O_{\partial Q} (t) &  =\bigcup_{\tau \in \partial Q}  O_{\tau}(t) \\
O_{k}(t) & = \bigcup_{\stackrel{\tau \in \partial Q}{\dim(\tau) \leq k}} O_{\tau}(t) \\
O_{Q}(t) & =  \log(t) \cdot Q - O_{\partial Q}(t) .
\end{align*}
We will also write $O^{\bC^{*}}_{\tau}(t) \subset \ctorus{n}$ for the inverse image of $O_{\tau}(t)$ under the logarithmic moment map, and by extension the meaning of $O^{\bC^{*}}_{k}(t)$ should be clear.

While the distance between $\partial \sQ_{t,s}$ and $Q$ grows logarithmically, we can in fact choose the semi-tropical amoeba degeneration such that, for any $s$, the distance between $\partial \sQ_{t,0}$ and $\partial \sQ_{t,s}$ decreases with $t$.

\begin{prop} \label{C-k-close}
For any integer $k$, and each $\epsilon$ sufficiently small, there is an appropriate choice of parameters such that for sufficiently large $t$, the $C^{k}$ distance between $ M_t $ and $ M_{t,s} $ is proportional to $t^{-\epsilon}$. \noproof
\end{prop}

Here, the $C^{k}$ norm is taken with the respect to the metric on $\ctorus{n}$ induced by the flat metric on $\bR^{n}$.  The proof of the above proposition follows from elementary estimates of the type performed in \cite{abouzaid}*{Section 3}.  We explain the main ideas:

First, recall that tropical geometry produces a polyhedra complex $\Pi_{\nu} \subset \bR^n$ (tropical amoeba) from a Laurent polynomial $W$ and a map $\nu \co A \to \bR$, where $A \subset \bZ^n$ is the subset of monomials which have non-zero coefficients.  For simplicity, we assume that the origin lies in $A$ with coefficient $-1$, that $\nu$ vanishes at the origin, and that all other coefficients are equal to $1$.  In particular, we shall write
\begin{equation} W_t = - 1 + \sum t^{-\nu(\alpha)} z^\alpha .\end{equation}
We assume that $\nu$ induces a maximal subdivision, in particular the components of the complement of $\Pi_{\nu}$ are labeled by the elements of $A$, and each monomial is represented.  Concretely, consider the function 
\[ u \mapsto \langle u, \alpha \rangle - \nu(\alpha) \]
on $\bR^n$. The component $C_{\alpha}$ is the one where the function corresponding to $\alpha$ dominates all others.  The most elementary result relating tropical and complex geometry is:
\begin{lem}
If $\cA_t$ is the amoeba of $M_t = W_t^{-1}(0)$, then, as $t \to + \infty$, the rescaled amoeba $\frac{1}{\log(t)} \cA_{t}$ converges in the Hausdorff topology to $\Pi_{\nu}$. \noproof
\end{lem}
The idea of the proof (and the proof itself) is extremely simple.  Consider a point which is $\epsilon \log(t)$ away from $\log(t) \Pi_{\nu}$.  It must lie in one of the components $\log(t) C_{\alpha}$.  Using the fact that in the coordinates coming from the identification of $\ctorus{n}$ with $T^{*} \bR^{n} /\bZ^n$, we have $z^{\alpha} = e^{\langle \alpha, u + i \theta \rangle}$, we easily compute that the monomial of $W_{t}$ labeled by $\alpha$ dominates all others by a factor of $e^{\log(t) \epsilon}$.  If $t$ is large enough, the remaining terms are too small for the sum to vanish.   In particular, no element of $W^{-1}_{t}(0)$ can project to a point whose distance to $\log(t) \Pi_{\nu}$ is larger than $\epsilon \log(t)$.

To see to what extent this estimate is optimal, we consider the example where
\[ W_t = -1 + t^{-1} z_1 + t^{-1} z_2 .\]
The tropical amoeba has a vertex at $(1,1)$ at which horizontal, vertical, and diagonal infinite edges meet.  In particular, $(\log(t), \log(t))$ is the vertex of the rescaling.  Note that the point
$( \log(t/2), \log(t/2) ) $
lies on the amoeba $\cA_t$.  The distance between this point and $\log(t) \Pi_{\nu}$ is independent of $t$ and equal to $\log(\frac{1}{2})$; in particular, the amoeba $\cA_t$ does not converge to $\Pi_{\nu}$.  However, if $M_{t,\nu}$ is the complex tropical hypersurface lying over $\log(t) \cdot \Pi_{\nu}$ then it is easy to prove a general result:
\begin{lem} \label{C-k-close-flat-case}
If $O_{n-2}$ contains the $ \epsilon$ neighbourhood of the $n-2$ skeleton of $\Pi_{\nu}$, then the distance in the complement of $O_{n-2}(t)$  between the amoeba $\cA_{t}$ and $\log(t)\Pi_{\nu}$ is bounded in norm by a constant multiple of $t^{-\epsilon}$.  Moreover, the distance in the $C^k$ topology between corresponding submanifolds of $\ctorus{n} - O^{\bC^{*}}_{n-2}(t)$ is also proportional to $t^{-\epsilon}$.  
\noproof
\end{lem}

Again, the proof is completely elementary. To simplify the explanation, we will assume  that we are working near $\log(t) C_{0}$, the region where the constant term dominates.  This is the only situation we care about in the desired application.

We illustrate $C^0$ convergence.  Let us restrict our attention to an $n-1$ cell of $\log(t) \Pi_{\nu}$ where the function defined by some monomials $\alpha$ has norm $1$.   Consider a point outside the $\delta$-neighbourhood of this cell, say in the direction of $\alpha$.  We know that, at such a point, the monomial corresponding to $\alpha$ has norm larger than $e^{\delta}$.  In particular, the difference 
\[ -1 + t^{-\nu(\alpha)} z^{\alpha} \]
has norm larger than $\delta$. Now, the norm of the remaining monomials is essentially bounded by $t^{-\epsilon}$ away from the $\epsilon \log(t)$ neighbourhood of the $n-2$-skeleton.  This means that if we pick $\delta$ to be $K t^{-\epsilon}$ for some fixed constant $K$, these remaining monomials are not large enough for the sum to vanish.  In particular, every point of $\cA_{t}$ which is more than $\epsilon \log(t)$ away from the $n-2$-skeleton of $\log(t) \Pi_{\nu}$ lies within $K t^{- \epsilon}$ of the $n-1$-skeleton.  In order to prove $C^{0}$-closeness in $\ctorus{n}$, it suffices to take into account the norm of both the real and imaginary parts of $W_{t,0}$.

The same analysis yields  $C^{k}$ bounds.  Indeed,  we're using a flat metric with respect to the $(u,\theta)$ variables, so if we write $z^{\alpha} = e^{ \langle \alpha, u + i \theta \rangle} $, we compute that the action of $dz^{\alpha}$ on tangent vectors is given by
\begin{equation} dz^{\alpha} ( \_ ) = \langle \alpha, \_ \rangle e^{ \langle \alpha, u + i \theta \rangle}  \end{equation}
In particular, the norm of $d z^{\alpha}$ is equal to $|\alpha| \cdot |z^{\alpha}|$, so the bounds on $C^0$ norms of each monomial imply the bounds on $C^1$ norms.  We conclude that away from the $\log(t) \epsilon$ neighbourhood of the $n-2$ skeleton, the leading order term of $dW_{t}$ has norm which is on the order of $t^{\epsilon}$ larger than the norm of any other term.  This proves $C^1$-closeness, and the reader may easily check that the same estimate holds for higher derivatives.

To extend the analysis beyond the maximal cells, we consider a cell $\sigma \subset \Pi_{\nu}$ of arbitrary dimension, and the polynomial $W^{\sigma}_t$ consisting of the monomials in $W_{t}$ corresponding to the maximal cells which ``meet at $\sigma$."  We write $M^{\sigma}_{t}$ for the zero level set of this function. Assuming that $O_{\tau}$ contains an $\epsilon$ neighbourhood of $\tau$, for every cell $\tau$, the same proof technique yields the following result:
\begin{lem} \label{lem:C^k-distance_to_model}
In the $C^k$ norm, the distance between $M_{t}$ and $M^{\sigma}_t$ is proportional to $t^{-\epsilon}$ on $O^{\bC^{*}}_{\sigma}(t)$ (the subset of $\ctorus{n}$ which projects (under the logarithmic moment map) to $O_{\sigma}(t)$).    \noproof
\end{lem}
The proof of Proposition \ref{C-k-close} should now be clear, since the function $W_{t,1}$ is essentially patching the various polynomials $W^{\sigma}_t$ over the skeleton of $\Pi_{\nu}$.  As a consequence of the fact that the $C^k$-norm of $W^{\sigma}_t$ is bounded by the $C^0$-norm of its monomial terms, we conclude
\begin{cor} The $C^k$-norms of $W_t$ and $W_{t,1}$ are uniformly bounded in $\log(t) C_{0}$  (the neighbourhood where the constant term dominates).  \end{cor}

\subsection{Construction a $C^k$-small symplectomorphism}
The next step of the construction of the category of tropical Lagrangian sections involved the choice of a symplectomorphism $\phi_{t,1}$ of $\ctorus{n}$ mapping $M_t$ to $M_{t,1}$.  We assume $\epsilon$ has been fixed as in Proposition  \ref{C-k-close}.
\begin{prop} \label{C-k-convergence}
For any integer $k$, and $\delta$ sufficiently small, there is an appropriate choice of parameters such that for sufficiently large $t$, the symplectomorphism
\[ \phi_{t,1} \co ( \ctorus{n}, M_t) \to (\ctorus{n}, M_{t,1}) \]
is supported in a neighbourhood of $M_t$ of size $\delta$ and has $C^{k}$ norm controlled by $\frac{t^{-\epsilon}}{\delta^k}$. Moreover, we can assume that $\phi_{t,1}$ preserves the positive real locus in $\ctorus{n}$.  \noproof
\end{prop}

To produce the desired symplectomorphism, we appeal to Weinstein's theorem for symplectic submanifolds ``with bounds."   First, we prove a bound on the size of the Weinstein neighbourhoods of $M_t$:

\begin{lem} \label{lem:weinstein_neighbourhood_fixed_size}
There exists a constant $\delta_0$ such that, for sufficiently large $t$, the $\delta_0$ neighbourhood of $M_{t}$ is embedded and is contained in a Weinstein neighbourhood of $M_t$.  Moreover, for $t$ sufficiently large, the curvature of $M_{t}$ is uniformly bounded.
\end{lem}
\begin{proof}
By the previous section,  the hypersurfaces $M_t $ and $M^{\sigma}_{t}$ are uniformly $C^k$-close in  $O^{\bC^{*}}_{\sigma}(t)$ so it suffices to prove the result for each $M^{\sigma}_{t}$.  Note that this is the product of a flat hyperplane with the zero level set of a general hyperplane in a complex torus of dimension $\codim(\sigma)$.  We can always write this family as the zero level set of
\begin{equation} \label{eq:standard_hyperplane}  -1 + \sum_{i=0}^{k}  t^{- \nu_i} z^{\alpha_i}  \end{equation}
where the vectors $\alpha_i$ are linearly independent and form a basis for the integral subspace of $\bZ^{n}$ that they span.  Let $A_{\sigma}$ denote the matrix whose rows are $\alpha_i$, and $\vnu_{\sigma}$ the column vector whose entries are $\nu_i$.  It is easy to check that translation (in the base $\bR^n$) by any vector $\vb$ which is solution to 
\[ \vb \cdot A_{\sigma} = b \vnu_{\sigma} \]
intertwines the hypersurfaces $M^{\sigma}_t$ and $M^{\sigma}_{t+b}$, and hence their Weinstein neighbourhoods.  The desired result follows.
\end{proof}

The proof of Proposition \ref{C-k-convergence} will now follow from general considerations.  Let $Y$ be a symplectic manifold equipped with a Riemannian metric, and $M$ a symplectic submanifold.  By Weinstein's results, there exists a neighbourhood of $M$ which is symplectomorphic to a neighbourhood of the zero section in the total space $E$ of the normal bundle of $M$.  Choose $\delta < 1 $ such that the $\delta$-neighbourhood of $M$ is contained in such a Weinstein neighbourhood.
\begin{lem}
If $M_{s}$ is a family of symplectic submanifolds for $s \in [0,1]$ with $M_0 = M$ whose distance from $M$ in the $C^{k}$ topology (for $k > 0$)  is bounded by $\epsilon$, and $\epsilon < \delta$ then there exists a symplectomorphism
\[ \phi_{s} \co (Y, M) \to (Y, M_s) \]
satisfying the following properties
\begin{enumerate}
\item The $C^{k}$ norm of $\phi_s$ is bounded by  $C \frac{\epsilon}{\delta^k}$ with proportionality constant depending only on the curvatures of $Y$, $M$, and $\omega$.
\item The support of $\phi_s$ is contained in a $\delta$ neighbourhood of $M_0$.
\end{enumerate}
\end{lem}
\begin{proof}[Sketch of proof:]
The assumptions of $C^{k}$-closeness imply that $M_s$ lies in the $\delta$ neighbourhood of $M$, and is moreover a graph of a $C^{k}$-small section of the normal bundle under the identification from Weinstein's theorem.  Choosing our constant $C$ sufficiently large, we can perform all estimates in the normal bundle.  Note that the restriction of the symplectic form to each $M_s$ is $C^k$-close to the pullback of the symplectic form from $M$.  Further, the two forms are clearly cohomologous.  In particular, applying the usual proof of Moser's stability theorem for symplectic forms, we find, along each $M_s$, a vector field $X_s$ with values in $TY$ whose $C^k$ norm is controlled by a universal constant times $\epsilon$, and such that the flow of $X_s$ defines a symplectomorphism from $M$ to $M_1$.

Using the symplectic form, we turn $X_s$ into a $1$-form on $M_s$ valued in $T^*Y$.  We find a family $H_s$ of functions which generate a hamiltonian flow mapping $M$ to $M_1$.  Near $M$, the $C^{k+1}$ norm of $H_s$ (and hence the $C^k$-norm of the flow) is still controlled by $\epsilon$.  However, we require that each $H_s$ vanish outside the $\delta$-neighbourhood of $M$ in order to be able to extend the flow to $Y$.  Choosing an appropriate cutoff function we can still ensure that the $C^k$ norm is controlled by $\frac{\epsilon}{\delta^k}$ up to a fixed multiplicative constant.\end{proof}
\begin{rem}
Note that we have much flexibility in the previous construction.  In particular, if $L$ is a Lagrangian submanifold of $Y$ such that the intersections $L_s = L \cap M_s$ are Lagrangian submanifolds whose $C^k$ distance is bounded by $\epsilon$, then we may choose vector fields $X_s$, hence a symplectomorphism $\phi_{1}$, which preserve $L$ and still satisfy the desired bounds.  \end{rem}

\subsection{Controlling Morse flows}
We now prove the first part of Proposition \ref{morse_floer_compactness}.  We first introduce yet another partition of the moment polytope $Q$ whose elements we will denote by $P_{\tau}$, and whose $\log(t)$ rescaling we will denote by $P_{\tau}(t)$.   We will require the following properties:
\begin{enumerate}
\item For each cell $\tau$, we have inclusions
\[ \bigcup_{\sigma \in \partial \tau} O_{\sigma} \subset \bigcup_{\sigma \in \partial \tau} P_{\sigma} \subset \bigcup_{\sigma \in \partial \tau} V_{\sigma}, \]
\item and, if $P_{k}$ is the union of all $P_{\tau}$ labeled by cells $\tau$ of dimension less than or equal to $k$, then $\sQ_{t,1} - P_{k}(t)$ is convex. 
\end{enumerate}

The easiest way to produce such a cover is to use supporting hyperplanes for each cell $\tau$, i.e. hyperplanes which intersect $Q$ precisely at $\tau$.  For example, at a vertex $v$, we can set $P_{v}$ to be the intersection of $Q$ with the half-space defined by a hyperplane parallel to a supporting  hyperplane at $v$.  Note that $O_{v}$ can be assumed to be arbitrarily small, so it is easy to  achieve the condition 
\[ O_{v} \subset P_{v} \subset V_{v} .\]
Next, one would define $P_{\tau}$ for an edge $\tau$ by taking the intersection of $Q$ with a half-space whose boundary is parallel to a supporting hyperplane at $\tau$ (we can define such a hyperplane to be the product of $\tau$ with a supporting hyperplane of the projection of $Q$ to the orthogonal complement).   Details for the remainder of the (inductive) construction are left to the reader.  As before, we will write $P_{\partial Q}(t)$ for the union of all the sets $P_{\tau}(t)$ for $\tau \subset \partial Q$, and  $P_Q(t)$ for its complement in $Q$.

Next, consider any function $f$ whose restriction to $\partial \sQ_{t,1} \cap P_{\sigma}(t)$ is of the form 
\begin{equation} \label{raise_zero_section} f_{\sigma} + \langle u , \psi \rangle  + \rho \end{equation}
where $f_{\sigma}$ depends only the directions tangent to $\sigma$, $\psi$ is a fixed integral vector, and for which there exists a constant $C$ such that the differential of $\rho$ satisfies 
\begin{equation} \label{eq:bounds_derivative_near_boundary}0  < \langle d \rho , v \rangle \leq  C |v| r_{\partial \cQ_{t,1}}  \end{equation}
whenever the vector $v$ is perpendicular to $\sigma$ and points towards the interior of $Q$.  In other words, the differential of $\rho$ is positive on the inward pointing cone at $\sigma$, and the positivity is controlled by the distance to the boundary.

\begin{lem} \label{lem:gradient_trajectories_stay_inside_general}
If the constant $C$ in Equation \eqref{eq:bounds_derivative_near_boundary} is sufficiently small, then every gradient trajectory of $f$ with endpoints in $P_{Q}(t)$ lies entirely within $P_{Q}(t)$.
\end{lem}
\begin{proof}
We will prove that no such gradient trajectory can intersect $P_{\sigma}(t)$ by considering the component of the flow in the plane perpendicular to $\sigma$.  By induction, we assume that if $k = \dim(\sigma)$ there is no gradient trajectories which moves from the complement of $P_{k-1}(t)$ into $P_{k-1}(t)$ and then exits into $P_{Q}(t)$.

By $C$ sufficiently small, we mean that $|d \rho| \ll 1$.  If $\psi$ is tangent to $\sigma$, then Equation \eqref{eq:bounds_derivative_near_boundary} implies that the gradient flow points away from the boundary so that no gradient trajectory can enter $P_{\sigma}(t)$ from $P_{\tau}(t)$ where $\tau$ has higher dimension.  Assuming the vector $\psi$ is not tangent to $\sigma$, the gradient flow of $f$ on $P_{\sigma}$ is $C^{1}$ close to the gradient flow of $f_{\sigma} + \langle u , \psi \rangle$, so that the component of $\grad(f)$ in the plane perpendicular to $\sigma$ is $C^1$ close to a linear flow.  In particular, the convexity of $\sQ_{t,1} - P_{k}(t)$ implies that a gradient flow line which enters $P_{\sigma}(t)$ either escapes to the boundary, or moves into $P_{\tau}(t)$ for $\tau \subset \partial \sigma$, at which point we use the inductive assumption to conclude that the gradient flow line cannot return to $P_{Q}(t)$.
\end{proof}

If $f$ defines a tropical Lagrangian section, then its restriction to $ V_{\sigma}(t_0)$ can be written as
\begin{equation} f^{t_0}_{\sigma} + \langle u , \psi \rangle  + \rho^{t_0} \end{equation}
where $\rho$ vanishes, together with its differential, at the boundary.  More precisely, sufficiently close to the boundary, $\rho^{t_0}$ is $C^2$-close to a constant multiple of the square of the distance to the boundary.  If this constant is positive, then $\rho^{t_0}$ satisfies condition \eqref{eq:bounds_derivative_near_boundary}. In particular,  we can deform any function $f^{t_0}$ defining a tropical Lagrangian section so that the term $\rho^{t_0} $  satisfies conditions \eqref{eq:bounds_derivative_near_boundary} on $P_{\sigma}(t_0)$ for any sufficiently small constant $C$, in which case the restriction of $f^{t_0} - v$ to $P_{\partial Q}(t_0)$ has no critical points for any integral covector $v$. More generally, we have:
\begin{lem}
If $f^{t_0}$ defines a tropical Lagrangian section over $\sQ_{t_0,1}$, and its restriction to $P_{\sigma}(t_0)$ satisfies Conditions \eqref{raise_zero_section} and \eqref{eq:bounds_derivative_near_boundary}, then there exists a family of functions $f^{t}$ defining tropical Lagrangian sections over $\sQ_{t,1}$ and which satisfy the following additional properties
\begin{enumerate}
\item The restriction of $f^{t}$ to $P_{Q}(t)$ is given by
\begin{equation} \label{eq:rescale_f} f^t (u)  =  \frac{\log(t)}{\log(t_0)} f\left(\frac{\log(t_0)}{\log(t)} u \right)    \end{equation}
\item The restriction of $f^{t}$ to $P_{\sigma}(t)$ can be decomposed as
\begin{equation} \label{raise_zero_section_t} f^t_{\sigma} + \langle u , \psi \rangle  + \rho^{t} \end{equation}
with $d \rho^{t}$ bounded independently of $t$, and satisfying the bound \eqref{eq:bounds_derivative_near_boundary}.
\end{enumerate}
\end{lem}
\begin{proof}
The expression \eqref{eq:rescale_f} defines a function on $ \frac{\log(t_0)}{\log(t)} \sQ_{t_0,1}$.  A look at the construction of $\sQ_{t,1}$ reveals that this rescaled copy of $\sQ_{t_0,1}$ is contained in $\sQ_{t,1}$.  In each element $V_{\sigma}(t)$ of our cover, there is a unique way to extend 
\[ f_{\sigma} \left(\frac{\log(t_0) u}{\log(t)} \right) \]
to the complement of $\frac{\log(t_0)}{\log(t)} \sQ_{t_0,1} $ preserving the property that the function is independent of the directions perpendicular to $\sigma$.  

To complete the construction, it suffices to explain how to choose $\rho^t$.  Again, the construction is inductive, so  we start with a neighbourhood of a vertex.  Start by extending $\rho^{t_0} \left(\frac{\log(t_0) u}{\log(t)} \right)$ to $P_{v}(t)$ by setting it to vanish outside the intersection with $ \frac{\log(t_0)}{\log(t)} \sQ_{t_0,1}$; this function is once differentiable.  Next, we pick any function $\rho'$ which agrees with the square of the distance to the boundary near $\partial \sQ_{t_0,1}$, and satisfies Equation \eqref{eq:bounds_derivative_near_boundary} for some positive number $C$.   Rescaling $\rho'$ so that it is sufficiently small, we may interpolate between $\rho'$ and $\rho^{t_0} \left(\frac{\log(t_0) u}{\log(t)} \right)$ using cutoff functions so that Equation \eqref{eq:bounds_derivative_near_boundary} still holds.  A $C^1$-small perturbation gives a smooth function with the desired property.

If $v$ is a vertex of an edge $\tau$, then the restriction of $\rho^{t}_{v}$ to the intersection of $P_{v}(t)$ with $V_{\tau}(t)$  satisfies the condition \eqref{eq:bounds_derivative_near_boundary} for the positive cone of $\tau$.  Convexity of this property allows us to extend it to a function $\rho^{t}_{\tau}$  on the open set $P_{\tau}(t)$, and the rest of the construction proceeds by induction.
\end{proof}

Note that our construction of the functions $\rho^{t}_{\sigma}$ does not guarantee a uniform lower bound in Equation \eqref{eq:bounds_derivative_near_boundary}  as we may have to rescale the intermediate function $\rho'$ to ensure that the cutoff still satisfied  \eqref{eq:bounds_derivative_near_boundary}.  However, we may ensure that the norm of $d \rho^{t}$ is uniformly bounded in $t$.  Applying Lemma \ref{lem:gradient_trajectories_stay_inside_general}, we conclude the the first part of Proposition \ref{morse_floer_compactness} if we choose $\delta^{t}_1$ such that the $\delta^{t}_1$ neighbourhood of $\partial \sQ_{t,1}$ is contained in $P_{\partial Q}(t)$.  In the next sections, we will perturb $f^t$ by $C^1$-small amounts near $\partial \sQ_{t,1}$ to ensure that the remaining parts of Proposition \ref{morse_floer_compactness} also hold; we are free to do this as long as the bound \eqref{eq:bounds_derivative_near_boundary} holds for some sufficiently small constant $C$.

\subsection{Compactness for $J_{t,s}$ holomorphic curves}

We will now modify the family $L^t$ in order to satisfy the second property listed in Proposition \ref{morse_floer_compactness}.  The idea is to ensure that there is a constant $\epsilon_t$ such the images of the Lagrangians $L^t$ agree with a line segment near the circle of radius $\epsilon_t$ in $\bC$.  We will then prove a compactness theorem analogous to the one proved in \cite{abouzaid}*{Lemma 2.8}, without requiring that the segment extend all the way to the origin.   

We start by considering the Laurent polynomials $W^{\sigma}_t$.  As in the proof of Lemma \ref{lem:weinstein_neighbourhood_fixed_size}, the fact that translation in $\bR^{n}$ intertwines the set of solutions of $W^{\sigma}_t$ for different values of $t$ implies that for any constant $c_1$, there is a uniform constant $2\epsilon$ independent of $\sigma$ and $t$ such that the parallel transport of the positive real locus in $M^{\sigma}_{t}$ along the straight-line segment in $\bC$ from the origin to $(-2\epsilon, 2c_1 \epsilon)$ is a graph over the zero-section, and its defining function satisfies Equation \eqref{eq:bounds_derivative_near_boundary}  for $C = 2c_1$.  This latter property holds because we know that the function defining this graph is $C^2$-close to a positive multiple of the square of the distance to the boundary, which manifestly satisfies the desired condition in a sufficiently small neighbourhood.

Applying Proposition \ref{C-k-convergence}  and Lemma \ref{lem:C^k-distance_to_model}, which together imply that, as $t$ grows, $M_{t,1}$ becomes $C^k$-close at every point to some hypersurface $M^{\sigma}_{t}$ we conclude:

\begin{lem} \label{lem:straight_near_circle}
For any fixed constants $ 0 < 2c_1 < C$, there exists a constant  $\epsilon$ so that the following properties hold for all  $t$ sufficiently large:
\begin{enumerate}
\item the support of $\phi_{t,s}$ lies entirely within the inverse image of the ball of radius $\epsilon$ under $W_{t,1}$, and
\item the parallel transport of $\partial L^t$ along a segment from the origin to $(-\epsilon, c_1 \epsilon)$ defines a section over a neighbourhood of $\partial \sQ_{t,1}$, and
\item the Lagrangian $L^{t}$ can be chosen to agree with the above parallel transport near the circle of radius $\epsilon$.  Further, the image of $L^{t}$ under $W_{t,1}$ lies in the left half plane, with  all interior intersection points  with the zero section mapping to the complement of the disc of radius $\epsilon$.
\end{enumerate}
\end{lem}
\begin{proof}
Since the first two parts follow immediately from the fact that as $t$ goes to infinity,  the $C^k$-norm of $\phi_{t,1}$ as well as its support may be assumed to shrink, we only discuss the last part.  The proof that the image of $L^{t}$ may be assumed to lie in the left half plane can be done in two steps.  First, we consider the standard hypersurface
\begin{equation} \label{eq:standard_poly} -1 + \sum_{i} z_{i} \end{equation}
The complement of the amoeba of this standard polynomial has a uniquely determined component whose inverse image in $\ctorus{n}$ is contained in the locus where the norm of all variables $z_i$ is bounded by $1$.  The boundary of this subset of $\ctorus{n}$ consists of complex numbers $z_i$ whose norms lie on the component of the boundary of the standard amoeba  corresponding to the positive real locus; this implies that the image of the boundary under the polynomial \eqref{eq:standard_poly} lies in the left half plane (in fact, it intersects the $y$-axis only at the origin).  Since interior points  are obtained by rescaling an appropriate boundary point by a positive real number bounded in norm by $1$, and all coefficients in \eqref{eq:standard_poly} are positive, this entire component maps to the left half plane.

Note that by acting by $\SL_{n}(\bZ)$ and by translations, we conclude that the analogous component of the complement of the amoeba of $M^{\sigma}_{t}$ has inverse image which also maps, under $W^{\sigma}_{t}$, to the left half-plane.  This implies the desired result for $W_{t,1}$ by using $C^{k}$-closeness to polynomials $W^{\sigma}_{t}$ in the domains $V_{\sigma}(t)$ (note that in $V_{Q}(t)$, $W_{t,1}$ is $C^k$-close to $-1$).
\end{proof}

We now prove that choosing $L_{i}^{t}$ to agree with the parallel transport of its boundary near a circle of radius $\epsilon$ implies a compactness theorem for $J_{t,s}$ holomorphic discs.

\begin{lem}  \label{lem:straight_near_circle_compact}
If $L^{t}$ is chosen as in Lemma \ref{lem:straight_near_circle}, then every $J_{t,s}$ holomorphic disc with marked points mapping to interior intersection points of $L^t$ with the zero section projects under $W_{t,1}$ to the complement of the disc of radius $\epsilon$.
\end{lem}
\begin{proof}
Assume by contradiction that there exists such a holomorphic disc $u$ which intersects the disc of radius $\epsilon$.  Note that $J_{t,s}$ agrees with $J$ near the inverse image of this circle, so that $W_{t,1}$ is $J_{t,s}$ holomorphic near the circle of radius $\epsilon$ for every $s \in [0,1]$.  Since all intersection points map outside the circle of radius $\epsilon$, the image of $u$ must intersect the circle of radius $\epsilon$; let $z$ be such an intersection point whose inverse image under $u$ contains an interior point of the disc.  Note that the degree of $u$ at $z$ (number of inverse images counted with signs and multiplicity) is necessarily positive by holomorphicity.

\begin{figure}[t]   \centering
 \input{lagrangian_straight_circle.pstex_t}
   \caption{}
   \label{fig:lagrangian_straight_circle}
\end{figure}

Consider the path $\gamma$ shown in Figure \ref{fig:lagrangian_straight_circle}.  By applying  the maximum principle to the composition of $W_{t,1}$ with the projection to the $x$-axis, we find that the image of any holomorphic curve under $W_{t,1}$ cannot intersect the complement of the disc of radius $\epsilon$ in the right half plane.  In particular, the inverse image of $z_{+}$ under the composition  $W_{t,1} \circ u$ is empty.  However, the difference between the degree of $u$ and $z$ and $z_{+}$ is given by the intersection number between $\gamma$ and the image of the boundary of the disc under $u$.  The boundary of the disc consists of a pair of intervals connecting the two intersections points along the zero section and $L^t$.  Since the images of the intersection points lie in the left half plane and away from the disc of radius $\epsilon$, we obtain two paths in the plane, with end points away from the disc of radius $\epsilon$ about the origin, and which agree respectively with the negative $x$-axis, and the straight line of slope $-c_1$ near the circle of radius $\epsilon$.  In particular, the intersection number between $\gamma$ and the boundary vanishes, which yields the desired contradiction since the degree of $u$ at $z$ cannot vanish.
\end{proof}

To prove Part (2) of Proposition \ref{morse_floer_compactness} we choose $ \delta_{2}^{t}$ so the inverse image under $W_{t,1}$ of the circle of radius $\epsilon$ in $L^t_i$ projects to the complement of the $\delta_{2}^{t}$ neighbourhood of $\partial \sQ_{t,1}$ under the logarithmic moment map.  Let us write
\begin{equation}   \sQ_{t,1} (\delta_{2}^{t}) \end{equation}
for the complement of this neighbourhood in $ \sQ_{t,1}$.

The previous result implies that only the interior points of a $J_{t,s}$-holomorphic curve can map outside to the complement of $ \sQ_{t,1} (\delta_{2}^{t}) $, and that the support of  $\phi_{t,s}$ is contained in this complement.  In order to exclude this possibility, we recall that the boundary of $\partial \sQ_{t,1}$ is convex, so we may choose $t$ large enough and $\delta_{2}^{t}$ small enough, so that  the boundary of $ \sQ_{t,1} (\delta_{2}^{t}) $ is also convex.

The image of any $J_{t,s}$-holomorphic curve under $\phi^{-1}_{t,s}$ is a $J$-holomorphic curve, where $J$ is the standard complex structure.  Since $\phi_{t,s}$ preserves the boundary of $ \sQ_{t,1} (\delta_{2}^{t}) $, it suffices to prove that no $J$-holomorphic curve in $\ctorus{n}$ can have a projection which escapes $ \sQ_{t,1} (\delta_{2}^{t}) $.    The maximum principle for holomorphic curves in $\bC^{n}$ immediately implies this (See Lemma \ref{lem:maximum_principle}).

\subsection{Compactness for $\epsilon J$ holomorphic curves}

We finally complete the proof of Proposition \ref{morse_floer_compactness} by proving that $\epsilon J$ holomorphic curves with boundary on $L^t$ cannot escape towards the boundary of $\partial \sQ_{t,1}$.  Our main technique will involve using the product decomposition near every cell $\tau$ of $\partial Q$, then the maximum principle for holomorphic curves with boundary in the factor orthogonal to $\tau$.

Let $\sV$ be a convex smooth hypersurface in $\bR^{n}$, and write $\sV^{\bC}$ for its inverse image in $\bC^{n}$ under the projection to the real coordinates $\bC^{n} \to \bR^{n}$.
\begin{lem} \label{lem:maximum_principle}
No holomorphic curve in $\bC^{n}$ with boundary on the zero-section can touch $\sV^{\bC}$ from the interior.  \end{lem}
\begin{proof}
If $h$ is a function on $\bC^{n}$ depending only the real variables whose $0$-level set is $\sV^{\bC}$, then the kernel of the restriction of $d^{c} h$ to $\sV^{\bC}$ agrees with the complex distribution.  The convexity of $\sV$ implies that $d (d^{c} h)$ is positive on this distribution, so $\sV^{\bC}$ is $J$-convex, and we can apply the maximum principle at the interior, see \cite{mcduff}.  Since, in addition, $d^{c} h$ vanishes on $\bR^{n} \subset \bC^{n}$, we have von Neumann boundary conditions, so we can also apply the maximum principle at the boundary.  
\end{proof}

We now return to the problem at hand, and pick any $\delta^{t}_{3}$ which is smaller than both $ \delta^{t}_{1} / 2 $ and $ \delta^{t}_{2} / 2$. The boundary of $\sQ_{t,1}$ is convex, so it is easy to pick a convex hypersurface $\sV \subset \sQ_{t,1}$ satisfying the following properties:
\begin{enumerate}
\item $\sV$ is contained in the $2 \delta^{t}_{3}$ neighbourhood of $ \partial \sQ_{t,1}$ but in the complement of the $\delta^{t}_{3}$ neighbourhood, and
\item The intersection of $\sV$ with $V_{\tau}(t)$ splits as an orthogonal product with a factor consisting of directions tangent to $\tau$, and the other given by a convex hypersurface $\sV_{\tau}$.
\end{enumerate}
We pick a function $h$ whose zero level set is $\sV$ and for which $0$ is a regular value, and write $\alpha = d^{c} h$ after extending $h$ to $\bC^{n}$.

Consider the Lagrangian $L^{t}_{0}$, which, on each open set $V_{\sigma}(t)$ is given as the graph of the differential of the function
\begin{equation}
f_{\sigma} + \langle \psi, \_ \rangle \end{equation}
as in Equation \eqref{raise_zero_section}.  
\begin{lem} \label{lem:decomposition_vanishing}
The $1$-form $d^{c} h$ vanishes on $L^{t}_{0}$.
\end{lem}
\begin{proof}
It suffices to prove this in each set $V_{\sigma}$ where both $\sV$ and $L^{t}_{0}$ split as Riemannian products.  Since $L^{t}_{0}$ is flat precisely in the directions where $\sV$ is curved the result easily follows.
\end{proof}

We will now pick $\rho^t$ carefully so that $d^{c} h$ still vanishes on $L^t$ near its intersection with $\sV$.  First, we observe the following
\begin{lem}
There is a choice of normal vector $\nu_{\sigma}$ for each face $\sigma \subset \partial Q$ satisfying the following conditions:
\begin{enumerate}
\item $\nu_{\sigma}$ lies in the inward cone at $\sigma$, and
\item if $\sigma \subset \partial \tau$, then $\nu_{\tau}$ agrees with the orthogonal projection of $\nu_{\sigma}$ along the tangent space of $\tau$.
\end{enumerate}
\end{lem}
\begin{proof}
Let $q_0$ be an arbitrary interior vertex of $Q$, and set $\nu_{\sigma}$ at a point $p \in \sigma$ to be the projection of the vector from $p$ to $q_0$ to the orthogonal complement of the tangent space of $\sigma$ at $p$.  The first property follows from the convexity of $Q$, and the second from the compatibility of orthogonal projection with restriction to subspaces.
\end{proof}

We will now require $\rho$ to satisfy the following property in each set $V_{\sigma}(t)$:
\begin{equation} \label{eq:split_rho}
\parbox{30em}{Near $\sV$, the projection of the gradient of $\rho^t$ to the plane orthogonal to $\sigma$ is constant and agrees with $\nu_{\sigma}$.
}
\end{equation}

We are forcing $d\rho^t$ to only vary in direction tangent to $\sigma$ near the intersection of $V_{\sigma}(t)$ with $\sV$.  One needs to check that there is indeed a function satisfying the above property (e.g. the corresponding $1$-form needs to be closed).  We can proceed by induction with respect to the dimension of the cell; whenever $\sigma$ is a vertex, we write $\rho_{\sigma}$ for a function whose gradient on $V_{\sigma}(t) \cap \sV$ is fixed to agree with $\nu_{\sigma}$.  At the end of the construction, $\rho_{\sigma}$ shall be the restriction of $\rho^t$ to $V_{\sigma}(t)$.  If $\sigma_1$ and $\sigma_{2}$ are the boundary vertices of an edge $\tau$, $\rho_{\sigma_i}$ defines a function on part of the boundary $V_{\tau}(t)$ whose restriction to $V_{\tau}(t) \cap \sV$ has normal component which agrees with $\nu_{\tau}$.  Consider any extension of  $\rho_{\sigma_i}$ to a function on $V_{\tau}(t)$ whose gradient has normal component $\nu_{\tau}$.  Since the space of such functions is convex, an appropriate linear combination of the functions coming from the two end points of $\tau$ gives a function $\rho_{\tau}$ which extends both $\sigma_{1}$ and $\sigma_{2}$ to a function on  
\begin{equation*} V_{\sigma_1}(t) \cup  V_{\tau}(t) \cap  V_{\sigma_2}(t) \end{equation*} satisfying the desired condition.   The reader may easily complete the construction for higher dimensional cells, and check that the resulting function $\rho^t$ can be chosen so that Condition \eqref{eq:bounds_derivative_near_boundary} holds.

The same decomposition argument as in Lemma \ref{lem:decomposition_vanishing} implies:
\begin{lem} 
If $\rho^t$ satisfies the condition described in Equation \eqref{eq:split_rho}, then the restriction of $d^{c}h$ to $L_{t}$ vanishes near $\sV$.  In particular, no $J$-holomorphic curve with boundary on the zero-section and $L_{t}$ may be tangent to $\sV$.  \noproof 
\end{lem}

Note that by picking $\rho^t$ small enough, we can ensure that the conclusions of the previous two sections still hold.  We have therefore completed the proof of the last part of Proposition \ref{morse_floer_compactness}.

\section{Moduli spaces of trees} \label{mod-space-trees}
In this appendix, we define a stratification into polyhedral cells of certain moduli spaces of trees introduced in Definition \ref{defin-shrubs}.  We begin by noting that the moduli space $\Stasheff_{d}$ admits a decomposition into subsets labeled by the topological type of a tree with $d+1$ leaves.  Let us write $[T]$ for the homeomorphism type of a tree and $\Stasheff_{d}([T])$ for the corresponding subset of $\Stasheff_{d}$.  Note that the lengths of edges defines a homeomorphism
\begin{equation} \label{eq:cell_decomp_stasheff}  \Stasheff_{d}([T]) \to (0,+\infty)^{r}\end{equation}
where $r$ is the number of non-trivial edges.  The closure $\tilde{\Stasheff}_{d}([T])$ of  $\Stasheff_{d}([T])$ in $\Stasheff_{d}$ is homeomorphic to $[0,+\infty)^{r}$, with the corner strata corresponding to topological types obtained by collapsing one or more edges of $[T]$.  Since every tree can be obtained by collapsing a trivalent tree, we conclude

\begin{lem}
The moduli space $\Stasheff_{d}$ admits a decomposition into cells which are homeomorphic to $[0,+\infty)^{d-2}$.  \noproof \end{lem}

Note that the compactification of the moduli space $\Stasheff_{d}$ corresponds to allowing edges to achieve infinite length.  We shall denote the closure of each cell $\Stasheff_{d}([T])$ in $\overline{\Stasheff}_{d}$ by $\overline{\Stasheff}_{d}([T])$.  Upon fixing an identification of $[0,\infty]$ with $[0,1]$ (say via the exponential of the negative), we conclude
\begin{lem}
The compactified moduli space $\overline{\Stasheff}_{d}$ admits a canonical decomposition into subsets $\overline{\Stasheff}_{d}([T])$ which are homeomorphic to cubes.   \noproof
\end{lem}
Note that this decomposition is compatible with the decomposition of the boundary of $\overline{\Stasheff}_{d}$ as the product of lower dimensional moduli spaces.  With this in mind, we can now describe the regularity condition on the moduli spaces of gradient trees $\overline{\Stasheff}_{d}$ stated in Definition \ref{defin:regularity_moduli_spaces_trees}.

First, we use the homeomorphisms \eqref{eq:cell_decomp_stasheff} to equip each cell $\tilde{\Stasheff}_{d}([T])$ with a smooth structure as a manifold with corners.

Recall that $\Stasheff(q, \vp)$ is the moduli spaces of gradient trees with inputs $\vp = (p_1, \ldots p_d)$ and output $q$.  Let us write $\tilde{\Stasheff}(q, \vp; [T])$ for those gradient trees whose topological type corresponds to a tree in $\tilde{\Stasheff}_{d}([T]) $. In order to prove that this is a smooth manifold, we consider the space
\[ \tilde{\cE}(q, \vp; [T]) = W^{s}(q) \times W^u(p_1, \ldots , p_d) \times \tilde{\Stasheff}_{d}([T]) ,\]
where $ W^u(p_1, \ldots , p_d) $ is the product of the descending manifolds of the critical points $\{ p_1, \ldots, p_d \}$.

Let us fix a tree $T \in \tilde{\Stasheff}_{d}$, and a basepoint $s \in T$.  Let $v^o_j$ be the finite end of each incoming edge $e_j$.  The geodesic $A$ between $v^o_i$ and $s$ determines a map
\[ \Phi^s_{j} \co W^u(p_j) \to Q \]
which is the composition of the gradient flows along the functions labeling the edges of $A$.  The gradient flow is taken for time equal to the length of the corresponding edge.  We can similarly define a map with source $W^{s}(q)$.  Given a choice of basepoints (which is unique up to homotopy since the total space of the tautological bundle over $\tilde{\Stasheff}_{d}$ is contractible), we obtain a smooth map
\begin{equation} \label{eq:evaluation_map_define_moduli_space_trees} \Phi_{[T]} \co \tilde{\cE}(q, \vp; [T]) \to Q^{d+1} .\end{equation}
The space $\tilde{\Stasheff}(q, \vp; [T])$ is inverse image of the diagonal $\Delta^{d+1} \subset  Q^{d+1}$  under this map.  Sard's theorem implies that this inverse image is a smooth manifold whenever the map $\Phi_{[T]}$ is transverse to the diagonal.  The diffeomorphism of $Q^{d+1}$ induced by gradient flow along the edges connecting different basepoints proves that this moduli space and its smooth structure are independent of the choice of basepoints.

Note that the moduli space $\Stasheff(q, \vp)$ is the union of the spaces $\tilde{\Stasheff}_{d}(q, \vp; [T])$ over their common boundary strata.
\begin{defin} \label{defin:regularity_moduli_trees}
We say that $\Stasheff(q, \vp)$ is {\bf regular} if all evaluation maps $\Phi_{[T]}$ are transverse to the diagonal.
\end{defin}
The proof of the next Lemma, which is omitted, is an exercise in Morse theory.  The key point is that whenever $\Phi_{[T]}$ is transverse for a fixed metric tree $T$, it is also transverse for nearby metric trees which can be obtained by resolving vertices of $T$ of valency greater than $2$. 
\begin{lem}
For generic choice of functions $\vf$ the moduli space $\Stasheff(\vf)$ is regular, and admits a compactification $\overline{\Stasheff}(\vf)$ into a compact manifold with boundary, admitting a decomposition into (smooth) compact manifolds with boundary. \noproof \end{lem}

The moduli space of shrubs is topologised in exactly the same way as the Stasheff moduli space.  In particular, it also admits a decomposition into smooth manifolds.  One way to see this is to observe
\begin{lem} \label{smooth_structure_interior}
There is surjective projection
\[ \Shrub_{d} \to \Stasheff_{d} \]
obtained by replacing all incoming edges with edges of infinite length.  The fiber $\tilde{\Shrub}_{d}([T])$ over each closed cell $\tilde{\Stasheff}_{d}([T])$ admits a unique smooth structure which makes this projection a fibration with fibre an open interval, and with trivialization given by the length of an external edge.
\end{lem}
\begin{proof}
The only part that requires checking is that the different trivializations determine compatible smooth structures.  This follows easily from the fact that the difference between such trivializations can be expressed in terms of the lengths of internal edges which are, indeed, smooth functions on the base.
\end{proof}
It is now clear that $\Shrub_{d}$ is not compact, and that its non-compactness comes from two phenomena; internal edges can expand to infinite length as in $\Stasheff_{d}$, and external edges can collapse to have length zero.
\begin{defin}
The moduli space of {\bf singular shrubs} $\overline{\Shrub}_{d}$ is the moduli space of singular ribbon trees with $d$ incoming edges which satisfy the following properties:
\begin{itemize}
\item The incoming edges have finite type, and
\item the outgoing edge has infinite type, and
\item all irreducible components are stable, and are either shrubs or are Stasheff trees (all external edges are infinite), and
\item the distance from any point to the incoming vertices lying above it is independent of the choice of incoming vertex.
\end{itemize}
\end{defin}
\begin{rem}
It is immediate from the above Definition that the irreducible components containing the incoming vertices are shrubs, and that all other irreducible components are Stasheff trees.
\end{rem}

The main result of this Appendix is
\begin{prop} \label{cpct_mod_space}
$\overline{\Shrub}_{d}$ is a compact manifold with having $\Shrub_{d}$ as its interior.  Moreover, it admits a decomposition into smooth manifolds with corners $\overline{\Shrub}_{d}[T]$ labeled by the topological type of the shrub.
\end{prop}
We first observe
\begin{lem}
By replacing all incoming vertices with vertices of infinite type and collapsing any unstable components, we obtain a projection map
\[\cp: \overline{\Shrub}_{d}([T]) \to  \overline{\Stasheff}_{d}([T]).\]
The fibre at a point $T \in \Stasheff_d$ is identified by the length function of an incoming edge $e$ with a closed interval $[a_e, + \infty]$. \noproof
\end{lem}

As we shall see presently, this lemma essentially determines the topology and smooth structure of each cell away from the inverse image of $\partial \overline{\Stasheff}_d$.  We shall call the part of the boundary lying over $\Stasheff_d$ the {\bf horizontal boundary}  and the part lying over the boundary of $\overline{\Stasheff}_{d}$ the {\bf vertical boundary} of $\overline{\Shrub}_{d}$.  In order to understand this later subset, we consider
\begin{defin}
A collection of successive vertices of a ribbon tree forms a {\bf branch} if the complement of the minimal tree containing them is connected.
\end{defin}

\begin{defin}
An {\bf order-preserving} partition of $\{1, \ldots, d\}$ is a partition into subsets $A_1, \cdots, A_{d'}$ such that $a_i < a_j$ if $a_i \in A_i$ and $a_j \in A_j$.  We will denote the set of such partition into $d'$ subsets by $\Part(d,d')$.
\end{defin}

Given a tree $T'$ with $d'$ incoming edges, and such a partition, we obtain a unique tree  $T$ with $d$ incoming edges,  such that those that belong to a subset that is not a singleton have vanishing length.  We have
\begin{lem}
The horizontal boundary of $\overline{\Shrub}_{d}([T])$ is covered by the images of the inclusions of
\[ \overline{\Shrub}_{d'}([T']) \times \Part(d,d') \to \overline{\Shrub}_{d}([T]),\]
and by the section of $\cp$ corresponding to all incoming edges having infinite length.

The vertical boundary of $\overline{\Shrub}_d([T])$ is covered by the images of the inclusions for each $1 \leq n \leq d$ of
\[ \coprod_{\sum_{1 \leq i \leq n} d_i = d} \overline{\Shrub}_{d_1}([T_1]) \times \cdots \times \overline{\Shrub}_{d_{n}}([T_n])  \times \overline{\Stasheff_{n}} ([T_0]) \to \overline{\Shrub}_{d}([T]).\] with the restriction that there exists some $d_i \neq 1$ .

\end{lem}
\begin{proof}
The section of $\cp$ takes a Stasheff tree with $d$ inputs and attaches half infinite intervals (isometric to $[0, + \infty)$) to all incoming vertices.  This covers the part of the horizontal boundary where some incoming edge (hence all incoming edges) reach infinite length. The remainder of the horizontal boundary corresponds to an edge collapsing to have length $0$.  

Given a tree $T$, we consider the largest subtree $T_0$ containing the outgoing vertex, but no binary vertex.  If such a tree is included properly in $T$, then its complement is a disjoint union of trees $T_i$ that occur in lower dimensional moduli spaces and $T$ is in the image of the inclusion
\[ \coprod_{\sum_{1 \leq i \leq n} d_i = d} \overline{\Shrub}_{d_1}([T_1]) \times \cdots \times \overline{\Shrub}_{d_{n}}([T_n])  \times \overline{\Stasheff_{n}}([T_0]) \to \overline{\Shrub}_{d}([T]).\]
If $d_i=1$ for all $i$, then each $\Shrub_{d_i}$ is a point and we are considering a component of the horizontal boundary.
\end{proof}

One way of stating the naturality of the above boundary inclusions is to observe
\begin{lem}
The following diagram commutes
\[
\xymatrix{ \overline{\Shrub}_{d_1}([T_1]) \times \cdots \times \overline{\Shrub}_{d_{n}} ([T_n]) \times \overline{\Stasheff_{n}} ([T_0]) \ar[d] \ar[r] &  \overline{\Shrub}_{d} ([T])\ar[d]\\ \overline{\Stasheff}_{d_1}([T_1]) \times \cdots \times \overline{\Stasheff}_{d_{n}} ([T_n]) \times \overline{\Stasheff_{n}} ([T_0])  \ar[r] &  \overline{\Stasheff_d}([T]), }
\] 
where the vertical arrows are given by the various projections $\cp$, and the bottom horizontal arrow is the usual gluing map which makes $\overline{\Stasheff_d}$ into a topological operad.
\end{lem}

We need a few more pieces of notation.  If $(d_1, \ldots, d_n)$ is a partition of $d$, let $D_i$ denote the partial sums 
\[D_i = \sum_{j \leq i} d_j.\]
For each tree $T$, let
\[ T(D_i, D_{i+1}) \]
denote the minimal connected subtree of $T$ containing all vertices whose labels range from $D_i$ to $D_{i+1}-1$.  Let $U(d_1, \ldots, d_n)$ denote the subset of $\overline{\Shrub}_{d}$ consisting of trees such that all $T(D_i,D_{i+1})$ are branches. Similarly, let $e(D_i, D_{i+1})$ denote the outgoing edge at the outgoing vertex of $T(D_i, D_{i+1})$ (In particular, $e(D_i, D_{i+1})$ is not contained in this tree).

We are now ready to give the
\begin{proof}[Proof of Proposition \ref{cpct_mod_space}]
The smooth structure on $\overline{\Shrub}_{d}[T]$ is already determined in the interior by Lemma  \ref{smooth_structure_interior}.  Along the horizontal boundary, it can defined by requiring that $\cp$ is a fibration, with trivialization given by the length function of an incoming edge along the component of the horizontal boundary covered by the images of
\[ \overline{\Shrub}_{d'}[T'] \times \Part(d,d') \to \overline{\Shrub}_{d}[T].\]
Along the other component of the horizontal boundary (where the length of incoming edges takes infinite length), we require that the smooth trivialization be given by the exponential of the negative length function of an incoming edge.  Again, this smooth structure is independent of which incoming edge is chosen since the length of incoming edges differ from each other by the lengths of the internal edges, which are smooth functions on the base.

It suffices therefore to define the topology (and smooth structure) along the vertical boundary.  We shall proceed  by induction; in the base case, both $\overline{\Shrub}_{1}$ and  $\Shrub_{1}$ are points.  Our induction hypothesis is that
\begin{itemize}
\item Whenever the length of an edge is finite, it defines a smooth function to $\bR$.
\item When the length of an edge reaches infinity, then the exponential of its negative is a smooth function.
\item The subsets of ${\Shrub}_{i}([T])$ consisting of homeomorphic trees with no vertices of valency greater than $4$ are open (such trees are not required to be combinatorially equivalent, i.e. adding binary vertices will not change the homeomorphism type).
\end{itemize}

This last condition might seem superfluous; its presence is a consequence of the fact that the fibre over a boundary point of $\overline{\Stasheff_n}$ can consist of rather complicated strata, and that we have not yet specified how these strata fit together.

By the previous lemma, every point $T$ on the vertical boundary lies in the image of 
\[ \overline{\Shrub}_{d_1} ([T_1]) \times \cdots \times \overline{\Shrub}_{d_{n}} ([T_n])\times \overline{\Stasheff_{n}} ([T_0])\to \overline{\Shrub}_{d}([T]).\] 
By induction we have a topology (and smooth structure as a manifold with corners) on the source of the above map.  We shall require these maps to be diffeomorphisms onto their images.  There is no obstruction to achieving this since the images of such maps intersect along products of lower dimensional moduli spaces whose smooth structure is already determined.

Consider the set $U(d_1, \ldots, d_n)$ defined in the paragraph preceding the beginning of this proof.  We impose the condition that these sets are open.  This condition is, by induction, compatible with all previous choices.  We define a map
\[\ft: U(d_1, \ldots, d_n) \to  \overline{\Shrub}_{d_1} ([T_1]) \times \cdots \times \overline{\Shrub}_{d_{n}} ([T_n])\times \overline{\Stasheff_{n}} ([T_0])\]
by replacing each edge $e(D_i, D_{i+1})$ by an edge of infinite length.  This map is clearly surjective.  Let $U(K, d_1, \ldots, d_n)$ denote the open subset of $\overline{\Shrub}_{d}([T])$ where all edges $e(D_i,D_{i+1})$ have length greater than $K$.  If $V_0(\frac{1}{K})$ is a fundamental system of neighbourhoods of $T_0$ in the product
\[  \overline{\Shrub}_{d_1} ([T_1]) \times \cdots \times \overline{\Shrub}_{d_{n}} ([T_n])\times \overline{\Stasheff_{n}} ([T_0]),\]
we define a fundamental system of neighbourhoods of $T_0$ in $\overline{\Shrub_d}([T])$ to be
\[ U_0\left(\frac{1}{K}\right) = \ft^{-1}\left(V_0\left(\frac{1}{K}\right) \right)\cap  U(K, d_1, \ldots, d_n).\]
Again, this is clearly compatible with previous choices of topology.  To define the smooth structure as a manifold with boundary, we observe that the projection map
\[  U_0\left(\frac{1}{K}\right) \to V_0\left(\frac{1}{K}\right) \]
is trivialized by the length function of an edge  $e(D_i,D_{i+1})$ with target $(K,+\infty]$.  We determine the smooth structure by again requiring that the exponential of the negative length is smooth.  Since by induction, $V_0\left(\frac{1}{K} \right)$ is diffeomorphic to the standard model of a neighbourhood of a corner, the same is true for its product with a closed interval.
\end{proof}

\section{Orientations in Morse theory} \label{coherent-orient}

In order to define an $A_{\infty}$ structure on $\bZ$-valued Morse complexes, we must orient the moduli spaces of gradient trees.  The first step is to coherently orient the Stasheff polyhedra.  This is, in fact, relatively easy to do using the well-known identification with the moduli space of discs.   The reader can consult Section (12g) of \cite{seidel-book} for details.  Note that such an orientation induces an orientation of each maximal cell $ {\Stasheff}_{d} ([T]) $.

We now orient the moduli spaces of gradients trees which appear in the definition of the higher products in the Morse pre-category.  Consider the manifold
\[ \cE(q, \vp) = W^{s}(q) \times W^u(p_1, \ldots , p_d) \times {\Stasheff}_{d} ,\]
Let us assume that orientations for all descending manifolds have been chosen. If we write $\lambda(N)$ for the top exterior power of the tangent bundle of $N$ as a graded vector space, then the induced orientation is given by the isomorphism
\begin{equation} \label{product-orientation} \lambda \left( \Stasheff(q, \vp)\right) \cong \lambda \left( \Delta^{d+1} \right) \otimes  \lambda^{-1} \left( Q^{d+1} \right)  \otimes \lambda \left( \cE(q, \vp) \right) ,\end{equation}
where the tangent space $ \cE(q, \vp) $ at a given point is defined to be the tangent space of the corresponding cell $\tilde{ \cE}(q, \vp; [T]) $ for $[T]$ a trivalent tree; the global orientability of the Stasheff polyhedra implies that the choice of a resolution of a non-trivalent tree does not affect the signs.

We will twist this orientation by
\begin{multline} \label{twist} \sigma(q, \vp)  = (m+1)  \cdot \\ \left(  \sum_{j=0}^{d} ( \deg(q) + \sum_{k \leq j} \deg(p_k) )  + \dim(\Stasheff(q, \vp)) \cdot (1+m+d+\deg(q) \right).\end{multline}

\begin{rem}
In Equation \eqref{product-orientation}, we used the inverse of  vector space $ \lambda^{-1} \left( Q^{d+1} \right) $,  Our convention for orienting them is as follows: If $v_1 \wedge \cdots \wedge v_n$ is an orientation of $\lambda V  $, then $ v_n^{-1} \wedge \cdots \wedge v_1^{-1}$ is the ``inverse" orientation of $\lambda^{-1} V$. 
\end{rem}

We can now complete the
\begin{proof}[Proof of Proposition \ref{signed_A_infty} ]
Let $\Stasheff(q, \vp)$ be a $1$-dimensional moduli spaces of gradient trees.  We know that its boundary consists of points corresponding to
\[ \Stasheff( q,  \vp[i]_1,  r ,\vp[d]_{i+d_2+1}) \times   \Stasheff( r, \vp[i+d_2]_{i+1})  ,\]
where
\begin{align*}
\deg(r) & = d_2 -2 + \sum_{ j = i+1}^{i+d_2} \deg(p_j), \end{align*}
 and the sequences $\vp[b]_{a}$ stand for $(p_a, \ldots, p_{b-1} , p_b)$.  Set $ d_1 = d-d_2+1 $. For simplicity, we assume that $d_2 \neq 1$, in other words that we're working entirely with stable trees.  The idea does not change in the case where $d_2=1$, though the notation is complicated by the presence of a factor associated to the group of automorphisms of moduli spaces of gradient trajectories between critical points.

It suffices to compare the induced orientation of the boundary with its orientation as a product.  Ignoring the factors $\sigma(q, \vp)$ for now, we note that an orientation of a tubular neighbourhood of the boundary is given by
\begin{align*}
& \bR\cdot \nu \otimes \lambda(\Delta^{d_1+1}) \otimes  \lambda^{-1} \left( Q^{d_1+1} \right)  \otimes \lambda \left( \cE(q, \vp[i]_1  , r , \vp[d]_{i+d_2+1} ) \right)  \\
\quad & \otimes  \lambda(\Delta^{d_2+1}) \otimes \lambda^{-1} \left( Q^{d_2+1} \right)  \otimes \lambda \left( \cE(r , \vp[i+d_2]_{i+1}) \right) ,\end{align*}
where $\nu$ is the outward normal vector at the boundary.  One can pass from this orientation to the orientation as an open subset of $\Stasheff(q, \vp)$ through simple steps.  We go through them one by one and record the corresponding change of orientation.  Note that all computations are done modulo $2$.
\begin{enumerate}
\item Recall that the $d_1+1$ factors $Q^{d_1+1}$ correspond to the output $q$ and the inputs $(p_1, \ldots, p_i, r, p_{i+d_2+1}, \ldots, p_d)$.  Writing $\lambda^{-1}(Q^{d_1+1})$ as the tensor power of the corresponding number of copies of $\lambda^{-1}(Q)$, we permute these lines so that the factor corresponding to $r$ is adjacent to  $\lambda(\Delta^{d_1+1})$.  Because we're working with the inverse of top exterior power of the tangent bundle, we have to slide the factor corresponding to $r$ past those corresponding to $p_{i+d_2+1}, \ldots, p_d $. This yields a sign given by the parity of
\begin{equation} m \cdot m (d - d_2 - i) = m (d - d_2 - i) . \end{equation}
Since the diagonal is oriented in the same way as $Q$, we have a natural identification of
\[ \lambda(\Delta^{d_1+1}) \otimes \lambda^{-1}(Q) \to \bR ,\]
and we can cancel these factors.
\item Similarly, we permute the factors of $Q^{d_2+1}$ so that the copy corresponding to $r$ appears next to $\Delta^{d_2+1}$, and cancel these two factors.  The permutation introduces a sign 
\begin{equation} m \cdot md_2 = m d_2 .\end{equation}
\item We decompose
\[ \cE(q , \vp[i]_{1}, r, \vp[d]_{i+d_2+1}) = W^{s}(q) \times W^{u}(p_{1}, \ldots, p_{i}) \times W^{u}(r) \times  W^{u}(p_{i+d_2+1}, \ldots, p_{d}) \times \Stasheff_{d_1} ,\]
and permute $W^u(r)$ with $Q^{d_1} \times  W^{s}(q) \times W^{u}(p_{1}, \ldots, p_{i})$.  This changes the orientation by
\begin{align} 
& (m-\deg(r) ) \left( m d_1 + \deg(q) + i m + \sum_{j=1}^{i} \deg(p_{j}) \right) \notag
\\ &  = (m-\deg(r) ) \left( m d_1 + \deg(q) + i (m+1) + \maltese_{i} \right).
 \end{align} 
\item We move the $W^{s}(r)$ factor of $ \cE(r , \vp[i+d_2]_{i+1}) $ past \[Q^{d_1} \times W^{s}(q) \times W^{u}(p_1, \ldots, p_i, p_{i+d_2+1}, \ldots, p_d) \times \Stasheff_{d_1} \times Q^{d_2}\] with a sign of
\begin{align}
& \deg(r) \left(  d_1 + \deg(q) + m  - \sum_{j=1}^{i} \deg(p_j) - \sum_{j=i+d_2+1}^{d} \deg(p_j) + m d_2 \right)  \notag \\
& = \deg(r)\left( d_1 + \deg(q) + m + d_1 + \deg(q) + \deg(r) + md_2 \right) \notag \\
& = \deg(r)(md_2 + m + 1 ).
\end{align}
\item We now permute $Q^{d_2}$ with the terms that separate it from $Q^{d_1}$.  The corresponding sign is
\begin{align} m d_2 \left( d_1 + \deg(q) + m (d_1-1)  - \sum_{j=1}^{i} \deg(p_j) + \sum_{j=i+d_2+1}^{d} \deg(p_j) \right) \notag \\
= m d_2 ( m (d_1 -1) + \deg(r)  ).
\end{align}
\item Recall that the factors of $Q^{d+1}$ correspond  to the sequence $(q,p_1, \ldots , p_d)$.  The terms of $Q^{d_1} \times Q^{d_2}$ are not in this same order.  Since we're taking the inverse of the canonical bundle, the desired permutation has sign
\begin{equation} m d_2 (i+1). \end{equation}
\end{enumerate}
Let us note that, after taking into account the signs which we introduced in the previous steps, our orientation of a neighbourhood of the boundary is now given by
\begin{multline*}\bR \cdot \nu \otimes \lambda(W^{u}(r)) \otimes \lambda(W^s(r)) \otimes \lambda^{-1}(Q^{d+1}) \otimes \lambda \left( W^s(q) \right) \otimes  \\ \lambda \left( W^{u}(p_1, \ldots, p_i, p_{i+d_2+1}, \ldots, p_d ) \right) \otimes \lambda \left( \Stasheff_{d_1} \right) \otimes \lambda \left( W^u(p_i+1, \ldots, p_{i+d_2}) \right) \otimes \lambda \left( \Stasheff_{d_2} \right) .\end{multline*}
We chose our orientation of $W^{u}(r)$ in such a way as to have an oriented isomorphism
\[ \lambda(W^{u}(r)) \otimes \lambda(W^s(r)) \cong \lambda(Q). \]
The next steps to arrive at our chosen orientation of $\Stasheff(q, \vp)$ are:
\begin{enumerate} \setcounter{enumi}{6}
\item Permute $ W^u(p_i+1, \ldots, p_{i+d_2})$ with $ W^{u}(p_{i+d_2+1}, \ldots, p_d )  \times \Stasheff_{d_1} $, and introduce a sign
\begin{align} 
& \left( d_2 m - \sum_{j = i+1}^{i+d_2} \deg(p_j) \right) \cdot \left(d_1 + (d-d_2-i)m + \sum_{j=i+d_2+1}^{d} \deg(p_j) \right) \notag \\
& = (d_2 (m+1) + \deg(r)) \left( i+ (d-d_2-i)m + \deg(q) + \deg(r) + \maltese_{i} \right) 
\end{align}
\item Permute the normal vector $\bR \cdot \nu$ to be adjacent to $\Stasheff_{d_1} \times \Stasheff_{d_2}$.  The corresponding sign is given by
\begin{equation} dm + \deg(q) + dm + \sum_{j=1}^{d} \deg{p_j} = d + 1 \end{equation}
\end{enumerate}

The total sign change so far is given by the parity of

\begin{multline*}
(\deg(q) + \maltese_i)(m+ d_2 (m+1)) \\
+ \deg(r)( m+1)d_2+  md_2d+ mi + md_2 + d_2 i + d + m+ 1\end{multline*}

There are two more sources of signs that we have neglected.
\begin{enumerate}\setcounter{enumi}{8}
\item There is a choice of coherent orientations on the Stasheff polyhedra such that the orientation of a neighbourhood of the boundary of $\Stasheff_d$ differs from the product orientation
\[ \bR\cdot \nu \times \Stasheff_{d_1} \times \Stasheff_{d_2} \]
precisely by 
\[ d_2 (d -i) + i + d_2 ,\]
as explained in Equation (12.25) of \cite{seidel-book}.  
\item The additional change of orientation $\sigma$ given by Equation \eqref{twist} introduces a sign equal to
\[ (m+1) ( (d_2 + 1) (\deg(q) + \maltese_{i}) + d_2 \deg(r) + dd_2 + d_2 + i)  +1+ m + d+ \deg(q) \]
as can be easily checked by comparing the contribution of the inputs and outputs to $\sigma$ in each of the three moduli spaces we're considering. 
\end{enumerate}

 The reader can easily check that the total sign change is, indeed, $\maltese_{i}$.  Since the (signed) sum of the number of points on the boundary of a $1$-dimensional manifold vanishes, we conclude that the $A_{\infty}$ equation holds.

\end{proof}

We note that the complete proof of Proposition \ref{cech-to-morse} also require orientations of moduli spaces of gradient trees, this time
\[\Shrub(p,\vC) ,\]
for a critical point $p$ and cycles $\vC = (C_1, \ldots, C_d)$.  In fact, such an orientation can be obtained through means entirely analogous to those used in the previous proof once we explain how to orient the moduli spaces of incoming-finite trees $\Shrub_{d}$.  We note that that the method we're about to use also gives orientations of the Stasheff polyhedra.

\begin{figure}
   \centering
 \input{cover_tree_chart.pstex_t}
   \caption{}
   \label{fig:cover_tree_chart}
\end{figure}

Consider a tree $T \in \Shrub_{d}$.  Let $\{v_i\}$ denote the incoming vertices of $T$.  Let $A_i$ be the descending arc containing $v_i$ which is maximal with respect to the property of containing only finite edges.  Let $\sE_{1}$ be the ordered collection of finite edges whose union is $A_{1}$, with the inverse order given by the natural orientation of the tree, i.e. the last element of $\sE_1$ is the incoming edge at $v_1$. Inductively define $\sE_i$ as the ordered collection of edges of $A_{i}$, excluding the incoming edges if $i \geq 2$, with the additional restriction that the edges appearing in $\sE_i$ do not appear in any preceding set $\sE_j$.  We write $|\sE_{i}|$ for the number of elements in such a collection.  See Figure \ref{fig:cover_tree_chart}.

\begin{lem}
If $T$ is trivalent, the lengths of the edges which appear in
\[ \bigcup_{i} \sE_{i} \]
give local coordinates for $\Shrub_{d}$ near $T$. \noproof
\end{lem}

\begin{figure}
   \centering
 \input{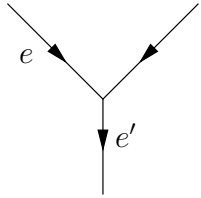}
 \caption{The path from $e$ to the node of $T$ takes a right turn at $e'$}
   \label{right-turn}
\end{figure}

Note that some of the collections $\sE_{i}$ may be empty.  We now use this canonical choice of coordinates to define an orientation.  First, we need
\begin{defin}
If $e$ is an edge of a trivalent tree $T$, then the {\bf dexterity} of $e$, denoted $r_{T}(e)$ is the number of right turns along the geodesic from $e$ to the outgoing vertex.
\end{defin}
In other words, if $A_e$ is the geodesic from $e$ to the outgoing vertex, then we can stand at every edge of $A_{e}$, looking towards the outgoing vertex, and decide whether the next edge of $A_{e}$ will correspond to a left turn or a right turn at the next vertex (See Figure \ref{right-turn}). Left and right are, of course, determined by the chosen cyclic ordering at each vertex.

To simplify notation, we write $de$ for $dl(e)$.  We define an orientation by
\begin{align*}
\Omega_{T}(v_i) & = \bigwedge_{e \in \sE_{i}} (-1)^{r_T(e)} de \\
\Omega_{T} & = \Omega_T(v_1) \wedge \cdots \wedge \Omega_T(v_d)
\end{align*}
where the edges that appear in the first expression are ordered along the path from the outgoing vertex to $v_i$.

\begin{lem}
The orientation $\Omega_{T}$ is independent of the topological type of the tree $T$.
\end{lem}
\begin{proof}
It suffices to consider neighbourhoods of trees with one quadrivalent vertex (these form hyperplanes in the moduli space of trees separating the open cells consisting of trivalent trees of a fixed topological type).  Crossing such a hyperplane, a branch $B(j,k)$ of the tree moves from one side of a trivalent vertex to the other.  Such a branch is uniquely determined by its incoming leaves $\{v_j, \ldots, v_k\}$.   See Figure \ref{fig:moduli_trees_hyperplane}.

\begin{figure}[htbp] 
   \centering
 \input{moduli_trees_hyperplane-right.pstex_t}
  \input{moduli_trees_hyperplane.pstex_t}   \input{moduli_trees_hyperplane-left.pstex_t} 
   \caption{}
   \label{fig:moduli_trees_hyperplane}
\end{figure}
Fix a tree with a unique quadrivalent vertex $v$.  Note that all our coordinates extend naturally through the hyperplane, except for the map to $\bR$ which is given by the length of the edge which is collapsed at $v$.  Let $T$ and $T'$ be nearby trees on either side of the hyperplane.  We fix a co-orientation of the hyperplane by assuming that the path from $v_{k+1}$ to the outgoing vertex includes the edge $\alpha  \subset T$ which is collapsed when passing through the hyperplane.  While the corresponding edge $\alpha' \subset T'$ leads, by necessity, to $v_{j-1}$.  In fact, we can extend $l(\alpha)$ smoothly to the other side of the hyperplane by declaring it equal to $- l(\alpha')$.  In particular, we have
\begin{align*}
d \alpha & = - d \alpha' \\
r_{T}(\alpha) & = 1 + r_{T'} (\alpha'). \end{align*}

Let $B(i,\ell)$ denote the branch of $T$ containing $v$ as its bottom-most vertex. It is obvious that for every edge that does not lie on $B(i,\ell)$, we have
\[r_{T}(e) = r_{T'}(e) .\]
We claim the same equality holds for all edge on $B(k+1,\ell)$ and for $B(j,k)$.  Indeed, we only remove a left turn from the paths with origin the edges of $B(k+1,\ell)$ to the outgoing vertex, while we permute the order of a left and right turn on corresponding paths from the edges of $B(j,k)$.  It suffices therefore to check that
\[ \Omega_{T}(v_i) \wedge \cdots \wedge \Omega_{T}(v_{j}) =   \Omega_{T'}(v_i) \wedge \cdots \wedge \Omega_{T'}(v_{j}) .\]
Observe $d\alpha$ is the first $1$-form which appears in our formula for $\Omega_T(v_j)$, while $d\alpha'$ appears in our formula for $\Omega_{T'}(v_{i})$.  Let us therefore split
\[ \Omega_{T'}(v_i) = \Omega^{1}_{T'}(v_i) \wedge d \alpha' \wedge \Omega^{2}_{T'} (v_i) .\]
The result now follows easily by observing that the dexterity of each edge in $B(i,j-1)$, changes by $1$ as $\alpha$ is collapsed, which exactly cancels the sign produced by commuting $d\alpha$ past 
\[ \Omega^2_{T}(v_i) \wedge \cdots \wedge \Omega_{T}(v_{j-1}).\]
\end{proof}

\section{$A_{\infty}$ pre-categories} \label{precat}
We recall some elementary definitions concerning $A_{\infty}$ pre-categories, due to Kontsevich and Soibelman \cite{KS}. 
\begin{defin}
An $A_{\infty}$ pre-category $\cA$ consists of the following data:
\begin{itemize}
\item A collection of objects $\Ob(\cA)$.
\item For each integer $n$, a collection of transverse sequences
\[ \Ob^{n}_{tr}(\cA) \subset  \Ob(\cA)^{n}\]
\item For each ordered pair $(X_0, X_1) \in  \Ob^{2}_{tr}(\cA)$, a $\bZ$-graded chain complex
\[ \cA(X_0, X_1) .\] 
\item For each transverse sequence of length $d+1$, $(X_0, X_1, \ldots, X_d)$, a degree $d-2$ map of graded vector spaces, called the $d$th higher product
\[ \xymatrix{\cA(X_{d-1}, X_d) \otimes \cdots \otimes \cA(X_0, X_1) \ar[rr]^{\mu_{d}} & & \cA(X_0, X_d).}\]
\end{itemize}
The defining data of an $A_{\infty}$ pre-category is required to satisfy the following compatibility conditions
\begin{itemize}
\item Every subsequence of a transverse sequence is transverse.
\item The higher products are required to satisfy the $A_{\infty}$ equation, i.e. Equation \eqref{A_infty-equation}.
\end{itemize}
\end{defin}
Kontsevich and Soibelman further defined the notion of quasi-isomorphisms in $A_{\infty}$ pre-categories, and used it to give a criterion of unitality of $A_{\infty}$ pre-categories.
\begin{defin}
A morphism $f \in \cA(X_0, X_1)$ is said to be a {\bf quasi-isomorphism} if, whenever $(Z,X_0, X_1)$, or $(X_0, X_1, Y)$ are transverse triples, the maps
\begin{align*} 
\xymatrix{\cA(X_1, Y) \ar[r]^{\mu_2(\_, f)} & \cA(X_0, Y)} \\
\xymatrix{\cA(Z, X_0) \ar[r]^{\mu_2(f, \_)} & \cA(Z,X_1) } \end{align*}
induce isomorphisms on homology.
\end{defin}
\begin{defin}
An $A_{\infty}$ pre-category is unital if for each object $X$, and each transverse sequence $S$, there exist objects $X_+$ and $X_-$ which are quasi-isomorphic to $X$ such that the sequence
\[ (X_-, S, X_+) \]
is transverse.
\end{defin}

In this setting, it is appropriate to replace the usual notion of full subcategory.
\begin{defin}
The inclusion of a full sub-pre-category $\cB$ in an $A_{\infty}$ pre-category $\cA$ is said to be {\bf essentially surjective} if every object of $\cA$ is quasi-isomorphic to an object of $\cB$.
\end{defin}
In order to simplify our notation, we deviate slightly from Kontsevich and Soibelman's definition of $A_{\infty}$ functor.
\begin{defin}
An $A_{\infty}$ functor $F$ between $A_{\infty}$ pre-categories $\cA$ and $\cB$ consists of the following data:
\begin{itemize}
\item An essentially surjective inclusion $\cA_{F} \subset \cA$.
\item A map
\[ F: \Ob(\cA_F) \to \Ob(\cB).\]
\item For each transverse sequence $(X_0, \ldots, X_d)$ of $\cA_F$, a degree $d-1$ map 
\[ F^d: \cA(X_{d-1}, X_d) \otimes \cdots \otimes  \cA(X_0, X_1) \to \cB(FX_0, FX_d).\]
\end{itemize}
This data must satisfy the following conditions
\begin{itemize}
\item The image of every transverse sequence of $\cA_F$ is transverse in $\cB$.
\item The maps $F^{d}$ satisfy the $A_{\infty}$ equation for functors 
\begin{align*} & \sum_{\substack{k \\ \l_1 + \ldots + l_k = d} } \mu^{B}_{k} \circ ( F^{l_1} \otimes \cdots \otimes F^{l_k}) = \\
& \quad \sum_{i + d_2 \leq d} (-1)^{\maltese_{i}} F^{d-d_2+1} \circ (\mathbb{1}^{d-i-d_2} \otimes \mu^{A}_{d_2} \otimes \mathbb{1}^{i}) .\end{align*}
\end{itemize}
\end{defin}
This definition differs from that of \cite{KS} where they require $\cA_F = \cA$.  However, this is irrelevant since this difference disappears upon considering quasi-equivalence classes of $A_{\infty}$ pre-categories.
\begin{defin}
An $A_{\infty}$ { \bf quasi-equivalence} is an $A_{\infty}$ functor which induces an isomorphism at the level of homology categories.

Two $A_{\infty}$ pre-categories $A$ and $B$ are {\bf equivalent} if there exists a sequence $\{ A_{i} \}_{i=0}^{n}$ of $A_{\infty}$ pre-categories and $A_{\infty}$ quasi-equivalences
\[ A = A_0 \to A_1 \leftarrow A_2 \to A_3 \leftarrow \cdots \to A_n =  B.\]
\end{defin}

\section{Cellular chains and simplicial cochains on manifolds with boundary} \label{singular-cochains-manifolds}

Let $Q$ be a smooth compact oriented manifold with boundary $\partial Q$.  We would like to rephrase the classical equivalence between the cup product on the simplicial cochains of a triangulation $\sP$ and the intersection pairing on cellular chains in a language more suitable to our goals in this paper.  A stronger result appears in \cite{mcclure} (the commutative structure is taken into account, cf. \cite{wilson}), but the techniques used seem incompatible with our specific setup. To make the definition of the cup product standard, we work with a special type of triangulation:

\begin{defin}
A triangulation $\sP$ is {\bf simplicial} if there is a complete ordering on the vertices of $\sP$ such that every cell of $\sP$ can be expressed uniquely as the  ``hull" of an ordered sequence of vertices $v_0 < \ldots < v_k$.  We write $[v_0 , \ldots , v_k] \subset \sP$ for such a cell.
\end{defin}

With this condition, one can use the full machinery of simplicial sets, though we shall not.  We now construct a category which encodes the cup product on relative simplicial cochains.

\begin{defin}
The pre-category $\Rel(\sP)$ has as objects functions 
\[ H \co \partial Q \to \bR .\]
A pair $H_0$ and $H_1$ of such functions is transverse if the difference $H = H_1 - H_0$ satisfies the following condition:  the inclusion of $H^{-1}([0, + \infty))$ into $\partial^+_{H} \sP$, the union of all cells where $H$ takes some positive value, is a deformation retract.

If $H_0$ and $H_1$ are transverse, the space of morphisms from $H_0$ to $H_1$ is given by
\[ \Rel(\sP)(H_0, H_1) = C^*(\sP, \partial^+_{H_1-H_0} \sP ).\]

A sequence $(H_0, H_1, \ldots, H_d)$ is transverse if all pairs $(H_i, H_j)$ with $i < j$ are transverse.  All higher products vanish, and the composition is defined through the inclusion $\partial^+_{H_2 - H_0} \sP \subset  \partial^+_{H_2 - H_1} \sP \cup \partial^+_{H_1 - H_0} \sP$ and the corresponding diagram:
\[ \xymatrix{ C^*(\sP, \partial^+_{H_2 - H_1} \sP ) \otimes C^*(\sP, \partial^+_{H_1 - H_0} \sP ) \ar[r] &  C^*(\sP,  \partial^+_{H_2 - H_1} \sP \cup  \partial^+_{H_1 - H_0} \sP) \ar[d] \\
& C^*(\sP, \partial^+_{H_2 - H_1}  \sP) .}\]
\end{defin}

For the sake of definiteness, we describe the chain model theory that we use.  Our treatment is standard in the sense that it can be extracted from an algebraic topology textbook \cites{hatcher,may}.  The chain complex of a simplicial triangulation is defined by
\[ C_{k}(\sP) \equiv \bigoplus_{[v_0 , \ldots , v_k]} \bZ \cdot [v_0 , \ldots , v_k].\]
Again, we reiterate the fact that our conventions are such that writing $[v_0 , \ldots , v_k]$ subsumes the fact that the vertices are listed in increasing order.  The differential is given as usual by
\[ \partial [v_0 , v_1 , \ldots , v_k] = \sum_{i=0}^{k} (-1)^{i}  [v_0 , \ldots, v_{i-1}, \hat{v_{i}}, v_{i+1}, \ldots , v_k] .\]
If $A \subset \sP$ is a subcomplex, we define the relative cochain complex 
\[ C^*(\sP, A) \equiv \{ \phi \co C_{*}(\sP) \to \bZ  \,\, | \,\,  \phi(\sigma) = 0 \textrm{ if } \sigma \in A \} .\]
For subcomplexes $A, B \subset \sP$, we define the cup product 
\begin{equation}C^*(\sP, A)  \otimes C^*(\sP, B) \to C^*(\sP, A \cup B)\end{equation}
by the formula
\begin{equation}\label{relative-cup}  \phi \cup \psi ([v_0, \cdots, v_k]) = \sum_{i=0}^{k} \phi([v_0, \cdots, v_i]) \cdot \psi([v_i, \cdots, v_k]) .\end{equation}
A priori, $\phi \cup \psi$ is a homomorphisms from $C_{*} (\sP)$ to $\bZ$, but, since every cell of $A \cup B$ lies either in $A$ or in $B$, it is easy to check that $\phi \cup \psi$ vanishes on all such cells and hence defines an element of $ C^*(\sP, A \cup B)$  (cf. \cite{hatcher}*{p. 209} where the less straightforward case of singular cohomology is discussed).  

We now proceed with a cellular version of this construction.  By a cellular subdivision $\sR$ of $Q$, we mean a covering by $n$-dimensional closed polytopes such that the intersection of any two polytopes is either empty or a common cell on the boundary.  Given such a subdivision, we define the corresponding cellular chain complex by
\[ C_{*}(\sR) \equiv \bigoplus_{C \in \sR} \bZ \cdot \Omega(C) ,\]
 where $\Omega(C)$ is the orientation line of the cell $C$.  The differential is obtained by restricting the orientation to the boundary.  Let us fix a subdivision $\check{\sP}$ dual to a simplicial triangulation $\sP$, and a corresponding constant $\epsilon >0$ which is much smaller than the diameter of all cells in $\check{\sP}$.  It will be convenient to use the following terminology:
 
 \begin{defin}
The {\bf $\epsilon$-essential part} $C - \partial C(\epsilon)$ of a cell $C \in \check{\sP}$ is the complement, in $C$, of an $\epsilon$-neighbourhood of $\partial C$.
\end{defin}

By definition, every cell $C$ of codimension $k$ in $\check{\sP}$ is dual to a cell of the form $[v_0, \ldots, v_k]$, and the link of $C$ is isomorphic to the simplex with vertices $(v_0, \ldots, v_k)$.  We will sometimes write $\sigma$ for a cell of $\sP$ and $\check{\sigma}$ for the dual cell.  Concretely, a point in the interior of $C$ lies on the closure of exactly $k$ top dimensional cells which are dual to the vertices $\{v_0, \ldots, v_k \}$.  We will have to understand cell subdivisions which are $C^1$ close to $\check{\sP}$.  To do this, it is helps to first understand the local model.

\begin{figure}   \centering
 \includegraphics{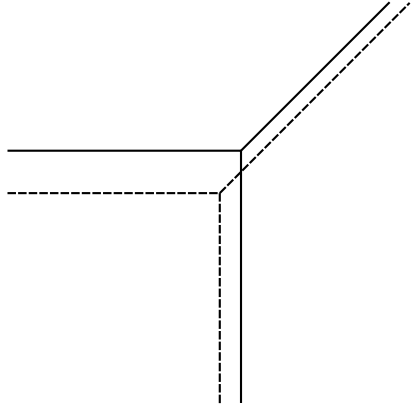}
   \caption{The one-skeleta of the polyhedral subdivisions $\Pi^2$ (full lines) and $\Pi^2(r)$ (dashed)}
   \label{simplicial_position_model}
\end{figure}

Let $\Pi^d$ be the cellular subdivision of $\bR^d$ into the $d+1$ cells, $d$ of which are of the form 
\[ \{ (x_1, \ldots, x_d) \, | \, 0, x_1, \ldots, x_d \leq x_j \}\]
for some $j$ between $1$ and $d$, and the remaining cell is the set 
\[ \{(x_1, \ldots, x_d) \, | \,  x_1, \ldots, x_d \leq 0\}.\]   It is easy to check that $\Pi^d$ is a model for the dual subdivision of a $d$-simplex $\Delta^d= [v_0, \ldots, v_d]$ with the cell where $0$  dominates dual to $v_0$, and the cell where $x_j$ dominates all the other terms dual to the vertex $v_j$ for $1 \leq j \leq d$.  Given any positive real number $r$, we may construct a new cellular subdivision of $\bR^d$ by translating $\Pi_d$ in the direction of the vector 
\begin{equation} \label{eq:formula_simplicial_perturbation} (-r/d, \ldots, -r/2, -r) .\end{equation}
We write $\Pi^d(r)$ for this new subdivision of $\bR^d$.  See Figure \ref{simplicial_position_model}.  If $\sigma$ is a cell of $\Delta^d$, we write $\check{\sigma}$ for the dual cell of $\Pi^d$ and $\check{\sigma}(r)$ for the dual cell of $\Pi^d(r)$.
\begin{lem}  \label{simplicial-model}
If  $\sigma =  [u_0, \ldots, u_k]$ and $\tau = [w_0, \ldots, w_l]$ are cells of $\Delta^d$ then their duals $\check{\sigma}$ and $\check{\tau}(r)$ have non-empty intersection if and only if $u_k \leq w_0$.
\end{lem}
\begin{proof}
Observe that  $\check{\sigma}$ and $\check{\tau}(r)$  intersect if and only if $\check{[u_i]}$ intersects $\check{[w_j]}(r)$ for each $i$ and $j$.  However, it is easy to check that this can only happen if $u_i \leq w_j$, for all $i$ and $j$, which because of our ordering conventions, is equivalent to $u_k \leq w_0$. \end{proof}

For any $r< r'$, the composition of translation by the vector of Equation \eqref{eq:formula_simplicial_perturbation} and a dilation of each variable maps $\Pi^d(r)$  to $\Pi$ and $\Pi^d(r)$ to $\Pi^d(1)$.  More generally, given an increasing collection $r_i$, the subdivisions $\Pi^d(r_i)$ are mutually transverse.  Moreover, there is a self diffeomorphism of $\bR^{d}$ mapping $\Pi^d(i)$ to $\Pi^d(r_i)$, which can be easily made to depend smoothly on the parameters $r_i$.  We conclude
\begin{lem} \label{lem:parametrise_family_intersections_cells}
If $\sigma^{i}$ is a sequence of cells of  $\Delta^d$, there exists a map
\begin{equation} \bigcap_{i} \check{\sigma}^{i}(i) \times \{ (r_1, \ldots, r_k) | r_1 \leq \cdots \leq r_k \} \to \bR^{d} \end{equation}
whose restriction to a fixed $(r_1, \ldots, r_k)$ is a diffeomorphism onto
\begin{equation*}  \bigcap_{i} \check{\sigma}^{i}(r_i). \end{equation*}  \noproof
\end{lem}  

To interpret the above results algebraically, we observe the intersections of all open cells in $\Pi^d$ and $\Pi^d(r)$ give a new polyhedral subdivision of $\bR^d$ which we write $\Pi^d \Cap \Pi^d(r)$.  Lemma  \ref{simplicial-model} implies that $\Pi^d \Cap \Pi^d (r)$ admits a canonical cellular map to $\Pi^d$ which collapses 
\[ \check{[v_i]} \cap \check{[v_j]}  (r)\]
to $\check{[v_i, v_j]}$.  It is easy to check that this map, which we've defined on the top simplices of $\Pi^d \Cap \Pi^d  (r)$, extends to a cellular map.  We write $\coll_{\Pi}$ for this map.

By composing the corresponding map on cellular chains with the intersection product
\[ C_{n-i}(\Pi^d) \otimes C_{n-j}(\Pi^d(r)) \to C_{n-i-j}( \Pi^d \Cap \Pi^d(r)) \] 
we obtain an intersection product valued in $ C_{n-i-j}(\Pi^d (r)) $.
\begin{lem} \label{model-compatible-cup-products}
The following diagram, where the vertical arrows are the canonical isomorphisms, is commutative
\[\xymatrix{C_{n-i}(\Pi^d) \otimes C_{n-j}(\Pi^d(r)) \ar[d] \ar[r] & C_{n-i-j}(\Pi^d (r))  \ar[d] \\
C^{i}(\Delta^d) \otimes C^{j}(\Delta^d) \ar[r]^{\cup} & C^{i+j}(\Delta^d) } \]
\end{lem}
\begin{proof}
This is essentially the content of Lemma \ref{simplicial-model}, once one observes that all intersections between cells of $\Pi^d$ and $\Pi^d(r)$ that are not recorded by the cup products on $C^*(\Delta)$ are collapsed by the map $\Pi^d \Cap \Pi^d (r)\to \Pi^d (r)$ to lower dimensional cells.
\end{proof}

This local model readily generalizes to our original setting with the simplex $\Delta(d)$ replaced by a triangulation $\sP$ of a manifold $Q$.

\begin{defin} \label{simplicial-position}
A pair of subdivisions $(\check{\sP}_{0}, \check{\sP}_{1})$ which are $\epsilon$ close to $\check{\sP}$ in the $C^1$ topology are said to be {\bf in simplicial position} if, whenever $C \in \check{\sP}$ is a cell of dimension $d$, there exists a coordinate chart near the $\epsilon$-essential part of $C$
\[(C - \partial C(\epsilon)) \times \bR^d \to Q \]
taking $(C - \partial C(\epsilon)) \times \Pi^d_0$ to $\check{\sP}$, and $(C - \partial C(\epsilon)) \times \Pi^d (r)$ to $\check{\sP}_1$, and which is compatible with the respective stratifications. 
\end{defin}

It is easy to check that constructing subdivisions in simplicial position is a completely local problem.  In fact, starting with a subdivision $\check{\sP}$, we use the ordering of the vertices of $\sP$ to order the top strata of $\check{\sP}$.  We may then construct a vector field which, on every cell of $\check{\sP}$ points in the direction of the nearby top cells which appears first in this order.  The image of $\check{\sP}$ under the flow of this vector field (for any sufficiently small time) is in simplicial position with respect to $\check{\sP}$.

 If $(\check{\sP}_{0}, \check{\sP}_{1})$  are in simplicial position, then the local copies of the polyhedral subdivisions $\Pi^d \Cap \Pi^d (r)$ patch together globally to a cellular decomposition $\check{\sP}_{0} \Cap \check{\sP}_{1}$ of $Q$.  The proof of the following global result is an immediate consequence of the local properties of $\Pi^d$ and $\Pi^d(r)$ recorded in Lemmas \ref{simplicial-model} and \ref{model-compatible-cup-products}:
\begin{prop} \label{properties-simplicial}
If $(\check{\sP}_{0}, \check{\sP}_{1})$ are in simplicial position, then cells $\check{\sigma}_0 \in \check{\sP}_{0}$ of and $\check{\tau}_1 \in \check{\sP}_{1}$ dual respectively to $\sigma =[u_0, \ldots, u_k] $ and $\tau = [w_0, \ldots, w_l]$ intersect if and only if 
\begin{itemize}
\item $u_k = w_0$ and the cell $\rho = [ u_0, \ldots, u_k = w_0, \ldots, w_l]$ exists, or
\item $u_i < w_j$ for all $i,j$ and the cell $\rho = [ u_0, \ldots, u_k, w_0, \ldots, w_l] $ exists.
\end{itemize}
Moreover, if $\check{\sP}_{0}$ and $\check{\sP}_{1}$  are $\delta$-close in the $C^1$ topology, then, in the first case, $\check{\sigma}_0 \cap \check{\tau}_1$ is $\delta$-close to $\check{\rho}_1$ in the $C^1$-topology (in particular, the cells are diffeomorphic), while $\check{\sigma}_0 \cap \check{\tau}_1$ lies in a $\delta$-neighbourhood of the lower dimensional cell $\check{\rho}_1$ in the second case.

Finally, collapsing all cells of  $\check{\sP}_{0} \Cap \check{\sP}_{1}$ of the second type defines a cellular map 
\[\coll_{\sP_1} \co \check{\sP}_{0} \Cap \check{\sP}_{1} \to \check{\sP}_1. \] \noproof
 \end{prop}

As in the local model, the intersection pairing therefore defines a map
\[ C_{n-i}(\check{\sP}_{0}) \otimes C_{n-j}(\check{\sP}_{1}) \to C_{n-i-j}(\check{\sP}_{0} \Cap \check{\sP}_{1}).\]
By composing the above map with the collapsing map $\check{\sP}_{0} \Cap \check{\sP}_{1} \to  \check{\sP}_{1}$, we obtain a product valued in  $C_{n-i-j}(\check{\sP}_{1})$.   We now have all the necessary pieces to introduce a cellular analogue of the category $\Rel(\sP)$.  

\begin{defin}
The category $\Cell(\check{\sP})$ has as objects pairs $(\check{\sP}_0, H_0)$ where $\check{\sP}_0$ is a cellular subdivision $C^1$-close to $\check{\sP}$, and $H_0$ is real valued function on $\partial Q$.

A pair of objects $(\check{\sP}_0, H_0)$ and $(\check{\sP}_1, H_1)$ is said to be transverse if $(\check{\sP}_0,\check{\sP}_1) $ are in simplicial position and $(H_0, H_1)$ is a transverse pair in $\Rel(\sP)$.  In this case, the space of morphisms between them is defined to be
\[  \Cell(\check{\sP})_*((\check{\sP}_0, H_0), (\check{\sP}_1, H_1)) \equiv C_{n-*}(\check{\sP}_1 - \partial^+_{H_1 - H_0} \check{\sP_1})  \]
where 
\[  \partial^+_{H_1 - H_0} \check{\sP}_1 = \{ C \in \check{\sP}_1 | C \cap \sigma \neq \emptyset, \textrm{ for some } \sigma \in \partial^+_{H_1 - H_0} \sP) \}.\]

A sequence of objects $\{ (\check{\sP}_i, H_i) \}_{i=0}^{d}$ is transverse whenever all pairs $(\check{\sP}_i,\check{\sP}_j)$ with $i < j$ are in simplicial position and the sequence $\{ H_i \}_{i=0}^d$ is transverse in $\Rel(\sP)$.  All higher products vanish, and the multiplication for a transverse triple $(\check{\sP}_0, H_0)$, $(\check{\sP}_1, H_1)$, and  $(\check{\sP}_2, H_2)$  
\[ C_{n-j}(\check{\sP}_{2} -  \partial^+_{H_2 - H_1} \check{\sP}_2)  \otimes C_{n-i}(\check{\sP}_{1} - \partial^+_{H_1 - H_0}\check{\sP}_1) \to C_{n-i-j}(\check{\sP}_{2} - \partial^+_{H_2 - H_0} \check{\sP}_2),\]
is defined by collapsing the intersection of cells in $ \check{\sP}_{1} $ and $ \check{\sP}_{2} $  to $ \check{\sP}_{2}$.

\end{defin}

\begin{prop} \label{cell-sing}
The forgetful map $(\check{\sP}_0, H_0) \to H_0$ extends to an $A_{\infty}$ functor
\[ \Cell(\check{\sP}) \to \Rel(\sP) \]
which is a quasi-equivalence.
\end{prop}
\begin{proof}
At the level of morphisms, the functor simply identifies a cellular chain with the corresponding dual cochain.  By construction, this is always an isomorphism of chain complexes.  It suffices therefore to check that the functor respects composition; this is the content of Lemma \ref{model-compatible-cup-products} locally, and can be readily seen to hold globally.
\end{proof}

\end{document}